\documentclass[USenglish]{article}

\usepackage[utf8]{inputenc}
\usepackage[small]{dgruyter}
\usepackage{microtype}

 \usepackage{latexsym}
 \usepackage{amsmath}
 \usepackage{amsfonts}
 \usepackage{amssymb}
 \usepackage{pdflscape}  
\usepackage{graphicx}
\usepackage[T1]{fontenc}
\usepackage[utf8]{inputenc}

\input{amssym.def}
\input{amssym.tex}

\def\theequation{\thesection.\arabic{equation}}
\newcommand\init{\setcounter{equation}{0}}

\newtheorem{theoreme}{Theorem }[section]
\newtheorem{proposition}[theoreme]{Proposition}
\newtheorem{lemma}[theoreme]{Lemma}
\newtheorem{definition}[theoreme]{Definition}
\newtheorem{corollary}[theoreme]{Corollary}
\newtheorem{remark}[theoreme]{Remark}

\newcommand{\beq}{\begin{equation}}
\newcommand{\eeq}{\end{equation}}
\newcommand{\bes}{\begin{subequations}}
\newcommand{\ees}{\end{subequations}}

\newcommand{\rrtimes}{{\rtimes}}
\def\bel{\begin{lemma}}
\def\eel{\end{lemma}}
\def\bet{\begin{theoreme}}
\def\eet{\end{theoreme}}
\def\bed{\begin{definition}}
\def\eed{\end{definition}}
\def\bep{\begin{proposition}}
\def\eep{\end{proposition}}
\def\ber{\begin{remark}}
\def\eer{\end{remark}}

\newcounter{smallarabics}
\newenvironment{arabicenumerate}
{\begin{list}{{\normalfont\textrm{(\arabic{smallarabics})}}}
  {\usecounter{smallarabics}\setlength{\itemindent}{0cm}
   \setlength{\leftmargin}{5ex}\setlength{\labelwidth}{4ex}
   \setlength{\topsep}{0.75\parsep}\setlength{\partopsep}{0ex}
   \setlength{\itemsep}{0ex}}}
{\end{list}}

\newcounter{smallroman}

\newcommand{\ben}{\begin{arabicenumerate}}
\newcommand{\een}{\end{arabicenumerate}}

\def\rr{{\mathbb R}}
\def\ss{{\mathbb S}}
\def\zz{{\mathbb Z}}
\def\cc{{\mathbb C}}
\def\nn{{\mathbb N}}

\newcommand{\bni}{\mathbin{\rotatebox[origin=c]{-90}{$\in$}}}
\newcommand{\leftmapsto}{\mathbin{\rotatebox[origin=c]{-180}{$\mapsto$}}}
\renewcommand\Re{\mathrm{Re}}
  \renewcommand\Im{\mathrm{Im}}
\newcommand\diff{{{\rm diff}}}
\newcommand\Diff{{{\rm Diff}}}
\def\I{{\rm\scriptscriptstyle I}}
\def\II{{\rm\scriptscriptstyle {I}{I}}}
\def\0{{\rm\scriptscriptstyle 0}}

\def\sph{{\rm sph}}

\def\heis{{\rm heis}}

\def\i{{\rm i}}

\def\fl{{\rm fl}}

\def\sch{{\rm sch}}
\def\Sch{{\rm Sch}}

\def\sgn{{\rm sgn}}
\def\id{{\iota}}

\def\e{{\rm e}}
\def\d{{\rm d}}

\def\cC{{\cal C}}

\def\cV{{\cal V}}
\def\cF{{\cal F}}

\def\cS{{\cal S}}

\def\cV{{\cal V}}
\def\cA{{\cal A}}

\def\cK{{\cal K}}

\def\cL{{\cal L}}

\def\cY{{\cal Y}}

\newcommand\tF{\tilde{F}}

\def\12{\frac{1}{2}}
\def\32{\frac{3}{2}}
\def\52{\frac{5}{2}}

\def\qed{$\Box$\medskip}
\def\proof{{\bf Proof.}\ \ }
\def\p{ \partial}

\def\x{\langle x \rangle}

\newcommand\gl{{\rm gl}}
\newcommand\agl{{\rm agl}}
\newcommand\GL{{\rm GL}}
\newcommand\AGL{{\rm AGL}}
\newcommand\aso{{\rm aso}}
\newcommand\ASO{{\rm ASO}}
\newcommand\AO{{\rm AO}}
\newcommand\so{{\rm so}}
\newcommand\SO{{\rm SO}}

\newcommand\dds{\,\partial_s}

\newcommand\ddr{\,\partial_r}
\newcommand\ddp{\,\partial_p}

\newcommand\ddw{\,\partial_w}

\newcommand\ddu[1]{\,\partial_{u_{#1}}}

\newenvironment{acknowledgments}{
  \addvspace{1.5\baselineskip}%
  \topsep=0pt\partopsep=0pt%
  \trivlist\item[\hspace{\labelsep}\bfseries\sffamily Acknowledgments.]
}{}

\begin{document}


  \author[1]{J. Derezi\'{n}ski}
  
  \affil[1]{Department of Mathematical Methods in Physics, 
Faculty of Physics,
 University of Warsaw, Pasteura 5, 02-093 Warszawa, Poland,
email jan.derezinski@fuw.edu.pl}
  
  \title{Group-theoretical origin
  of symmetries of hypergeometric class equations and functions}
  \runningtitle{Group-theoretical origin
  of symmetries of hypergeometric class equations and functions}
  \abstract{We show that  properties of hypergeometric class equations
and functions become transparent if  we   derive them from appropriate 2nd order differential equations with constant coefficients. More precisely, properties of the hypergeometric and Gegenbauer equation
can be  derived  from generalized symmetries of the Laplace equation in 4, resp. 3 dimension. Properties of the confluent, resp. Hermite equation can be derived from generalized symmetries of the heat equation in 2, resp. 1 dimension. Finally, the theory of the ${}_1F_1$ equation (equivalent to the Bessel equation) follows from the symmetries of the Helmholtz equation in 2 dimensions.
All  these symmetries become very simple when viewed on the level of the 6- or 5-dimensional ambient space.\\
Crucial role is played by the Lie algebra of generalized symmetries of these 2nd order PDE's, its Cartan algebra, the set of roots and the Weyl group.
Standard hypergeometric class functions
 are special solutions of these  PDE's  diagonalizing the Cartan algebra.
 Recurrence relations of these functions correspond to the roots. Their discrete symmetries correspond to the elements of the  Weyl group. }
  \keywords{hypergeometric equation, confluent equation, Hermite equation, Bessel equation, Lie groups, Lie algebras, conformal invariance, Laplace equation}
  \classification{2010 MSC: 33C80.}

\maketitle


\begin{center}\sf\large  Table of contents
\end{center}
\begin{enumerate}\sf
\item Introduction
  \item Hypergeometric class equations
  \item Pseudo-Euclidean spaces
    \item (Pseudo-)orthogonal group
\item Conformal invariance of the Laplacian
\item Laplacian in 4 dimensions and the hypergeometric equation
\item Laplacian in 3 dimensions and the Gegenbauer equation
\item The Schr\"odinger Lie algebra and the heat equation
\item Heat equation in 2 dimensions and the confluent equation
\item Heat equation in 1 dimension and the Hermite equation
\item The Helmholtz equation in 2 dimensions and the ${}_0F_1$  equation
  \end{enumerate}

\section{Introduction}\label{s1}\init

These lecture notes  
are
devoted to the properties of
the following equations:

\bigskip

\noindent the {\em Gauss hypergeometric equation}, called also the {\em ${}_2F_1$ equation},
\beq \Big(w(1-w)\p_w^2+\big(c-(a+b+1)w\big)\p_w-ab\Big)F(w)=0;
\label{hyp1}
\eeq
the {\em Gegenbauer equation}
\beq \Big((1-w^2)\p_w^2-({a}+{b}+1)w\p_w-{a}{b}\Big)F(w)=0;
\label{hyp2}\eeq
 {\em Kummer's confluent equation}, called also the {\em${}_1F_1$ equation}, 
\beq\Big(w\p_w^2+(c-w)\p_w-a\Big)F(w)=0;
\label{hyp3}\eeq
the {\em Hermite equation}
\beq \Big(\p_w^2-2w\p_w-2{a}\Big)F(w)=0;
\label{hyp4}\eeq
and the {\em ${}_0F_1$ equation} (equivalent  to the better known {\em Bessel
  equation}, see eg. \cite{De})
\beq \Big(w\p_w^2+c\p_w-1\Big)F(w)=0.\label{hyp5}\eeq
Here, $w$ is a complex variable,
 $\partial_w$ is the differentiation with respect to $w$,
 and $a,b,c$ are arbitrary complex parameters.

 These equations are typical representatives of the so-called
 {\em hypergeometric class equations} \cite{NU}.
 (Nikiforov and Uvarov call them {\em hypergeometric type equations}; following
 \cite{SL}, we prefer in this context to use the word {\em class}, reserving {\em type} for narrower families of equations).
 We refer the reader to Sect. \ref{s2}, where  we discuss the terminology concerning
 hypergeometric class equations and functions that we use.

The equations (\ref{hyp1})--(\ref{hyp5})
and their solutions belong to the most natural objects of mathematics and often appear in applications \cite{Fl,MF,WW}.

The aim of these notes is to elucidate the mathematical   structure of a large  class of  identities satisfied by hypergeometric class equations and functions.
We believe that our approach brings order and transparency to this subject, usually considered to be complicated and messy.

We will restrict ourselves to {\em generic
  parameters} $a,b,c$. We will not discuss special properties of  two distingushed classes of parameters, when additional identities are true:
\ben\item the {\em polynomial case} (which corresponds to negative integer values of $a$);
\item the {\em degenerate case} (which 
   corresponds to integer values of $c$).
  \een

The notes are to a large extent  based on  \cite{De} and \cite{DeMaj}, with some additions and improvements.

\subsection{From 2nd order PDE's with constant coefficients to hypergeometric class equations}\label{sub1.1}

In our approach, each of the equations  (\ref{hyp1})--(\ref{hyp5}) is
derived from a certain
{\em complex  2nd order 
  PDE  with constant coefficients}.
 The identities satisfied by this PDE and their solutions
 are very straightforward---they look obvious and symmetric.
 After an  appropriate change of variables,
we derive (\ref{hyp1})--(\ref{hyp5}) and  identities satisfied by their solutions. They look  much more  complicated and messy.

We will argue that the main source of these identities are {\em generalized symmetries} of the parent  PDE. Let us briefly recall this concept.

Suppose that we are given an equation
\beq \cK f=0,\label{repo0}\eeq
where $\cK$ is a linear differential operator. 
Let $g$ be a Lie algebra and  $G$ a group
equipped with pairs of representations
\bes\begin{align}
  g&\ni B\mapsto B^\flat,B^\#,\label{repo2}\\
    G&\ni \alpha\mapsto \alpha^\flat,\alpha^\#,\label{repo1}
\end{align}\ees
where  (\ref{repo2}) has its values in   1st order differential operators and
(\ref{repo1}) in  point transformations with multipliers.
We say that (\ref{repo2}) and (\ref{repo1}) are generalized symmetries of (\ref{repo0}) if
\bes\begin{align}
 B^\flat \cK&=\cK  B^\#,\label{trans2}\\
  \text{resp.}\qquad     \alpha^\flat \cK&=\cK  \alpha^\#.\label{trans1}
\end{align}\ees
Note that (\ref{trans2}), resp. (\ref{trans1}) imply that $B^\#$ and $\alpha^\#$ preserve the space of solutions of (\ref{repo0}).

We will omit the word ``generalized'' if $B^\#=B^\flat$ and $\alpha^\flat=
\alpha^\#$.

 We can distinguish 3 kinds of PDE's with constant coefficients in complex domain. Below we list these PDE's, together with the   Lie algebra and group of their generalized symmetries:
\ben
\item The {\em Laplace equation} on $\cc^n$
  \beq \Delta_n f=0,\quad n>2.\label{lapla1}\eeq
  The {\em orthogonal Lie algebra} and {\em group} in $n{+}2$ dimensions, denoted $\so(n{+}2,\cc)$, resp. $\mathrm{O}(n{+}2,\cc)$, both  acting conformally in $n$ dimensions. (For $n=1,2$ there are additional conformal symmetries).
 \item The  {\em heat equation} on $\cc^{n-2}\oplus\cc$:
   \beq(\Delta_{n{-}2}+2\partial_s)f=0.\label{lapla2}\eeq
 The 
   {\em Schr\"odinger Lie algebra} and {\em group}
in  $n{-}2$ dimensions, denoted
  $\sch(n{-}2,\cc)$, resp.  $\Sch(n{-}2,\cc)$.
\item The  {\em Helmholtz equation} on $\cc^{n-1}$,
  \beq (\Delta_{n{-}1}-1)f=0.\label{lapla3}\eeq
  The {\em affine orthogonal  Lie algebra} and {\em group}  in $n{-}1$ dimensions, denoted  $\aso(n{-}1,\cc)$, resp. $\AO(n{-}1,\cc)$.
     \een
(The reason for the strange choice of dimensions in (\ref{lapla2}) and (\ref{lapla3}) will be explained later).

     The basic idea of our approach is as follows. Let us start from the  equation (\ref{repo0}), where $\cK$ is appropriately chosen from among
(\ref{lapla1}), (\ref{lapla2}) and (\ref{lapla3}).     
     In the Lie algebra
     of its generalized symmetries  we   fix
a certain maximal commutative  algebra, which we will call the ``Cartan algebra''. Operators that are eigenvectors of the  adjoint action of the ``Cartan algebra'' will be called ``root operators''.

In the group of generalized symmetries we fix a subgroup, which we call the ``Weyl group''. It is chosen in such a way, that its adjoint action fixes the 
 ``Cartan algebra''.

Note that in some cases the Lie algebra of symmetries is simple, and then the names {\em Cartan algebra}, {\em root operators} amd {\em Weyl symmetries} correspond to the standard names. In other cases the Lie algebra is not semisimple, and then the names are less standard -- this is the reason for the quotation marks that we use above. In the sequel we drop the quotation marks.

Let us fix a basis of the Cartan algebra $N_1,\dots,N_k$.
Suppose that the dimension of the underlying space is by $1$ greater than the dimension of the Cartan algebra. Then we  introduce new variables, say
$w,u_1,\dots,u_k$ such that $N_i=u_i\partial_{u_i}$.

Substituting a function of the form
\beq f=u_1^{\alpha_1}\cdots u_k^{\alpha_k}F(w),\label{subo1}\eeq
to the equation (\ref{repo0}), and using
\beq
N_i u^{\alpha_i}=\alpha_iu^{\alpha_i}\eeq
we obtain
the equation
\beq
\cF_{\alpha_1,\dots,\alpha_k}F=0,\label{subo2}\eeq
 which coincides with one of the equations (\ref{hyp1})--(\ref{hyp5}). The eigenvalues of the Cartan operators become the parameters of this equation.

Root operators shift the Cartan elements, typically by $1$ or $-1$ (like the well-known creation and annihilation operators). Therefore, root operators inserted into the relations (\ref{trans2}) lead to {\em transmutation relations} for
(\ref{hyp1})--(\ref{hyp5}).

Similarly, elements of the Weyl  group
 permute  Cartan elements or change their signs.
 Therefore, Weyl symmetries inserted into 
 (\ref{trans1}) leads to  {\em discrete symmetries} of
(\ref{hyp1})--(\ref{hyp5}).

Of course, one can apply (\ref{trans1}) to elements of $G$ other than Weyl symmetries, obtaining interesting integral and addition identities for hypergeometric class functions.  They  are, however, outside of the scope of these notes.

There are  five 2nd order PDE with constant coefficients where we can  perform this procedure. They are all listed in the following table:

\[\begin{array}{ccccc}
\hbox{PDE}&\begin{array}{c}\hbox{Lie}\\ \hbox{algebra}\end{array}
&\begin{array}{c}\hbox{dimension of}\\ \hbox{Cartan algebra}\end{array}
&\begin{array}{c}\hbox{discrete}\\ \hbox{symmetries}\end{array}&
\hbox{equation}\\[1ex]
\hline\\[1ex]
\Delta_4&\so(6,\cc)&3&\hbox{cube}& {}_2F_1; \\[1.5ex]
\Delta_3&\so(5,\cc)&2&\hbox{square}&\hbox{Gegenbauer};\\[1.5ex]
\Delta_2+2\p_t&\sch(2,\cc)&2&\zz_2\times
\zz_2&
{}_1F_1\hbox{ or }{}_2F_0;\\[1.5ex]
\Delta_1+2\p_t&\sch(1,\cc)&1&\zz_4&
\hbox{Hermite};\\[1.5ex]
\Delta_2-1&\aso(2,\cc)&1&\zz_2&{}_0F_1.
\end{array}\]

Note that some other 2nd order PDE's have too few  variables to be in the above list: this is the case of $\Delta_1$ and $\Delta_2$.
Others have too many variables: one can try to perform the above procedure,
however it leads to a differential  equation in more than one variable.

\subsection{Conformal invariance of the Laplace equation}

The key tool of our approach is the conformal invariance of
the Laplace equation. Let us sketch a derivation of this invariance. For simplicity we restrict our attention to the complex case, for which
we do not need to distinguish between various signatures  of the metric tensor.

In order to derive the conformal invariance of the Laplacian on $\cc^n$, or on other  complex manifolds with maximal conformal symmetry, it is convenient  to start  from the so-called {\em ambient space} $\cc^{n+2}$, where the actions of  $\so(n{+}2,\cc)$  and $\mathrm{O}(n{+}2,\cc)$ are obvious. In the next step these actions are restricted to the {\em null quadric}, and finally to the {\em projective null quadric}. Thus the dimension of the manifold goes down from $n{+}2$ to $n$. The null quadric can be viewed as a line bundle over the projective null quadric. By choosing an appropriate {\em section} we can identify the projective null quadric, or at least its open dense subset, with the flat space $\cc^n$ or some other complex manifolds with a complex Riemannian structure, e.g. the product of two spheres. 
The Lie algebra $\so(n{+}2,\cc)$  and the group $\mathrm{O}(n{+}2,\cc)$ act conformally on these manifolds.

What is more interesting, the above construction
leads to a definition of an invariantly defined operator, which we denote $\Delta^\diamond$, transforming functions on the null quadric homogeneous of degree $1-\frac{n}{2}$ onto functions homogeneous of degree $-1-\frac{n}{2}$. After fixing a section, this operator can be identified with the conformal Laplacian on the corresponding complex Riemannian manifold of dimension $n$. For instance, one obtains the Laplacian $\Delta_n$ on $\cc^n$. 
The  representations of   $\so(n{+}2,\cc)$ and $\mathrm{O}(n{+}2,\cc)$ on the level of the ambient space were true symmetries of $\Delta_{n+2}$.  After the reduction to $n$ dimensions,  they become generalized symmetries of the conformal Laplacian.

The fact that conformal transformations of the Euclidean
space are generalized symmetries of the Laplace equation was
apparently known already to Lord Kelvin.
Its explanation in terms of the null quadric first appeared
in \cite{Boc}, and is discussed e.g. in \cite{CGT}.
The reduction of $\Delta_{n+2}$ to $\Delta_n$ mentioned above, is based on a beautiful idea of
 Dirac in \cite{Dir}, which was later rediscovered e.g. in \cite{HH,FG}---see a discussion by  Eastwood \cite{East}.

The construction indicated above gives a rather special class of \hbox{(pseudo-)}Rie\-mannian manifolds---those having a conformal group of maximal dimension, see e.g. \cite{EMN}.
However, conformal invariance can be generalized to 
arbitrary \hbox{(pseudo-)}Rie\-mannian manifolds. In fact, the
Laplace-Beltrami operator plus an appropriate multiple of the scalar
curvature, sometimes called the {\em Yamabe Laplacian}, is invariant in
a generalized sense with respect to
conformal maps, see e.g.
\cite{Tay,Or}.


\subsection{The Schr\"odinger  Lie algebra and Lie algebra as generalized symmetries of the Heat equation}
\label{subsect-heat}

The  heat equation 
(\ref{lapla2}) possesses a large Lie algebra and group of generalized symmetries, which in the complex case, as we already indicated, we denote by
 $\sch(n{-}2,\cc)$ and $\Sch(n{-}2,\cc)$. Apparently, they were known already to
 Lie \cite{L}. They were rediscovered (in the
  essentially equivalent context of
 the free
 Schr\"odinger equation) by Schr\"odinger \cite{Sch}. They were  then studied e.g. in \cite{Ha,Ni}.

By adding an additional variable, one can consider the heat equation as the Laplace equation acting on functions with an exponential dependence on one of the variables. This allows us to express  generalized symmetries
of (\ref{lapla2}) by generalized symmetries of (\ref{lapla1}). They can be identified as a subalgebra  of $\so(n{+}2,\cc)$, resp. a subgroup of
$\mathrm{O}(n{+}2,\cc)$ consisting of elements commuting with a certain distinguished element of $\so(n{+}2,\cc)$.

\subsection{Affine orthogonal group and algebra as symmetries of the Helmholtz equation}
\label{subsect-helm}

Recall that the affine orthogonal group $\AO(n{-}1,\cc)$ is generated by rotations and translations of $\cc^{n-1}$. 
It is obvious  that elements of $\AO(n{-}1,\cc)$
commute with the Helmholtz operator $\Delta_{n-1}-1$. The same is true concerning the affine orthogonal Lie algebra $\aso(n{-}1,\cc)$. Therefore, they are symmetries of the Helmholtz equation  (\ref{lapla3}).

The Helmholtz equation is conceptually simpler than that of the Laplace and heat equation,
because all generalized symmetries are true symmetries.

Note that one can embed the symmetries of the Helmholtz equation in conformal symmetries of the Laplace equation, similarly as was done with the heat equation. In fact,
$\aso(n{-}1,\cc)$ is a subalgebra of $\so(n{+}2,\cc)$, and  $\AO(n{-}1,\cc)$ is a subgroup of $\mathrm{O}(n{+}2,\cc)$.

\subsection{Factorization relations}
\label{Factorization relations}

Another important class of identities satisfied by
hypergeometric class operators are
{\em factorizations} \cite{IH}. They come in pairs. They are identities of the form
\bes\begin{align}
   \cF_1&=\cA_-\cA_++c_1,\label{facto1-}\\
   \cF_2&=\cA_+\cA_-+c_2,\label{facto-}
   \end{align}\ees
where $\cA_+$, $\cA_-$ are 1st order differential operators, $c_1$, $c_2$ are numbers and $\cF_1$, $\cF_2$ are operators coming from
 (\ref{hyp1})---(\ref{hyp5}) with slightly shifted parameters.

 The number of  such factorizations is the same as the number of roots of the Lie algebra of generalized symmetries. They can be derived from certain identities in the {\em enveloping algebra}. They are closely related to the {\em Casimir operators} of its subalgebras.

Factorizations imply  {\em transmutation relations}.
In fact, it is easy to see that
(\ref{facto-}) and (\ref{facto1-}) imply
\bes\begin{align}
  \cA_-\cF_2&=(\cF_1+c_2-c_1)\cA_-,\label{darb1}\\
  \cA_+\cF_1&=(\cF_2+c_1-c_2)\cA_+.\label{darb2}
  \end{align}
\ees

Note that (\ref{darb1}) implies that
the operator $\cA_-$ maps the kernel of $\cF_2$ to the kernel of
$\cF_1+c_2-c_1$. Similarly, (\ref{darb2}) implies that the operator
$\cA_+$ maps the kernel of $\cF_1$ to the kernel of
$\cF_2+c_1-c_2$. The above construction is usually called
 the  {\em Darboux transformation}.

\subsection{Standard solutions of hypergeometric class equations}

So far we discussed only identities satisfied by the operators corresponding to the equations
(\ref{hyp1})---(\ref{hyp5}). The approach discussed in these notes is also helpful in deriving and classifying the identities for their solutions.

The equations (\ref{hyp1})---(\ref{hyp5}) have at least 1 and at most 3 singular points on the Riemann sphere. One can typically find two  solutions with a simple behavior at each of these points. We call them {\em standard solutions}. (If it is a regular--singular point, then the solutions are given by convergent power series, otherwise we have to use other methods to define them). The discrete symmetries  map standard solutions on standard solutions. The best known example of this method of generating solutions is {\em Kummer's table} \cite{Ku}, which lists various possible expressions for solutions of the hypergeometric equation.

\subsection{Recurrence relations of hypergeometric class functions}

All transmutation relations have the form
\beq
\cA\cF_1=\cF_2\cA,\label{recu1}\eeq
where $\cA$ is a first order differential operator and $\cF_1$, $\cF_2$ is a pair of hypergeometric class operators of the same type. Typically, some parameters of $\cF_2$ differ  from the corresponding parameters of $\cF_1$ by $\pm1$. Clearly, if a function $F_1$ solves $\cF_1F_1=0$, then $\cA F_1$ solves
$\cF_2\cA F_1=0$.

It turns out  that if $F_1$ is a standard solution of $\cF_1$, then $\cA F_1$ is proportional to one of standard solutions of $\cF_2$, say $F_2$. Thus we obtain
an identity
\beq
\cA F_1=a F_2, \label{recu}\eeq
called a {\em recurrence relation}, or a {\em contiguity relation}.

The recurrence relation (\ref{recu}) is fixed by the transmutation relation (\ref{recu1}) except for the coefficient $a$. In practice it is not difficult to determine $a$.

\subsection{From wave packets to integral representations}

Hypergeometric class functions possess integral representations, where integrands are elementary functions. We show that  integral representations come from  certain natural solutions of the parent 2nd order PDE, which at the same time are eigenfunctions of Cartan operators. It will be convenient to have a name for this kind of solutions---we will  call them {\em wave packets}.

Let us describe how to construct wave packets for the Laplace equation.
It is easy to see that each function depending only on variables from an isotropic subspace  is harmonic, that is, satisfies the Laplace equation. By assuming that the function is homogeneous in appropriate variables we can make sure that it is an eigenfunction of Cartan operators.

Unfortunately, the above class of functions is too narrow for our purposes. There is still another construction that can be applied: we can rotate a function and integrate it (``smear it out'')
with respect to a weight. This procedure does not destroy the harmonicity. By choosing the weight appropriately, we can make sure that the resulting wave packet is an eigenfunction of Cartan perators.
(The ``smearing out'' is essentially a generalization of the Fourier (or Mellin) transformation to the complex domain.)

 After substituting special coordinates to a wave packet, we obtain a
 function of the form (\ref{subo1}) with $F$ solving (\ref{subo2}), and 
having the form of an integral of an elementary function.

Wave packets for the heat and Helmholtz equation can be derived from wave packets for the Laplace equation.

\subsection{Plan of the lecture notes}

In Sect.  \ref{s2} we give a concise introduction to hypergeometric class equations and functions. One can view this section as an extension of the introduction, concentrated on the terminology and classification of equations and functions we consider in these notes.

  The remaining sections  can be divided into two categories.
  The first category consists of
  Sects \ref{s3}, \ref{s5} and  \ref{s8}. They have a general character and are devoted to basic geometric analysis in any dimension.
  The most important one among them is Sect. \ref{s5}, devoted to the conformal invariance of the Laplace equation.
  Of comparable importance is  Sect. \ref{s8},
 where
 the Schr\"odinger Lie algebra and  group  are introduced.
In Subsect. \ref{Harmonic functions}---\ref{wave3} we explain how to construct ``wave packets''. 
No special functions appear in Sects \ref{s3}, \ref{s5} and  \ref{s8}. They can be read independently of the rest of the notes.

The second category consists  of
Sects \ref{s6}, \ref{s7}, \ref{s9}, \ref{s10} and \ref{s11}. They are  devoted to a detailed analysis of equations (\ref{hyp1}), (\ref{hyp2}),
(\ref{hyp3}), (\ref{hyp4}), resp. (\ref{hyp5}). Typically,
each section starts with
the ambient space corresponding to
the 2nd order PDE from the left column of the table in Subsect. \ref{sub1.1}.
In the ambient space these symmetries are very easy to describe. 
Then we reduce the dimension and introduce special coordinates, which leads to the equation in the right
column of the table.

We made serious efforts to make Sects \ref{s6}, \ref{s7}, \ref{s9}, \ref{s10} and \ref{s11} as parallel as possible.
there is a one to one correspondence between subsections in all these 5 sections.
 We try to use a uniform terminology and analogous conventions. This makes our text somewhat repetitive---we believe that this is  helpful to the reader. Note also that these sections are to a large extent independent of one another.

We use various (minor but helpful) ideas to make our presentation as short and transparent
as possible. One of them is the  use of two kinds of parameters. The parameters that appear in  (\ref{hyp1}), (\ref{hyp2}),
(\ref{hyp3}), (\ref{hyp4}), and (\ref{hyp5}),  denoted $a,b,c$, are called {\em classical parameters}. They are convenient when one defines  ${}_kF_m$ functions by power series. However, in most of our text we prefer to use a different set of parameters,  denoted by Greek letters  $\alpha,\beta,\mu,\theta,\lambda$. They are much more convenient when we  describe  symmetries.

Another helpful idea is a consistent use of {\em split coordinates}
in $\cc^n$ or $\rr^n$.
In these coordinates root operators and Weyl symmetries have an especially simple form.

The notes are full of long lists of identities. We are convinced that most of them are easy to   understand and appreciate  without much effort. Typically, they are highly symmetric and parallel to one another.

We hesitated whether to use the complex or real setting for these notes.
The complex setting was e.g. in \cite{DeMaj}.
It offers undoubtedly some simplifications: there is no need to consider various signatures of the scalar product.
However, the complex setting can also be problematic: analytic functions are often  multivalued, which causes issues with some global constructions.
Therefore, in these notes, except for the introduction, we use  the real setting as the basic one. At the same time we keep in mind that all our formulas have  obvious analytic continuations to  appropriate complex domains.

In most of our notes, we do not make explicit the signature of the scalar product in our notation for Lie algebras and groups.
E.g. by writing  $\so(n)$ we mean $\so(q,p)$ for some $n=q+p$ or  $so(n,\cc)$.
Specifying each time the signature would be overly pedantic, especially since we usually want to complexify all objects, so that the signature loses its importance.


\subsection{Comparison with literature}

The literature about hypergeometric class functions is enormous---after all it is one of the oldest subjects of mathematics. Let us mention e.g. the books \cite{BE,SL,AAR,EMOT,Ho,MOS,NIST,R,WW}.

The relationship of special functions  to Lie groups and algebras was noticed  long time ago. For instance,
the papers by Weisner \cite{We1,We2} from the 50's describe
Lie algebras
associated with  Bessel and Hermite functions.

The idea of studying hypergeometric class equations with help of Lie
algebras
 was developed further by Miller. His early book \cite{M1}
considers mostly small Lie algebras/Lie groups, typically
$\mathrm{sl}(2,\cc)$/$\mathrm{SL}(2,\cc)$ and their contractions, and applies
them to obtain various identities about hypergeometric class functions. These Lie
algebras have
1-dimensional
Cartan algebras and a single pair of roots. This kind of analysis
is able to explain only a single pair of transmutation relations for each equation.
To explain bigger families of transmutation relations one
needs larger Lie algebras.

A Lie algebra strictly larger than $\mathrm{sl}(2,\cc)$ is $\so(4,\cc)$.
There exists a large literature on the relation of the
 hypergeometric equation with $\so(4,\cc)$
and its
  real forms,
see eg. \cite{KM,KMR}. This Lie algebra is however still too small to
account for all symmetries of the hypergeometric equation---its Cartan
algebra is only 2-dimensional, whereas the equation has three parameters.

An explanation of symmetries of the
Gegenbauer equation in terms of $\so(5,\cc)$ and of the hypergeometric
equation in terms of $\so(6,\cc)\simeq \mathrm{sl}(4,\cc)$ was first given by Miller,
see \cite{M4}, and especially \cite{M5}.

Miller and Kalnins wrote a series of papers where they studied the
symmetry approach to
separation of variables for various 2nd order partial differential
equations, such as the Laplace and wave equation, see eg. \cite{KM1}.
A large part of this research
is summed up in the book by Miller \cite{M3}. As an important
consequence of this study, one obtains detailed information about
symmetries of hypergeometric class equations.

The main tool that we use
to describe properties of hypergeometric class functions
are {\em generalized
 symmetries} of 2nd order linear PDE's. Their theory is described in another
book by Miller \cite{M2}, and further developed in \cite{M3}.

A topic that is extensively treated in the literature on the relation
of special functions to group theory, such as \cite{V,Wa,M1, VK}, is
derivation of various
addition formulas. Addition formulas
say that a certain special function can be
written as a sum, often infinite, of some related functions.
As we mentioned above, they are outside of the scope of this text---we concentrate on the simplest identities.

The relationship of Kummer's table with the group of symmetries of a
cube (which is the Weyl group of $\so(6,\cc)$) was discussed in \cite{LSV}.
A recent paper, where symmetries of the hypergeometric equation play
an important role is \cite{Ko}.

The use of transmutation relations as a tool
  to derive recurrence relations
for hypergeometric class functions is well known and can be found
eg. in the book by Nikiforov-Uvarov \cite{NU}, in the books by Miller
\cite{M1} or in older works such
as \cite{Tr,We1,We2}.

There exist various generalizations of hypergeometric class
 functions. Let us mention the
class of $\mathcal{A}$-hypergeometric functions, which provides a
natural generalization of the usual hypergeometric
function to many-variable situations \cite{Be,Bod}.
Saito \cite{Sa} considers
 generalized symmetries in the framework of $\mathcal{A}$-hypergeometric functions.

Another direction of generalizations of hypergeometric functions is the family of Gel’fand-Kapranov-Zelevinsky hypergeometric functions \cite{G,GKZ}. Similar constructions were explored by  Aomoto and others \cite{A,AK,M-H}. The main idea is to generalize  integral representations of hypergeometric functions, rather than hypergeometric equations. There exist  also interesting confluent versions of these functions \cite{KHT}.

A systematic presentation and derivation of  symmetries  of hypergeometric class equations and functions  from
2nd order PDE's with constant coefficients was given in  \cite{De} and \cite{DeMaj}.
These papers consistently use Lie-algebraic parameters, describe transmutation relations, discrete symmetries and factorizations. \cite{De} describes integral representations and recurrence relations.
\cite{DeMaj} concentrates on the study of hypergeometric class operators, leaving out the properties of hypergeometric class functions.

These lecture notes are to a large extent based on \cite{De} and \cite{DeMaj}.
There are some corrections  and minor changes of conventions. There are also some   additions.  A systematic derivation
of all integral representations from  ``wave packets'' in higher dimensions
seems to be new.

There are a number of topics related to the hypergeometric class equation that we do not touch. Let us mention the question
whether hypergeometric functions can be expressed in terms of
algebraic functions. This topic, in the context of
 $\mathcal{A}$-hypergeometric functions was considered eg. in the interesting
 papers \cite{Be,Bod}.

 We  stick to a rather limited class of equations and functions
(\ref{hyp1})--(\ref{hyp5}).  They have a surprisingly rich structure, which  often seems to be lost
 in more general classes. Nevertheless, it is natural to ask how far one can generalize  the ideas of these notes to other  equations and functions, such as higher hypergeometric functions, multivariable hypergeometric functions, Heun functions, $q$-hypergeometric functions, Painlev\'{e}{} equations.

\begin{acknowledgments}
  The support of the National Science Center under the grant UMO-2014/15/B/ST1/00126 is gratefully acknowledged.
  The author thanks P.~Majewski for collaboration at \cite{DeMaj}.
  He is also grateful to A.~Latosi\'{n}ski, T.~Koornwinder, M.~Eastwood, S.-Y. Matsubara-Heo and Y. Haraoka for useful remarks.
\end{acknowledgments}

\section{Hypergeometric class equations}
\label{s2}
\init
In this short section we fix our terminology concerning
hypergeometric class  equations and functions. 

\subsection{Remarks on notation}

We use  $\partial_w$ for the operator of differentiation in the variable $w$.
We will understand that the operator $\p_w$ acts on the whole expression on its right:
\beq
\partial_w f(w)g(w)= \partial_w \big(f(w)g(w)\big).\label{req0}\eeq
If we want to restrict the action of $\partial_w$ to the term immediately to the right, we will write $ f(w)_{,w}$, or simply  $f'(w)$.

We use lhs and rhs as the abbreviations for the left hand side and right hand side.

\subsection{Generalized hypergeometric series}

For $a\in\cc$ and $n\in\nn$ we define the {\em Pochhammer symbol}
\[(a)_j:=a(a+1)\cdots(a+j-1).\]

For $a_1,\dots,a_k\in\cc$, $c_1,\dots,c_m\in\cc\backslash\{0,-1,-2,\dots\}$,
we define the {\em  ${}_kF_m$
 generalized hypergeometric series}, or for brevity the  {\em  ${}_kF_m$ series}:
 \beq{}_kF_m(a_1,\dots, a_k;c_1,\dots,c_m;w):=
\sum_{j=0}^\infty\frac{(a_1)_j\cdots(a_k)_jw^j}{(c_1)_j\cdots(c_m)_jj!}.
\label{hiper}\eeq
By the d'Alembert criterion,
\ben\item  if $m+1>k$,   the series  (\ref{hiper}) is convergent for $w\in\cc$;
\item  if $m+1=k$,  the series  (\ref{hiper}) is convergent for $|w|<1$;
\item  if $m+1<k$,  the series  (\ref{hiper}) is divergent, however sometimes
a certain function can be naturally associated with  (\ref{hiper}). 
\een
The corresponding analytic function will be called the  {\em ${}_kF_m$ function}.


The zeroth order term of the series (\ref{hiper})  is $1$. A different normalization of (\ref{hiper}) is often useful:
\begin{align}
{}_k{\bf F}_m(a_1,\dots, a_k;c_1,\dots,c_m;w)&:=
\frac{{}_kF_m(a_1,\dots, a_k;c_1,\dots,c_m;w)}{\Gamma(c_1)\cdots\Gamma(c_m)}
\nonumber\\
&=\sum\limits_{j=0}^\infty
\frac{(a_1)_j\cdots(a_k)_jw^j}{\Gamma(c_1+j)\cdots\Gamma(c_m+j)j!}.
\label{hiper1}\end{align}
In (\ref{hiper1}) we do not have to restrict the values of
$c_1,\dots,c_m\in\cc$.

\subsection{Generalized hypergeometric equations}

\bet The  ${}_kF_m$ function (\ref{hiper}) solves the dfferential equation
\begin{equation}
  \begin{array}{rl}&
(c_1+w\partial_w)\cdots(c_m+w\partial_w)\partial_wF(a_1,\dots, a_k;c_1,\dots,c_m;w)\\[1ex]
=
&(a_1+w\partial_w)\cdots(a_k+w\partial_w)F(a_1,\dots, a_k;c_1,\dots,c_m;w).
\end{array}
\label{equa}\end{equation}
\eet

\proof
We check that both  the left and right hand side of (\ref{equa}) are equal to
\[a_1\cdots a_kF(a_1+1,\dots, a_k+1;c_1,\dots,c_m;w).\]
\qed

 We will call (\ref{equa})  the {\em ${}_k\cF_m$ equation}. It has the order   $\max(k,m+1)$.
Below we list all  ${}_kF_m$ functions with equations of the order at most $2$.

\begin{itemize}
\item {\bf  The ${}_2F_1$ function or the Gauss hypergeometric  function}
 \[\begin{array}{l}
 F(a,b;c;w)=\sum\limits_{n=0}^\infty
 \frac{(a)_n(b)_n}{n!(c)_n}w^n.
\end{array}\]
The series is convergent for $|w|<1$, and it extends to a multivalued
function on a covering of $\cc\backslash\{0,1\}$.
It is a solution of
 the {\em Gauss hypergeometric equation} or the {\em ${}_2\cF_1$ equation}
\label{hipp}
\[\left(w(1-w)\p_w^2+(c-(a+b+1)w)\p_w-ab\right)f(w)=0.\]

\item {\bf The ${}_1F_1$ function or Kummer's confluent function}

 \[\begin{array}{l}
 F(a;c;w)=\sum\limits_{n=0}^\infty
 \frac{(a)_n}{n!(c)_n}w^n.\end{array}\]
The series is convergent for all $w\in\cc$.
It is a solution of {\em Kummer's
confluent equation} or the {\em ${}_1\cF_1$ equation}
 \[\big(w\p_w^2+(c-w)\p_w-a\big)f(w)=0.\]

\item {\bf The ${}_0F_1$ function}

 \[\begin{array}{l}
F(-;c;w)= F(c;w)=\sum\limits_{n=0}^\infty
 \frac{1}{n!(c)_n}w^n
.\end{array}\]

The series is convergent for all $w\in\cc$.
It is a solution of  the
{\em ${}_0\cF_1$ equation} (related to the Bessel equation)
 \[(w\p_w^2+c\p_w-1)f(w)=0.\]

\item {\bf The ${}_2F_0$ function}

 For $\arg w\neq0$ we define
\[F(a,b;-;w):=\lim_{c\to\infty}F(a,b;c;cw).\]
It extends to an analytic function on the universal cover of
$\cc\backslash\{0\}$ 
with a branch point of an infinite order at 0.
It has the following divergent but asymptotic expansion:
\[
F(a,b;-;w)\sim\sum_{n=0}^\infty\frac{(a)_n(b)_n}{n!}w^n,
\ |\arg w-\pi|<\pi-\epsilon,\quad \epsilon>0.
\]
It is a solution of the ${}_2\cF_0$ equation
\[\left(w^2\p_w^2+(-1+(a+b+1)w)\p_w+ab\right)f(w)=0.\]
By a simple transformation described in Subsect. \ref{The ${}_2F_0$ function}
it is equivalent to the ${}_1\cF_1$ equation.

\item {\bf The  ${}_1F_0$  function or the power function}

\[\begin{array}{l}F(a;-;w)=(1-w)^{-a}=
\sum\limits_{n=0}^\infty
 \frac{(a)_n}{n!}w^n
\end{array}\]

It solves
\[\big((w-1)\p_w-a\big)f(w)=0.\]

\item {\bf The ${}_0F_0$ function or the  exponential function}

\[\begin{array}{l}F(-;-;w)=\e^w=
\sum\limits_{n=0}^\infty
 \frac{1}{n!}w^n.
\end{array}\]
It solves
\[(\p_w-1)f(w)=0.\]
\end{itemize}

\subsection{Hypergeometric class equations}

\label{sub-req}


Following \cite{NU}, equations of the form
\beq\left(\sigma(w)\p_w^2+\tau(w)\p_w+
\eta\right) f(w)=0,\ \ \
\label{req}\eeq
where
\bes\begin{align} \sigma&\text{ is a polynomial of degree }\leq2,\\
\tau&\text{ is a polynomial of degree }\leq1,\\
\eta&\text{ is a number,}\end{align}\label{conb}\ees
will be  called {\em hypergeometric class equations}.
Solutions of
(\ref{req})
will go under the name of {\em hypergeometric class functions}.   Operators
$\sigma(w)\p_w^2+\tau(w)\p_w+
\eta$ with $\sigma,\tau,\eta$ satisfying (\ref{conb}) will be called {\em hypergeometric class operators}.

Let us  review basic classes of hypergeometric class equations. 
We will always  assume that $\sigma(w)\neq0$.
Every class
 will be simplified
 by dividing by a constant and, except for (\ref{hy8}), by
 an affine change of the complex variable $w$.
 \medskip

\noindent
{\bf The ${}_2\cF_1$
  or  Gauss hypergeometric equation}
\beq
\left(w(1-w)\p_w^2+(c-(a+b+1)w)\p_w-ab\right)f(w)=0.\label{hy1}\eeq
\noindent{\bf The ${}_2F_0$  equation}
\beq\left(w^2\p_w^2+(-1+(1+a+b)w)\p_w+ab\right)f(w)=0.\label{hy2}\eeq
\noindent
{\bf The ${}_1\cF_1$ 
  or Kummer's  confluent equation}
 \beq
(w\p_w^2+(c-w)\p_w-a)f(w)=0.\label{hy3}\eeq
\noindent
{\bf The ${}_0\cF_1$ equation
}
 \beq
(w\p_w^2+c\p_w-1)f(w)=0.\label{hy4}\eeq
\noindent
{\bf The Hermite equation}
  \beq
(\p_w^2-2 w\p_w-2a)f(w)=0.\label{hy5}\eeq
\noindent
{\bf 2nd order Euler equation}
\beq
\left(w^2\p_w^2+bw\p_w+a\right)f(w)=0.\label{hy6}\eeq
\noindent
{\bf 1st order Euler equation for the derivative}
 \beq
(w\p_w^2+c\p_w)f(w)=0.\label{hy7}\eeq
\noindent
{\bf  2nd order equation with constant coefficients}
\beq
(\p_w^2+c\p_w+a)f(w)=0.\label{hy8}\eeq

Note that the equations (\ref{hy6}), (\ref{hy7}) and (\ref{hy8}) are elementary.
The remaining ones
(\ref{hy1}),  (\ref{hy2}),  (\ref{hy3}),  (\ref{hy4}) and  (\ref{hy5})
are the subject of these lecture notes.
 This is why they  are contained in  the list (\ref{hyp1})--(\ref{hyp5}) given at the beginning of these notes.
(Actually, (\ref{hy2}) is not explicitly mentioned in this list, however it is equivalent to (\ref{hy3}), so that these two equations are treated together).
 This list contains  also\\
\medskip
 
\noindent
{\bf The 
Gegenbauer equation}
\beq\left((1-w^2)\p_w^2-(a+b+1)w\p_w-ab\right)f(w)=0,\eeq
which can be reduced to a subclass of ${}_2\cF_1$ equations by a simple affine transformation. Its distinguishing property  is the invariance with respect to the reflection. The Gegenbauer equation has  special properties, which justify its separate treatment.

\section{(Pseudo-)Euclidean spaces}
\label{s3}
\init

In this section  we introduce basic terminology and notation related to  Lie algebras and groups acting on functions on $\rr^n$ or, more generally, on  manifolds. Lie algebras will be usually represented as 1st order differential operators. Lie groups will typically act as point transformations times  multipliers.

We will discuss
various operators   related to   (pseudo-)orthogonal Lie algebras and  groups. 
In particular, we will introduce a convenient notation to describe their Cartan algebras, root operators and  Weyl groups.
We will also discuss briefly the Laplacian and the Casimir operator.

We will show how to some special classes of harmonic functions---solutions of the Laplace equation. Of particular importance will be solutions that at the same time are eigenfunctions of the Cartan algebra. This construction will involve a contour integral, which can be viewed as a modification of the  Fourier or Mellin transformation. These solutions will be informally called {\em wave packets}.

Finally, in the last subsection we will show how to construct a harmonic function in $n{-}1$ dimension from a harmonic function in $n$ dimensions.


\subsection{Basic notation}

We  will write $\rr^\times$ for  $\rr\backslash\{0\}$,
$\rr_+$ for $]0,\infty[$ and $\rr_-$ for $]-\infty,0[$. We write $\cc^\times$ for $\cc\backslash\{0\}$.

We will  treat $\rr^n$ as a (real) subspace of $\cc^n$.
If possible, we will often extend functions from  real domains to holomorphic functions on complex domains.

In the following two subsections,
$\Omega,\Omega_1,\Omega_2$ are open subsets of $\rr^n$, or more generally, manifolds.

Often it is advantageous to  consider a similar formalism where $\Omega,\Omega_1,\Omega_2$ are open subsets of $\cc^n$, or more generally, complex manifolds. We will usually stick to the terminology typical for the real case. The reader can easily translate it to the complex picture, if needed.


\subsection{Point transformations with  multipliers}

 Let $\alpha:\Omega_1\to\Omega_2$ be a diffeomorphism.
The {\em transport} of functions  by the map
$\alpha$ will be also denoted by $\alpha$.\footnote{An alternative notation used often in mathematical literature for
the transport by $\alpha$ is  $\alpha_*$ or $(\alpha^*)^{-1}$.} More precisely, for
$f\in C^\infty(\Omega_1)$ we define ${\alpha}f\in C^\infty(\Omega_2)$ by
\[({\alpha}f)(y):=f({\alpha}^{-1}(y)).\]

If $m\in C^\infty (\Omega_2)$, then we have a map
$m{\alpha}: C^\infty(\Omega_1)\to C^\infty(\Omega_2)$ given by
\beq (m{\alpha} f)(y):=m(y)f({\alpha}^{-1}(y)).\label{lie1}\eeq
Transformations of the form (\ref{lie1}) will be called {\em point transformations with a multiplier}.

Clearly, transformations of the form (\ref{lie1}) with $\Omega=\Omega_1=\Omega_2$ and $m$ everywhere nonzero
form a group.

\subsection{1st order differential operators}

A vector field $X$ on $\Omega$ 
will be identified with the differential operator
\[Xf(y)=\sum_i X^i(y)\partial_{y^i}f(y),\ \ f\in C^\infty(\Omega),\]
where $X^i\in C^\infty(\Omega)$, $i=1,\dots,n$.
More generally, we will often use 1st order differential operators
\beq (X+M)f(y):=\sum_i X^i(y)\partial_{y^i}f(y)+M(y)f(y),\label{lie2}\eeq
where  $M\in C^\infty(\Omega)$.
Clearly, the set of operators of the form (\ref{lie2}) is a Lie algebra.

Let $\alpha:\Omega_1\to\Omega_2$ be a diffeomorphism.
If $X$ is a vector field on $\Omega_1$, then  ${\alpha}(X)$ is the vector field on $\Omega_2$  defined as
\[{\alpha}(X):=\alpha X  \alpha^{-1}.\]

\subsection{Affine linear transformations}

The  {\em general linear  group} is denoted $\GL(\rr^n)$.
It has a natural extension $\AGL(\rr^n):=\rr^n\rtimes \GL(\rr^n)$ called the {\em affine general linear group}. $(w,\alpha)\in\AGL(\rr^n)$ acts on $\rr^n$ by
\[\rr^n\ni y\mapsto w+\alpha y\in\rr^n.\]


 The permutation group $S_n$ can be naturally identified with a subgroup of  $\GL(\rr^n)$.
If $\pi\in S_n$, then 
\[(\pi y)^i:=y^{\pi_i^{-1}}.\]
On the level of functions, we have
\[\pi f(y^1,\dots,y^n)=f(y^{\pi_1},\dots,y^{\pi_n}).\]

The Lie algebra $\gl(\rr^n)$ represented by vector fields on $\rr^n$ is spanned by $y^i\partial_{y^j}$.

The Lie algebra 
  $\agl(\rr^n):=\rr^n\rtimes \gl(\rr^n)$ is spanned by $\gl(\rr^n)$ and by $\partial_{y_i}$.

A special element of $\gl(\rr^n)$ is  the {\em generator of dilations}, known also as the {\em Euler vector field},
\beq
A_n:=\sum_{i=1}^ny^i\partial_{y^i}.\label{euler}\eeq

We will often use the complex versions of the above groups, with $\rr$ replaced with $\cc$. We will  write $\GL(n)$ and $\gl(n)$, where the choice of the field  follows from the context.


\subsection{(Pseudo-)orthogonal group}

A {\em pseudo-Euclidean space} is $\rr^n$ equipped  with
  a symmetric nondegenerate $n\times n$ matrix  $g=[g_{ij}]$.
$g$ defines the   {\em scalar product} of vectors $x,y\in\rr^n$ and the {\em square} of a vector $x\in\rr^n$:
  \[\langle x|y\rangle:=\sum_{ij}x^i g_{ij}y^j,\quad
  \langle x|x\rangle=\sum_{ij}x^i g_{ij}x^j.\]

The matrix $[g^{ij}]$ will denote the inverse of $[g_{ij}]$.

 We will denote by $\ss^{n-1}(R)$ the   {\em   sphere
 in $\rr^n$ of
squared radius $R\in\rr$}:
\beq
\ss^{n-1}(R):=\{y\in\rr^n\ :\ \langle y|y\rangle=R\}.\label{sphere}\eeq
We will write $\ss^{n-1}:=\ss^{n-1}(1)$.

Actually, $\ss^{n-1}$  is  the  usual sphere only for the Euclidean signature. For non-Euclidean spaces it is a hyperboloid. Usually we will keep a uniform notation for all signatures. Occasionally, if we want to stress that $\ss^{n-1}$ has a specific signature,
it will be  denoted $\ss^{q,p-1}$,
where the signature of the ambient space is $(q,p)$ (see (\ref{signa})).

We also introduce   the {\em null quadric}
\beq \cV^{n-1}:=\ss^{n-1}(0)\backslash\{0\}.\label{null}\eeq

The {\em (pseudo-)orthogonal} and  the {\em special  (pseudo-)orthogonal group of $g$}
is defined as
\begin{eqnarray*}
{\rm O}(g)&:=&\{{\alpha}\in \GL(n)\ :\ \langle {\alpha}y|{\alpha}x\rangle=\langle y|x\rangle,\ y,x\in \rr^n\},\\
\SO(g)&:=&\{\alpha\in{\rm O}(g)\ :\ \det\alpha=1\}.\end{eqnarray*}
We also have
the {\em affine (special) orthogonal group} $\AO(g):=\rr^n\rtimes \O(g)$,
$\ASO(g):=\rr^n\rtimes \SO(g)$.

It is easy to see that the {\em pseudo-orthogonal Lie algebra}, represented by vector fields on $\rr^n$, can be defined by
\[\so(g):=\{B\in \gl(n)\ :\ B\langle y|y\rangle =0\}.\]
For $i,j=1,\dots,n$, define
\[B_{ij}:=\sum_k(g_{ik}y^k\p_{y^j}-g_{jk}y^k\p_{y^i}).\]
 $\{B_{ij}\ :\ i<j\}$ is a basis of $\so(g)$.
Clearly, $B_{ij}=-B_{ji}$ and $B_{ii}=0$.

The {\em affine pseudo-orthogonal Lie algebra}   $\aso(g):=\rr^n\rtimes \so(g)$ is spanned by $\partial_{y^i}$ and  $\so(g)$.

We will often use the complex versions of the above groups and Lie algebras.
In the real formalism we have to distinguish between various signatures of $g$---in the complex formalism there is only one signature and we can drop the prefix {\em pseudo}.

\subsection{Invariant operators}

Consider a pseudo-Euclidean space $\rr^n$.
 We define the Laplacian and the Casimir operator
\begin{eqnarray*}\Delta_n&:=&\sum_{i,j=1}^ng^{ij}\partial_{y^i}\partial_{y^j},\\
\cC_n&:=&\12\sum_{i,j,k,l=1}^ng^{ik}g^{jl}B_{ij}B_{kl}.
\end{eqnarray*}
 The above definitions do not depend on the choice of a basis.
$\Delta_n$ commutes with 
 $\AO(g)$ and  $\aso(g)$.
$\cC_n$ commutes with $\mathrm{O}(g)$ and  $\so(g)$.

Note the identity
\begin{eqnarray}
\langle y|y\rangle\Delta_n
&=&A_n^2+(n-2)A_n+\cC_n
,\label{polar}\end{eqnarray}
where $A_n$ is defined in (\ref{euler}).

\subsection{Orthonormal coordinates}

Suppose that $q+p=n$. Every scalar product of signature $(q,p)$ can be brought to the form
\beq
\langle y|y\rangle=-\sum_{i=1}^q y_i^2+\sum_{j=q+1}^{q+p}y_j^2.
\label{signa}\eeq
$\so(g)$ has a basis consisting of
\bes\begin{align}
  B_{ij}&=-y_i\partial_{y_j}+y_j\partial_{y_i},&&1\leq i<j\leq q\\
  B_{ij}&=y_i\partial_{y_{j}}+y_j\partial_{y_i}&&1\leq i\leq q,\quad q< j\leq n;\\
  B_{ij}&=y_i\partial_{y_j}-y_j\partial_{y_i},&&q< i<j\leq n.
\end{align}\ees
The Laplacian and the Casimir operator are
\begin{align}
  \Delta_n&=-\sum_{1\leq i\leq q} \partial_{y_i}^2+\sum_{q<j\leq n}\partial_{y_j}^2,\\
  \cC_n&=\sum_{1\leq i<j\leq q}B_{ij}^2+\sum_{q< i<j\leq n}B_{ij}^2
  -\sum_{\begin{array}{l}
   \scriptstyle   1\leq i\leq q,\\ \scriptstyle q< j\leq n
  \end{array}}B_{ij}^2.
\end{align}

We will rarely use orthonormal coordinates.

In the context of the signature $(q,p)$ the standard notation for the orthogoanl groups/Lie algebras is
 $\mathrm{O}(q,p)$,  $\AO(q,p)$,  $\so(q,p)$, $\aso(q,p)$.  We will however often use the notation $\mathrm{O}(n)$,  $\AO(n)$,  $\so(n)$, $\aso(n)$, without specifying the   signature of the quadratic form, and even allowing for an arbitrary choice of the field ($\rr$ or $\cc$).

\subsection{Split coordinates}

Suppose that $2m=n$.   $(m,m)$ will be called the {\em split signature}.
If the scalar product has such a signature, we can find coordinates such that
\beq\langle y|y\rangle=\sum_{i=1}^m2y_{-i}y_i.
\label{spli1}\eeq
We will say that
(\ref{spli1}) is a scalar product in {\em split coordinates}.

$\so(2m)$ has  a basis consisting of
\bes\begin{eqnarray} N_i:=B_{-ii}
&=&-y_{-i}\p_{y_{-i}}+y_{i}\p_{y_{i}}, \ \ \ j=1,\dots,m,\label{posa1}\\
B_{ij}&=&y_{-i}\p_{y_{j}}-y_{-j}\p_{y_{i}},\ \ \ 1\leq|i|<|j|\leq m.\label{posa2}\end{eqnarray}\ees

The subalgebra of $\so(2m)$ spanned by (\ref{posa1})
 is maximal commutative. It is called the {\em Cartan algebra} of $\so(2m)$.
(\ref{posa2}) are its {\em
   root operators}. They satisfy
\begin{eqnarray*}
  [N_k,B_{ij}]&=&-(\sgn(i)\delta_{k,|i|}+\sgn(j)\delta_{k,|j|})B_{ij}.
\end{eqnarray*}
The Laplacian and the Casimir operator are
\begin{align}
  \Delta_{2m}&=\sum_{i=1}^m 2\partial_{y_{-i}}\partial_{y_i},\\
    \cC_{2m}&=\sum_{1\leq |i|<|j|\leq m}B_{ij}B_{-i-j}-\sum_{i=1}^m N_i^2.
\label{casi1}\end{align}

Suppose now that $2m+1=n$.   In this case, $(m,m+1)$ will be called the {\em split signature}.
Every scalar product of such signature  can be brought to the form
\beq \langle y|y\rangle=y_0^2+\sum_{i=1}^m2y_{-i}y_i.
\label{spli2}\eeq
We will say that
(\ref{spli2}) is a scalar product in {\em split coordinates}.

 $\so(2m+1)$ has then
 a basis consisting of the above described basis of $\so(2m)$ and
\begin{eqnarray}
B_{0j}&=&
y_{0}\p_{y_{j}}-y_{-j}\p_{y_{0}},\ \ \ |j|=1,\dots,m.\label{posa3}
\label{posa22}\end{eqnarray}
The  additional roots satisfy
\begin{align}
  [N_k,B_{0j}]&=-\sgn(j)\delta_{k,|j|}B_{0j}.
\end{align}
The subalgebra  spanned by (\ref{posa1})
 is still maximal commutative in $\so(2m+1)$. It is called a {\em Cartan algebra} of $\so(2m+1)$.

We have
\begin{align}
  \Delta_{2m+1}&=\partial_{y_0}^2+\sum_{i=1}^m 2\partial_{y_{-i}}\partial_{y_i},\\
  \cC_{2m+1}&=\sum_{|i|=1}^mB_{0i}B_{0-i}+
  \sum_{1\leq |i|<|j|\leq m}B_{ij}B_{-i-j}-\sum_{i=1}^m N_i^2.\label{casi2}
\end{align}

In the real case we will most often consider the split signature, both in even and odd dimensions. In both real and complex cases we will  usually prefer split coordinates.
We will often write (\ref{spli1}) and  (\ref{spli2}) in the form
\beq\langle y|y\rangle=\sum_{|i|\leq m}y_{-i}y_i
\eeq
where it is understood that
$i\in\{-m,\dots,-1,1,\dots,m\}$ in the even case and
$i\in\{-m,\dots,-1,0,1,\dots,m\}$ in the odd case.

\subsection{Weyl group}
\label{Weyl symmetries}

In this subsection we introduce a certain finite subgroup of $\mathrm{O}(n)$, which will be called the {\em Weyl group}. We will also introduce a notation for elements of these groups.
The reader is referred to Subsects \ref{so(6) in 6 dimensions} and
\ref{sub-geg1}, for examples of application of  this notation.
 We will assume that the signature is split and split coordinates have been chosen.

Consider first dimension $2m$.
Permutations of $\{-1,\dots,-m\}\cup\{1,\dots,m\}$ that preserve the pairs
$\{-1,1\},\dots\{-m,m\}$ define elements of $\mathrm{O}(2m)$. They form a group, that we will call  denote $D_m$. It is isomorphic to $\zz_2^m\rtimes S_m$. It is the  Weyl group of $\mathrm{O}(2m)$.

The {\em flip} interchanging $-i,i$ will be denoted $\tau_i$.
 The flips $\tau_i$, with $i=1,\dots,m$, generate a subgroup of $D_m$ isomorphic to $\zz_2^m$.

To every $\pi\in S_m$ there corresponds an element of $D_m$
denoted $\sigma_\pi$, that permutes  pairs $(-i,i)$.
We have
\begin{eqnarray}\sigma_\pi f(y_{-1},y_{1},\dots,
y_{-m},y_m)&:=&
f(y_{-\pi_1},y_{\pi_1},\dots,
y_{-\pi_m},y_{\pi_m}).\label{per1}\end{eqnarray}

Let $\epsilon=(\epsilon_1,\dots,\epsilon_m)$ and
$\epsilon_1,\dots,\epsilon_m\in\{1,-1\}$. We will write
$\epsilon\pi$ as the shorthand for $\epsilon_1\pi_1,\dots,\epsilon_m\pi_m$.
We will
use the notation
\beq
\sigma_{\epsilon\pi}:=
\sigma_\pi\prod\limits_{\epsilon_j=-1}\tau_j
.\eeq

We have
 \begin{eqnarray*}
\sigma_{\epsilon\pi}B_{ij}\sigma_{\epsilon\pi}^{-1}
&= B_{\epsilon_i\pi_i,\epsilon_j\pi_j};\quad
\sigma_{\epsilon\pi} N_j\sigma_{\epsilon\pi}^{-1}&=\epsilon_jN_{\pi_j}.
\end{eqnarray*}

Using $\rr^{2m+1}=\rr\oplus\rr^{2m}$, we  embed
$D_m$ in ${\rm O}(2m+1)$.
We also introduce  $\tau_0\in {\rm O}(2m+1)$
given by
\beq \tau_0f(y_0,y_{-1},y_1,\dots,y_{-m},y_m)
:=f(-y_0,y_{-1},y_1,\dots,y_{-m},y_m).\label{pasa}\eeq
Clearly, $\tau_0$ commutes with $D_m$.
The group $B_m$ is defined as the group generated by $D_m$ and
$\tau_0$. It is isomorphic to $\zz_2\times\zz_2^m\rtimes S_m$.
It is the Weyl group of $\mathrm{O}(2m+1)$.
  
We set
\[\tau_{\epsilon\pi}:=\tau_0\sigma_{\epsilon\pi}.\]
We  have
\begin{eqnarray*}
\tau_{\epsilon\pi} B_{0j}\tau_{\epsilon\pi}^{-1}
=-B_{0,\epsilon_j\pi_j},\ \ &
\tau_{\epsilon\pi} B_{ij}\tau_{\epsilon\pi}^{-1}=B_{\epsilon_i\pi_i,\epsilon_j\pi_j},&\ \ 
\tau_{\epsilon\pi} N_j\tau_{\epsilon\pi}^{-1}=\epsilon_jN_{\pi_j}.\end{eqnarray*}


\subsection{Harmonic functions}
\label{Harmonic functions}

Suppose that $\rr^{n}$ is equipped with a scalar product.
We say that a function $F$ on $\rr^{n}$ is {\em harmonic} if
\beq \Delta_n F=0.\eeq

\bep Let $e_1,\dots e_k\in\rr^{n}$ satisfy
\[\langle e_i|e_j\rangle=0, \quad1\leq i,j\leq k.
\]
In other words, assume that $e_1,\dots,e_k$ span an {\em isotropic subspace} of $\rr^{n}$.
Let $f$ be a function of $k$ variables. Then
\[F(z):=f\big(\langle e_1|z\rangle,\dots,\langle e_k|z\rangle\big)\]
is harmonic.\label{isotro}\eep

For instance, consider $\rr^{n}$ with a split scalar product, where $n=2m$ or $n=2m+1$. Then any function $f(y_1,\dots,y_m)$ is harmonic, for instance
\beq
F_{\alpha_1,\dots\alpha_m}:=y_{1}^{\alpha_1}\cdots y_{m}^{\alpha_{m}},\label{fro}\eeq
which in addition satisfies
\beq  N_jF_{\alpha_1,\dots\alpha_m}=\alpha_jF_{\alpha_1,\dots\alpha_m}.\label{eigen}\eeq

Harmonic functions satisfying in addition the eigenvalue equations (\ref{eigen}) will play an important role in our approach.
Unfortunately, functions of the form (\ref{fro}) constitute a rather narrow class. We  need more general harmonic functions, which we will call {\em wave packets}. They are obtained by smearing a rotated
(\ref{fro})  with an appropriate weight, so that it is an eigenfunction of Cartan operators.
 This construction will be explained in the Subsect. 
\ref{wave1}--\ref{wave3}. It is essentially a version of the Fourier (or Mellin) transformation, possibly  with  a deformed complex contour of integration.

Note that the aim of Subsects \ref{wave1} and \ref{wave2} is to provide motivation, based on the concept of the Fourier transformation, for Subsect. \ref{wave3}, which contains the construction that will be used in what follows.

\subsection{Eigenfunctions of angular momentum I}
\label{wave1}

Suppose that $\rr^{n}=\rr^2\oplus\rr^{n-2}$, where we write
$z=(x,y,z')\in\rr^{n}$ and
\[\langle x,y,z'|x,y,z'\rangle=x^2+y^2+\langle z'|z'\rangle.\]
Set
\[N_1:=-\i(x\p_y-y\p_x).\]

Let  $m\in\zz$. Consider a function
$f(x,y,z')$. Then
 \begin{equation}
   F_m(x,y,z'):=\frac1{2\pi}\int_0^{2\pi} f(\cos\phi x-\sin\phi y,\sin\phi x+\cos\phi y,z')\e^{-\i m\phi}\d\phi,\label{porr1}
 \end{equation}
 \begin{equation}
 \label{porr1.}
 \text{satisfies }\quad N_1 F_m(x,y,z')= mF_m(x,y,z').\end{equation}
 Note that if $f$ is harmonic, then so is $F_m$. This construction is essentially the Fourier transformation.

Introduce complex coordinates
\beq
z_{\pm1}:=\frac{1}{\sqrt2}(x\pm\i y).\label{porr2}\eeq
We will write $f(z_{-1},z_1,z')=f(x,y,z')$, $F_m(z_{-1},z_1,z')=F(x,y,z')$.
The operator $N_1$ takes the familiar form
    \beq N_1=-z_{-1}\partial_{z_{-1}}+z_1\partial_{z_1},\label{kaj}\eeq
    and the metric becomes
    \beq\langle z_{-1},z_1,z'|z_{-1},z_1,z'\rangle=2z_{-1}z_1+\langle z'|z'\rangle.\label{porr3.}\eeq
    Then (\ref{porr1}) and (\ref{porr1.}) can be rewritten as
\begin{align}
  F_m(z_{-1},z_1,z')&:=\frac1{2\pi\i}\int_\gamma f(\tau^{-1} z_{-1},
  \tau z_1,z')\tau^{-m-1}\d\tau,\label{porr1a}\\
N_1F_m(z_{-1},z_1,z')&= mF_m(z_{-1},z_1,z'),\end{align}
where $\gamma$ is the closed contour
$[0,2\pi[\ni\phi\mapsto \tau=\e^{\i\phi}$.

    
\subsection{Eigenfunctions of angular momentum II}
\label{wave2}

We again consider $\rr^{n}=\rr^2\oplus\rr^{n-2}$, but  we change the signature of the metric. We assume that the scalar product   is given by
\beq\langle z_{-1},z_1,z'|z_{-1},z_1,z'\rangle=2z_{-1}z_1+\langle z'|z'\rangle.\label{porr3}\eeq
We start from a function $f(z_{-1},z_1,z')$. We would like to construct   an eigenfunction of $N_1$ with a generic eigenvalue $\alpha$, and not only with an integer eigenvalues as (\ref{porr1a}). To do this  we repeat a similar procedure as in the previous subsection. Now, however, we need to integrate over a  half-line, so we need conditions at the ends: we assume that
\beq f(\tau^{-1}z_{-1},\tau z_1,z')\tau^{-\alpha}\Big|_{\tau=0}^{\tau=\infty}=0.\eeq
We set
\begin{align}
F_\alpha&:=\frac1{2\pi\i}
\int_0^\infty f(\tau^{-1}z_{-1},\tau z_1,z')\tau^{-\alpha-1}\d\tau.
\end{align}
Then, with $N_1$ given by (\ref{kaj}),
\beq
N_1 F_\alpha(z_{-1},z_1,z')=\alpha F_\alpha(z_{-1},z_1,z').\label{porr1a-}\eeq

Indeed,
\begin{align*}
&  \partial_\tau
  f(\tau^{-1}z_{-1},\tau z_1,z')\tau^{-\alpha}\\=&
  -\alpha
  f(\tau^{-1}z_{-1},\tau z_1,z')\tau^{-\alpha-1}\\
  &-\tau^{-2}z_{-1}\partial_1f(\tau^{-1}z_{-1},\tau z_1,z')\tau^{-\alpha}
  +z_1\partial_2f(\tau^{-1}z_{-1},\tau z_1,z')\tau^{-\alpha}\\
  =&
 \big( -\alpha-z_{-1}\p_{z_{-1}}+z_{1}\p_{z_{1}}\big)
  f(\tau^{-1}z_{-1},\tau z_1,z')\tau^{-\alpha-1}.
\end{align*}
Hence
\begin{align}
  0&=\frac{1}{2\pi\i}\int_0^\infty\d\tau
 \partial_\tau
  f(\tau^{-1}z_{-1},\tau z_1,z')\tau^{-\alpha}\,=\,(-\alpha+N_1)F_\alpha.
  \end{align}
  
Note that $F_\alpha$ is  the Mellin transform of $\tau\mapsto
f(\tau^{-1}z_{-1},\tau z_1,z')$. If $f$ is harmonic, then so is $F_\alpha$.

\subsection{Eigenfunctions of angular momentum III}
\label{wave3}

Assume now that $z_{-1},z_1,z'$ are complex variables and
$f$ is holomorphic. Then we can formulate a result that includes (\ref{porr1.}) and (\ref{porr1a-}),  allowing for a greater flexibility of the choice of the contour of integration:

\bep\label{cac0}
Suppose that $]0,1[\ni s\overset{\gamma}\mapsto \tau(s)$ is a contour on the Riemann surface of
    \[\tau\mapsto
    f(\tau^{-1}z_{-1},\tau z_1,z')\tau^{-\alpha}\] that satisfies
    \beq
     f(\tau^{-1}z_{-1},\tau z_1,z')\tau^{-\alpha}\Big|_{\tau(0)}^{\tau(1)}=0.\eeq
Then
\beq
F_\alpha:=\frac1{2\pi\i}
\int_\gamma f(\tau^{-1}z_{-1},\tau z_1,z')\tau^{-\alpha-1}\d\tau
\eeq
solves
\[N_1 F_\alpha=\alpha F_\alpha.\]
\eep

\proof We repeat the arguments of the previous subsection, where we replace $[0,\infty[$ with $\gamma$. \qed

\subsection{Dimensional reduction}
\label{Dimensional reduction}

In this subsection we describe how to construct harmonic functions in $n-1$ dimensions out of a harmonic function in $n$ dimensions.

Suppose that $\rr^{n}$ is equipped with the scalar product
\[\langle z_{-1},z_{1},z'| z_{-1},z_{1},z'\rangle_{n}
=2z_{-1}z_{1}+\langle  z'| z'\rangle_{n-2}.\]
As usual,
we write
\begin{align}
  N_{1}&=-z_{-1}\p_{z_{-1}}+z_{1}\p_{z_{1}},\\
    \Delta_{n}&=2\p_{z_{-1}}\p_{z_{1}}+\Delta_{n-2}.
\end{align}
Introduce new variables and the Laplacian in $n{-}1$ dimensions.
\begin{align}
  z_0:&=\sqrt{2z_{-1}z_{1}},\quad
  u:=\sqrt{\frac{z_{1}}{z_{-1}}},\\
  \Delta_{n-1}&:=\partial_{z_0}^2+\Delta_{n-2}.
  \end{align}
  In the new variables,
  \begin{align}
        N_1&=u\p_u,\\
    \Delta_n&=\partial_{z_0}^2+\frac1{z_0}\partial_{z_0}
    -\frac1{z_0^2}(u\partial_u)^2+
    \Delta_{n-2}. \end{align}
Consequently,
\begin{align}
  z_0^{\frac12} \Delta_n z_0^{-\frac12}&=-\frac{1}{z_0^2}\Big(N_1-\frac12\Big)
  \Big(N_1+\frac12\Big)+\Delta_{n-1}.
  \end{align}
  Therefore, if we set
  \begin{align}
   F_\pm(z_0,u,z')&=u^{\pm\frac12}z_0^{-\frac12}f_\pm(z_0,z'),\end{align}
  then
  \begin{align}
    N_1F_\pm&=\pm\frac12F_\pm,\\
    z_0^{\frac12}u^{\mp\frac12}\Delta_nF_\pm&=\Delta_{n-1} f_\pm.
    \end{align}
 Hence, the $n-1$-dimensional Laplace equation
  $\Delta_{n-1} f=0$ is essentially equivalent to the $n$-dimensional Laplace equation   $\Delta_{n} F=0$ restricted to the eigenspace of $N_1=\pm\frac12$.

\section{Conformal invariance of the Laplacian}
\label{s5}
\init

Conformal manifolds are manifolds equipped with a conformal stucture---a  pseudo-Euclidean metric defined up to a positive multiplier. Conformal transformations are transformations that preserve the conformal structure.

The main  objects of this section 
are {\em projective null quadrics}. They possess a natural conformal structure with an exceptionally large group of conformal transformations.
 In fact, on the
 $n+2$ dimensional pseudo-Euclidean {\em ambient space} we have the obvious action of the pseudo-orthogonal Lie algebra and group. This action is inherited by the $n+1$ dimensional {\em null quadric} $\cV$, and then by its $n$-dimensional projectivization $\cY$.
 One can view $\cY$ as the base of the line bundle $\cV\to\cY$. By choosing a section $\gamma$ of this bundle we can equip $\cY$ with a pseudo-Riemannian structure. Choosing various sections defines metrics that differ only by a positive multiple---thus $\cY$ has a natural conformal structure. If the signature of the ambient space is $(q+1,p+1)$, then the signature of $\cY$ is $(q,p)$.

 We discuss a few examples of pseudo-Riemannian manifolds conformally equivalent to $\cY$ or to its open dense subset. The main example is the  flat pseudo-Euclidean space. Another example is  the product of two spheres $\ss^q\times\ss^p$, which is conformally equivalent to the entire
 $\cY$ of signature $(q,p)$.
 
 Especially simple and important are the low dimensional cases:
in $1$ dimension $\cY\simeq\ss^1$ and in $2$ dimensions $\cY\simeq\ss^1\times\ss^1$.
One should however remark that the dimensions 1 and 2 are somewhat special---
in these dimensions the  full  conformal Lie algebra is infinite dimensional, and the above construction gives only its subalgebra.

Conformal transformations are generalized symmetries of
the Laplacian. One can see this with help of a beautiful argument that goes back to Dirac.
Its first step is the construction of a certain geometrically defined operator denoted $\Delta_{n+2}^\diamond$, that transforms functions on $\cV$ homogeneous of degree $1-\frac{n}{2}$ into functions homogeneous of degree $-1-\frac{n}{2}$.
After fixing a section $\gamma$ of the line bundle
 $\cV\to\cY$, we can identify the somewhat abstract operator
$\Delta_{n+2}^\diamond$ with a concrete operator $\Delta_{n+2}^\gamma$ acting on fuctions on $\gamma(\cY)$. This operator turns out to be the Yamabe
Laplace-Beltrami operator for the corresponding pseudo-Riemannian structure.

On the $n+2$-dimensional ambient space the Laplacian  $\Delta_{n+2}$
 obviously
 commutes with the pseudo-orthogonal Lie algebra and group. On the level of $\gamma(\cY)$
 this commutation  becomes a transmutation of $ \Delta_{n+2}^\gamma$  with two different representations---one corresponding to the
 degree $1-\frac{n}{2}$, the other corresponding to the  degree $-1-\frac{n}{2}$.

 At the end of this section we consider in more detail the conformal action
 of the pseudo-orthogonal Lie algebra and group
  corresponding to the degree of homogeneity $\eta$
on  the flat pseudo-Euclidean space.
In particular, we compute the representations for all elements of the pseudo-orthogonal Lie algebra. For the pseudo-orthogonal group, we compute the representations of Weyl symmetries.

\subsection{Pseudo-Riemannian manifolds}

We say that a manifold $\cY$ is {\em pseudo-Riemannian} if it is equipped with a nondegenerate symmetric covariant 2-tensor
\[\cY\ni y\mapsto g(y)=[g_{ij}(y)]
,\]
called the {\em metric tensor}. For any vector field $Y$ it defines a function
$g(Y,Y)\in C^\infty(\cY)$:
\[\cY\ni y\mapsto g(Y,Y)(y):=g_{ij}(y)Y^i(y)Y^j(y).\]

Let $\alpha$ be a diffeomorphism of $\cY$.
As is well known,  the tensor $g$ can be transported by $\alpha$. More precisely, $\alpha^*(g)$ is defined by
\[\alpha^*(g)( Y,Y):=g\big(\alpha(Y),\alpha(Y)\big),\]
where $Y$ is an arbitrary vector field.
We say that $\alpha$ is {\em isometric} if $\alpha^* g=g$.

Let $X$ be a vector field.
The Lie derivative in the direction of
$X$ can be applied to the tensor $g$. More precisely, $\cL_Xg$ is defined by
\[(\cL_Xg)(Y,Y):=g\big([X,Y],Y\big)+
g\big(Y,[X,Y]\big).\]
  We say that a vector field $X$ is {\em Killing} if
  $\cL_Xg=0$.

  \subsection{Conformal manifolds}
  We say that the metric tensor $ g_1$ is {\em conformally equivalent to} $g$ if there exists a positive function
$m\in C^\infty(\cY)$      such that
      \[m(y)g(y)= g_1(y).\]
Clearly, the conformal equivalence is an equivalence relation in the set of metric tensors.      
We say that a manifold $\cY$ is equipped with a {\em conformal structure}, if it is equipped with an equivalence class of conformally equivalent metric tensors.

We say that a diffeomorphism  $\alpha$ is  {\em conformal} if for some metric tensor $g$ in the conformal class of $\cY$, $\alpha^* g$ is conformally equivalent to $g$. Clearly, this is equivalent to saying that for all $g$ in the conformal class of $\cY$, $\alpha^* g$ is conformally equivalent to $g$.

We say that a vector field $X$  is {\em conformal Killing} if for any metric tensors from the conformal class of $\cY$
there exists a smooth function
  $M\in C^\infty(\cY)$
  such that
  \beq \cL_Xg=Mg.\label{confi}
  \eeq
  Clearly, if (\ref{confi}) is true for one metric tensor $g$ from
  the conformal class of $\cY$, it is true for all  metric tensors conformally equivalent to $g$.
  \subsection{Projective null quadric}

  Consider a pseudo-Euclidean vector space $(\rr^{n+2},g)$
  of signature $(q+1,p+1)$, which we will call the {\em ambient space}. Recall that
  \[\cV^{n+1}:=\{z\in\rr^{n+2}\ :\ \langle z|z\rangle=0,\quad z\neq0\}.\]
  is the {\em null quadric}. For simplicity, we will often write $\cV$ for $\cV^{n+1}$.

The scaling, that is the action of  $\rr^\times$,  preserves $\cV$. Let $\cY:=\cV/\rr^\times$ be the
  {\em projective null quadric}.
  We obtain a {\em line bundle}
 $\cV\to \cY$
with the  base $\cY$ and the fiber $\rr^\times$.

  Let  $\cY_i$ be an open subset of $\cY$ and $\cV_i$ be the corresponding open subset of $\cV$. Let
  \[\cY_i\ni y\mapsto \gamma_i(y)\in\cV_i\]
  be a {\em section} of the bundle $\cV_i\to\cY_i$, that is a smooth map satisfying $y=\rr^\times\gamma_i(y)$.  Let $g_{\gamma_i}$ be the metric tensor $g$ restricted to $\gamma_i(\cY_i)$ transported to $\cY_i$.

  It is easy to prove the following fact:
  
  \bep Let $\gamma_i$, $i=1,2$, be  sections of $\cV_i\to\cY_i$. Then $g_{\gamma_i}$ are metrics on $\cY_i$ of signature $(q,p)$. The metrics $g_{\gamma_1}$ and $g_{\gamma_2}$ restricted to $\cY_1\cap\cY_2$ are conformally equivalent.
 \label{wqeq} \eep

Prop. \ref{wqeq}  equips $\cY$  with a conformal structure.

Choosing a section in the bundle $\cV\to\cY$  endows $\cY$ with the structure
of a  pseudo-Riemannian manifold.
For some special sections we obtain in particular various 
{\em symmetric spaces}
together with an explicit description of their conformal structure.
 In following subsections we present a few  examples of this construction.

Instead of $\cY$ one can consider $\tilde\cY:=\cV/\rr_+$. We obtain a bundle $\cV\to\tilde\cY$ with fibre $\rr_+$, which has similar properties as the bundle $\cV\to\cY$. It is a double covering of $\cY$, which means that we have a canonical $2-1$  surjection
$\tilde\cY\to\cY$.

Let $\gamma$ be a section of $\cV\to\cY$. Every $y\in\cY$ equals
$\rr^\times\gamma(y)$, and hence it is the disjoint union
of $\tilde y_+:=\rr_+\gamma(y)$ and $\tilde y_-:=\rr_-\gamma(y)$. Clearly $\{\tilde y_+,\tilde y_-\}\subset\tilde\cY$ is the preimage of $y$ under the canonical covering. Let us set
\beq\tilde\gamma(\tilde y_+):=\gamma(y),\quad
\tilde\gamma(\tilde y_-):=-\gamma(y).\eeq
Then $\tilde\gamma$ is a section of the bundle $\cV\to\tilde\cY$. With help of $\tilde\gamma$
we can equip $\tilde\cY$ with a metric $\tilde g_{\tilde\gamma}$. Obviously,  if $\cY$ is equipped with the metric
 $g_{\gamma}$, the canonical surjection 
$\tilde\cY\to\cY$ is isometric.

We would like to treat $\cY$ as the principal object, since it has a direct generalization to the complex case. However, for some purposes $\tilde\cY$ is preferable.

\subsection{Projective null quadric as a compactification of a pseudo-Euclidean space}
\label{Projective null quadric as a compactification of a (pseudo-)Euclidean space}

Consider a pseudo-Euclidean space $(\rr^n,g_n)$ of signature $(q,p)$  embedded in  the pseudo-Euclidean space $(\rr^{n+2},g_{n+2})$ of signature $(q+1,p+1)$.
We assume that
the square of a vector $(z',z_-,z_+)\in\rr^{n+2}=\rr^n\oplus\rr^2$   is
\[\langle z',z_-,z_+|z',z_-,z_+\rangle_{n+2}:=
\langle z'|z'\rangle_n+2z_+z_-.\]

Set
\[\cV_0:=\{(z',z_-,z_+)\in\cV\ :\ z_-\neq0\},\quad \cY_0:=\cV_0/\rr^\times.\]
$\cY_0$ is dense and open in $\cY$.

We have a bijection and a section
\beq\cY_0\ni\rr^\times \begin{bmatrix}y\\1\\-\frac{\langle y|y\rangle_n}{2}\end{bmatrix}\leftrightarrow \underset{\rr^n}{\underset{\bni}{y}}\mapsto
\begin{bmatrix}y\\1\\-\frac{\langle y|y\rangle_n}{2}\end{bmatrix}\in\cV_0.\label{bije}\eeq
Thus $\rr^n$ is identified with $\cY_0$. The metric on $\cY_0$ given by the above section coincides with
 the original metric on $\rr^n$. We have thus embedded $\rr^n$ with its conformal structure as a dense open subset of  $\cY$.

\subsection{Projective null quadric as a sphere/compactification of a hyperboloid}
\label{Projective null quadric as a sphere/compactification of a hyperboloid}
Consider a Euclidean space $(\rr^{n+1},g_{n+1})$ embedded in a pseudo-Euclidean space $(\rr^{n+2},g_{n+2})$ of signature $(1,n+1)$. We assume that the square of a vector $(z',z_0)\in\rr^{n+1}\oplus\rr=\rr^{n+2}$ is
\[\langle z',z_0|z',z_0\rangle_{n+2}=\langle z'|z'\rangle_{n+1}-z_0^2.
\]
Recall that
\[\ss^n:=\{\omega\in\rr^{n+1}\ :\ \langle \omega|\omega\rangle=1\}
\]
is the unit sphere
of dimension $n$.

We have a bijection and a section
\beq\cY\ni\rr^\times \begin{bmatrix}\omega\\1\end{bmatrix}
  \leftrightarrow \underset{\ss^n}{\underset{\bni}{y}}\mapsto \begin{bmatrix}\omega\\1\end{bmatrix}\in\cV.\label{bije2+}\eeq
Thus $\ss^n$ is identified with $\cY$. 
    The metric on $\cY$ given by the above section  coincides with
    the usual metric on $\ss^n$. 

$\tilde\cY$ is in this case simply the disjoint sum of two copies of $\ss^n$.
    
    The above construction can be repeated with minor changes for a general signature. Indeed, let the signature of $(\rr^{n+1},g_{n+1})$ be $(q,p+1)$,
    so that the signature of $(\rr^{n+2},g_{n+2})$ is  $(q+1,p+1)$.
 Set
\[\cV_0:=\{(z',z_0)\in\cV\ :\ z_0\neq0\},\quad \cY_0:=\cV_0/\rr^\times.\]
We have then the bijection and section
\beq\cY_0\ni\rr^\times \begin{bmatrix}\omega\\1\end{bmatrix}
  \leftrightarrow \underset{{\ss}^{q,p}}{\underset{\bni}{\omega}}\mapsto \begin{bmatrix}\omega\\1\end{bmatrix}\in\cV_0.\label{bije2++}\eeq
Note that now instead of the unit Euclidean sphere  we have  the unit hyperboloid of signature $(q,p)$, which has been identified with $\cY_0$, a dense open subset of $\cY$.

 \subsection{Projective null quadric as the Cartesian product of spheres}

 Consider now the space $\rr^{n+2}$ of signature $(q+1,p+1)$.
 The square of a vector $(\vec t,\vec x)=(t_0,\dots,t_q,x_0,\dots,x_p)$
 is defined as
 \beq
 \langle \vec t,\vec x|\vec t,\vec x\rangle:=
 -t_0^2-\cdots-t_q^2+x_0^2+\cdots+x_p^2.\label{signa1}\eeq

 Note that $\ss^q\times\ss^p$ is contained in $\cV$. It is easy to see that the map
\beq
\cY\ni \rr^\times
 (\vec\rho,\vec\omega)
\leftmapsto (\vec\rho,\vec\omega)\in
\ss^q\times\ss^p\subset \cV.\label{paio}\eeq
is a double covering. Indeed, we easily see that the map is onto and
 \[\rr^\times(\vec\rho,\vec\omega)=\rr^\times(-\vec\rho,-\vec\omega)
 .\]
 Thus
 \[\cY\simeq\ss^q\times\ss^p/\zz_2,\qquad
 \tilde\cY\simeq\ss^q\times\ss^p.\]


The map (\ref{paio}) can be interpreted as   a  section of $\cV\to \tilde\cY$.  
The corresponding metric tensor on $\cY$ is minus the standard metric tensor on $\ss^q$ plus the standard metric tensor on $\ss^p$. Its signature is $(q,p)$.

Again, similarly as in the previous subsection, the above construction can be generalized. Indeed,
replace (\ref{signa1}) with
 \begin{align*}
 \langle \vec t,\vec x|\vec t,\vec x\rangle:=& -t_0^2-\cdots-t_{q_1}^2+t_{q_1+1}^2+\cdots+t_{q_1+p_1}^2\\&+x_0^2+\cdots+x_{p_1}^2
 -x_{p_1+1}^2-\cdots x_{p_2+q_2}^2.\end{align*}
We then obtain a map
\beq
\cY\ni \rr^\times
 (\vec\rho,\vec\omega)
\leftmapsto (\vec\rho,\vec\omega)\in
\ss^{p_1,q_1}\times\ss^{q_2,p_2}\subset \cV.\label{paios}\eeq
Unlike (\ref{paio}), the map (\ref{paios}) is in general not onto---it doubly covers only an open dense subset of $\cY$.

\subsection{Dimension $n=1$}

Consider now the dimension $n=1$ in more detail.
The ambient space is
 $\rr^3$ with the split scalar product
\[\langle z|z\rangle=z_0^2+2z_{-1}z_{+1}.\]

The 1-dimensional projective quadric is isomorphic to $\ss^1$ or, what is the same, the 1-dimensional projective  space:
\[\cY^1\simeq\ss^1\simeq\rr\cup\{\infty\}=P^1\rr.\]

Indeed, it is easy to see that
\[\phi:\rr\cup\{\infty\}\to \cY^1\]
 defined by
\begin{eqnarray*}
 \phi(s)&:=&\Big(s,1,-\frac12s^2\Big)\rr^\times,\quad s\in\rr;\\
\phi(\infty)&:=&(1,0,0)\rr^\times
\end{eqnarray*}
is a homeomorphism.

The group $\mathrm{O}(1,2)$ acts on $P^1\rr$ by homographies (M\"obius transformations).

The Lie algebra $\so(1,2)$ is spanned by
\[B_{0,1},\ B_{0,-1},\ N_1,\]
with the commutation relations
\begin{eqnarray}\nonumber
[B_{0,1},B_{0,-1}]&=&N_1,\\
\nonumber
[B_{0,1},N_1]&=&B_{0,1},\\{}
[B_{0,-1},N_1]&=&-B_{0,-1}.\nonumber\label{comrel}
\end{eqnarray}
Appying (\ref{casi2}) with $m=1$ we obtain its Casimir operator:
\begin{subequations}\begin{eqnarray}
\cC_3&=&2B_{0,1}B_{0,-1}-N_1^2-N_1\\
&=&2B_{0,-1}B_{0,1}-N_1^2+N_1.
\end{eqnarray}\label{casimir1}\end{subequations}

\subsection{Dimension $ n=2$}

Consider finally the dimension $n=2$ in the signature $(1,1)$. The ambient space is
 $\rr^4$ with the split scalar product
\[\langle z|z\rangle=2z_{-1}z_{+1}+2z_{-2}z_{+2}.\]

The 2-dimensional projective quadric is isomorphic to the product of two circles:
\[\cY^2\simeq P^1\rr\times P^1\rr
.\]
Indeed, define
\[\phi:\big(\rr\cup\{\infty\}\big)\times \big(\rr\cup\{\infty\}\big)\to\cY^2\]
by
\bes\begin{eqnarray}
\phi(t,s)&:=&(-ts,1,t,s)\rr^\times,\label{kak1}\\
\phi(\infty,s)&:=&(-s,0,1,0)\rr^\times,\label{kak2}\\
\phi(t,\infty)&:=&(-t,0,0,1)\rr^\times,\label{kak3}\\
\phi(\infty,\infty)&:=&(-1,0,0,0)\rr^\times,\label{kak4}
\end{eqnarray}\ees
where $t,s\in\rr$. We easily check that $\phi$ is a homeomorphism.
In fact, rewriting  (\ref{kak1}) as
\begin{eqnarray*}
\phi(t,s)&=&\Big(-s,\frac1t,1,\frac{s}{t}\Big)\rr^\times\label{kak2a}\\
&=&\Big(-t,\frac1s,\frac{t}{s},1\Big)\rr^\times\label{kak3a}\\
&=&\Big(-1,\frac{1}{ts},\frac1s,\frac1t)\rr^\times,\label{kak4a}
\end{eqnarray*}
we see the continuity of $\phi$ at (\ref{kak2}), (\ref{kak3}), resp. (\ref{kak4}).


The Lie algebra $\so(2,2)$ is spanned by
\[N_1,\ N_2,\ B_{1,2},\ B_{1,-2},\ B_{-1,2},\ B_{-1,-2}.\]
Appying (\ref{casi2}) with $m=2$ we obtain its Casimir operator:
\[\cC_4=2B_{1,2}B_{-1,-2}+2B_{1,-2}B_{-1,2}-N_1^2-N_2^2-2N_1.\]

As is well known, $\so(2,2)$ decomposes into a direct sum
of two copies of $\so(1,2)$. Concretely,
\[\so(2,2)=\so^+(1,2)\oplus\so^-(1,2),\]
where $\so^+(1,2)$, resp. $\so^-(1,2)$, both
isomorphic to $\so(1,2)$, are
 spanned by
\begin{eqnarray*}
B_{1,2}, \ B_{-1,-2},\ N_1+N_2;&\text{resp.}&
B_{1,-2}, \ B_{-1,2},\ N_1-N_2.
\end{eqnarray*}
They have the commutation relations
\begin{eqnarray*}
\Big[\frac{B_{1,2}}{\sqrt2},\frac{B_{-1,-2}}{\sqrt2}\Big]=\frac{N_1+N_2}{2},&&\Big[\frac{B_{1,-2}}{\sqrt2},\frac{B_{-1,2}}{\sqrt{2}}\Big]=\frac{N_1-N_2}{2},
\\
{}\Big[\frac{N_1+N_2}{2},\frac{B_{-1,-2}}{\sqrt2}\Big]=\frac{B_{-1,-2}}{\sqrt2},&&\Big[\frac{N_1-N_2}{2},\frac{B_{-1,2}}{\sqrt2}\Big]=\frac{B_{-1,2}}{\sqrt2},
\\
{}\Big[\frac{N_1+N_2}{2},\frac{B_{1,2}}{\sqrt2}]=-\frac{B_{1,2}}{\sqrt2};&&\Big[\frac{N_1-N_2}{2},\frac{B_{1,-2}}{\sqrt2}\Big]=-2\frac{B_{1,-2}}{\sqrt2}.
\end{eqnarray*}
The corresponding Casimir operators are
\begin{eqnarray*}
\cC_3^{+}&=&B_{1,2}B_{-1,-2}-\frac14(N_1+N_2)^2-\frac12N_1-\frac12N_2\\
&=&B_{-1,-2}B_{1,2}-\frac14(N_1+N_2)^2+\frac12N_1+\frac12N_2,\\
\cC_3^{-}&=&B_{1,-2}B_{-1,2}-\frac14(N_1-N_2)^2-\frac12N_1+\frac12N_2\\
&=&B_{-1,2}B_{1,-2}-\frac14(N_1-N_2)^2+\frac12N_1-\frac12N_2.
\end{eqnarray*}
Thus
\[\cC_4=2\cC_3^{+}+2\cC_3^{-}.\]

In the enveloping algebra of $\so(2,2)$ the operators $\cC_3^{+}$ and
$\cC_3^{-}$ are distinct. They satisfy $\alpha(\cC_{-})=\cC_{+}$ for
$\alpha\in{\rm O}(2,2)\backslash\SO(2,2)$, for instance for
$\alpha=\tau_i$, $i=1,2$.

However, inside the associative algebra of differential operators on $\rr^4$ we have the identity
\[B_{1,2}B_{-1,-2}-B_{-1,2}B_{1,-2}=N_1N_2+N_1,\]
which implies
\[\cC_3^{+}=\cC_3^{-}\]
inside this algebra.
Therefore, represented in the algebra of differential operators we have
\begin{subequations}
\begin{eqnarray}
\cC_4&=&4B_{1,2}B_{-1,-2}-(N_1+N_2)^2-2N_1-2N_2\\
&=&4B_{-1,-2}B_{1,2}-(N_1+N_2)^2+2N_1+2N_2\\
&=&4B_{1,-2}B_{-1,2}-(N_1-N_2)^2-2N_1+2N_2\\
&=&4B_{-1,2}B_{1,-2}-(N_1-N_2)^2+2N_1-2N_2.
\end{eqnarray}\label{casimir2}
\end{subequations}

\subsection{Conformal invariance of the projective null quadric}
\label{s-diamond1}

 Obviously, $\mathrm{O}(n+2)$ and $\so(n+2)$ preserve $\cV$.
 They commute with the scaling (the action of $\rr^\times$). Therefore,
 we obtain the action on $\cY=\cV/\rr^\times$, which we denote as follows:
\bes\begin{eqnarray}
\so(n+2)\ni B&\mapsto &
 B^{\diamond},\label{popo1a.}\\
{\rm O}(n+2) \ni {\alpha}&\mapsto&
{\alpha}^{\diamond}.\label{popo2a.} \end{eqnarray}\ees
Clearly, the vector fields $B^\diamond$ are conformal Killing  and the diffeomorphisms
$\alpha^\diamond$ are conformal.


Let $\eta\in\cc$.
We define  $\Lambda_+^\eta(\cV)$
to be the set of smooth
functions on $\cV$ (positively) homogeneous of degree $\eta$, that is,
 satisfying
\[f(ty)=t^\eta f(y),\quad t> 0, \quad y\in\cV.\]

Clearly,  $B\in \so(n+2)$ and ${\alpha}\in {\rm O}(n+2)$ preserve $\Lambda_+^\eta(\cV)$. We will denote by $B^{\diamond,\eta}$, resp. $\alpha^{\diamond,\eta}$
the restriction of $B$, resp. $\alpha$ to  $\Lambda_+^\eta(\cV)$.
Thus we have representations
\bes\begin{eqnarray}
\so(n+2)\ni B&\mapsto &
 B^{\diamond,\eta},\label{popo1a}\\
{\rm O}(n+2) \ni {\alpha}&\mapsto&
{\alpha}^{\diamond,\eta},\label{popo2a} \end{eqnarray}\ees
acting on $\Lambda_+^\eta(\cY)$.

Clearly,
$\Lambda_+^0(\cV)$ can be identified with $C^\infty(\tilde\cY)$. Moreover,
(\ref{popo1a.}), resp. (\ref{popo2a.}) coincide with
(\ref{popo1a}), resp. (\ref{popo2a})   for $\eta=0$.

If $\eta\in\zz$ one can use another concept of homogeneity.
We define $\Lambda^\eta(\cV)$ to be
   the set of smooth
functions on $\cV$ satisfying
\[f(ty)=t^\eta f(y),\quad t\neq 0, \quad y\in\cV.\]
The properties of
$\Lambda^\eta(\cV)$ are similar to $\Lambda_+^\eta(\cV)$, except that 
$\Lambda^0(\cV)$ can be identified with $C^\infty(\cY)$.

\subsection{Laplacian on  homogeneous functions}
\label{s-diamond2}

The following theorem
according to Eastwood \cite{East} goes back to Dirac \cite{Dir}. We find it
curious because it allows in some
situations to restrict a {\em second order} differential operator to a
submanifold.

\bet Let $\Omega\subset\rr^{n+2}$ be an open conical set.
Let $ K\in C^\infty(\Omega)$ be homogeneous of degree
$1-\frac{n}{2}$
 such that \[ K\Big|_{\cV\cap\Omega}=0.\]
Then
\[\Delta_{n+2}  K\Big|_{\cV\cap\Omega}=0.\]
\label{wer}\eet


Before we give two proofs of this theorem,
 let us describe some of its consequences.

Let  $k\in\Lambda_+^{1-\frac{n}{2}}(\cV)$.
We can always find  $\Omega$,  a conical neighborhood of $\cV$,
and $ K\in\cA(\Omega)$  homogeneous of degree
$1-\frac{n}{2}$ such that
\[k= K\Big|_{\cV}.\]
Note that $\Delta_{n+2} K$ is homogeneous of degree
$-1-\frac{n}{2}$.
 We set
 \beq
 \Delta_{n+2}^\diamond k:=\Delta_{n+2} K\Big|_{\cV}.\label{restic}\eeq
By Theorem \ref{wer},
the above definition (\ref{restic}) does not depend on the choice of  $\Omega$ and  $ K$.
We  have thus defined a map
\beq
\Delta_{n+2}^\diamond:\Lambda_+^{1-\frac{n}{2}}(\cV)\to
\Lambda_+^{-1-\frac{n}{2}}(\cV).\label{diam}\eeq

Obviously,
\begin{subequations}\label{dada}
\begin{eqnarray} B\Delta_{n+2}&=&
\Delta_{n+2} B,\ \ \ B\in \so(n+2),\label{so33.}\\
{\alpha}\Delta_{n+2}&=&
\Delta_{n+2} {\alpha},\ \ \ {\alpha}\in
      {\rm O}(n+2).\label{SO33.}\end{eqnarray}
\end{subequations}

Restricting (\ref{dada})
to $\Lambda_+^{1-\frac{n}{2}}(\cV)$ we obtain
\begin{subequations}\label{dada1}
  \begin{eqnarray} B^{\diamond,-1-\frac{n}{2}}\Delta_{n+2}^\diamond&=&
\Delta_{n+2}^\diamond B^{\diamond,1-\frac{n}{2}},\ \ \ B\in \so(n+2),\label{so33}\\
{\alpha}^{\diamond,-1-\frac{n}{2}}\Delta_{n+2}^\diamond&=&
\Delta_{n+2}^\diamond {\alpha}^{\diamond,1-\frac{n}{2}},\ \ \ {\alpha}\in
      {\rm O}(n+2).\label{SO33}\end{eqnarray}
\end{subequations}

\noindent {\bf 1st proof of Thm \ref{wer}.}
We use the decomposition
    $\rr^{n+2}=\rr^n\oplus\rr^2$ described in Subsect.
\ref{Projective null quadric as a compactification of a (pseudo-)Euclidean space}, with the distinguished coordinates denoted $z_-,z_+$.
We denote the square of a vector, the Laplacian, the Casimir, resp. the generator of dilations
on $\rr^{n+2}$ by $R_{n+2}$, $\Delta_{n+2}$, $\cC_{n+2}$, resp. $A_{n+2}$. Similarly,
we denote the square of a vector, the Laplacian, the Casimir, resp. the generator of dilations
on $\rr^{n}$ by $R_n$, $\Delta_{n}$, $\cC_n$ resp. $A_{n}$.
We will also write
\[N_{m+1}:=z_+\partial_{z_+}-z_-\partial_{z_-}.\]

We have
\begin{align*}
  R_{n+2}&=
    R_n+2z_+z_-,\\
    \Delta_{n+2}&=\Delta_{n}+2\partial_{z_+}\partial_{z_-}
    ,\\
A_{n+2}&=A_{n}+z_+\partial_{z_+}+z_-\partial_{z_-}
.\end{align*}

 The following identity is a consequence  of (\ref{polar}):
\begin{eqnarray}\nonumber
R_n\Delta_{n+2}&=&
R_n\Delta_n+\big(R_{n+2}-2z_+z_-\big)2\partial_{z_+}\partial_{z_-}
\\\nonumber
&=&\cC_n+\Bigl(A_n-1+\frac{n}{2}\Bigr)^2-
\Bigl(\frac{n}{2}-1\Bigr)^2\\\nonumber
&&+R_{n+2} 2\partial_{z_+}\partial_{z_-}-
(z_+\partial_{z_+}+z_-\partial_{z_-})^2+
N_{m+1}^2\\\notag
&=&R_{n+2}  2\partial_{z_+}\partial_{z_-}
\\\notag&&+
\Bigl(A_n-1+\frac{n}{2}-z_+\partial_{z_+}-z_-\partial_{z_-}
\Bigr)
\Bigl(A_{n+2}-1+\frac{n}{2}\Bigr)\\
&&-
\Bigl(\frac{n}{2}-1\Bigr)^2+\cC_n+N_{m+1}^2.
\label{deq}\end{eqnarray}
    $\bigl(\frac{n}{2}-1\bigr)^2$ is a scalar.
      $\cC_n$ and $N_{m+1}^2$
  are polynomials  in elements of $\so(n+2)$, which are tangent to
$\cV$.
  Therefore,
all  operators in  the last line of (\ref{deq}) can be
restricted to $\cV$. The operator $ A_{n+2}-1+\frac{n}{2}$ vanishes
on functions in $\Lambda_+^{1-\frac{n}{2}}(\Omega)$.
The operator $R_{n+2} 2\partial_{z_+}\partial_{z_-}
$ is zero when restricted to $\cV$ (because $R_{n+2}$ vanishes on $\cV$).

Therefore,  if $ K$ is homogeneous of degree
$1-\frac{n}{2}$ vanishing on $\cV$, then
$R_n\Delta_{n+2}  K$ vanishes on $\cV$. We are free to choose different coordinates which give  different $R_n$'s. Therefore we can conclude that
$\Delta_{n+2}  K$ vanishes on $\cV$.
\qed

\begin{corollary}
Using the operator $\Delta_{n+2}^\diamond$, we can write
  \begin{eqnarray}
R_n\Delta_{n+2}^\diamond&=&
-\Bigl(\frac{n}{2}-1\Bigr)^2+\cC_n^{\diamond,1-\frac{n}{2}}+\big(N_{m+1}^{\diamond,1-\frac{n}{2}}\big)^2.
\label{deq1a}\end{eqnarray}
\end{corollary}

\noindent {\bf 2nd proof of Thm \ref{wer}.}
We use the decomposition $\rr^{n+2}=\rr^{n+1}\oplus\rr$ with the distinguished
  variable denoted by ${z_0}$, as in Subsect.
  \ref{Projective null quadric as a sphere/compactification of a hyperboloid}.
We denote the square of a vector, the Laplacian, the Casimir, resp. the generator of dilations
on $\rr^{n+1}$ by $R_{n+1}$, $\Delta_{n+1}$, $\cC_{n+1}$, resp. $A_{n+1}$. 
We have
\begin{align*}\label{2nd}
R_{n+2}&=R_{n+1}+z_0^2,\\
  A_{n+2}&=A_{n+1}+z_0\partial_{z_0},\\
\Delta_{n+2}&=\Delta_{n+1}+\partial_{z_0}^2.\end{align*}

We have the following identity
\begin{eqnarray}\nonumber
R_{n+1}\Delta_{n+2}&=&
R_{n+1}\Delta_{n+1}+\big(R_{n+2}-z_0^2\big)\partial_{z_0}^2\\\nonumber
&=&\cC_{n+1}+\Bigl(A_{n+1}+\frac{n-1}{2}\Bigr)^2-
\Bigl(\frac{n-1}{2}\Bigr)^2\\\nonumber
&&+R_{n+2}\partial_{z_0}^2- \Big(z_0\partial_{z_0}-\frac12\Bigr)^2+\Bigl(\frac12\Bigr)^2\\\notag
&=&R_{n+2} \partial_{z_0}^2+\Bigl(A_{ n+1}+\frac{n}{2}-z_0\p_{z_0}
\Bigr)
\Bigl(A_{n+2}+\frac{n}{2}-1\Bigr)\\
&&-\Bigl(\frac{n}{2}-1\Bigr)\frac{n}{2}+\cC_{n+1}.
\label{deq3}\end{eqnarray}
 Then we  argue similarly as in the 1st proof. \qed

\begin{corollary}
Using the operator $\Delta_{n+2}^\diamond$, we can write 
\begin{eqnarray}
R_{n+1}\Delta_{n+2}^\diamond&=&
-\Bigl(\frac{n}{2}-1\Bigr)\frac{n}{2}+\cC_{n+1}^{\diamond,1-\frac{n}{2}}.
\label{deq3a}\end{eqnarray}
\end{corollary}

\subsection{Fixing a section}

For nonzero $\eta$, in order to identify functions from $\Lambda_+^\eta(\cV)$ with functions on $\tilde\cY$ we need to fix a section of the line bundle $\cV\to\tilde\cY$.
Let us describe this in detail.

Let $\cV_0$ be an open homogeneous subset of $\cV$ and $\tilde\cY_0:=\cV_0/\rr_+$. Consider a section $\gamma:\tilde\cY_0\to\cV_0$. 
We then have the obvious identification
 $\psi^{\gamma,\eta}:\Lambda_+^\eta(\cV_0)\to C^\infty(\tilde\cY_0)$: for
$k\in \Lambda_+^\eta(\cV_0)$ we set
\beq\big(\psi^{\gamma,\eta} k\big)(y):=k\big( \gamma(y)\big),\ \ y\in\tilde\cY_0
.\label{psi}\eeq

The map $\psi^{\gamma,\eta}$ is bijective and we  can
introduce its
inverse, denoted  $\phi^{\gamma,\eta}$, defined  for any
$f\in C^\infty(\tilde\cY_0)$ by
\beq \big(\phi^{\gamma,\eta}f\big)\big(s\gamma(y)\big)=s^\eta f(y),\ \ s\in\rr_+,
\ \ y\in\tilde\cY_0.\label{phi}\eeq

Let
 $B\in \so(n+2)$ and $\alpha\in\mathrm{O}(n+2)$. As usual, $B$ and $\alpha$ are interpreted as transformations acting on functions  on $\rr^{n+2}$. Both $B$ and $\alpha$ preserve $\Lambda_+^\eta(\cV_0)$. Therefore, we can
define
\bes\begin{eqnarray}
B^{\gamma,\eta}&:=&\psi^{\gamma,\eta} B\phi^{\gamma,\eta},\label{poiu}\\
{\alpha}^{\gamma,\eta}&:=&\psi^{\gamma,\eta} {\alpha}\phi^{\gamma,\eta}.\label{poiu2}\end{eqnarray}\ees
$B^{\gamma,\eta}$ is a 1st order differential operator on $\tilde\cY_0$.
$\alpha^{\gamma,\eta}$ maps
$C^\infty\big(\tilde\cY_0\cap(\alpha^\diamond)^{-1}(\tilde\cY_0)\big)$
onto $C^\infty\big(\tilde\cY_0\cap \alpha^\diamond(\tilde\cY_0)\big)$.

It is easy to see that for any $B\in \so(n+2)$ and $\alpha\in\mathrm{O}(n+2)$
there exist $M_B\in C^\infty(\tilde\cY_0)$
and $m_\alpha
\in C^\infty\big(\tilde\cY_0\cap \alpha^\diamond(\tilde\cY_0)\big)$
such that 
\bes\begin{align}
  B^{\diamond,\eta}f(y)&=B^\diamond f(y)+\eta M_B(y) f(y),\label{popo1a..}\\
{\alpha}^{\diamond,\eta}f(y)&=m_\alpha^\eta(y)\alpha^\diamond f(y).\label{popo2a..} \end{align}\ees

We define also
\beq
\Delta_{n+2}^{\gamma}:=\psi^{\gamma,-1-\frac{n}{2}} \Delta_{n+2}^\diamond\phi^{\gamma,1-\frac{n}{2}}
.\eeq
This is a second order differential operator on $\tilde\cY_0$. It satisfies
\bes\begin{eqnarray} B^{\gamma,-1-\frac{n}{2}}\Delta_{n+2}^\gamma&=&
\Delta_{n+2}^\gamma B^{\gamma,1-\frac{n}{2}},\ \ \ B\in \so(n+2),\label{so33-}\\
{\alpha}^{\gamma,-1-\frac{n}{2}}\Delta_{n+2}^\gamma&=&
\Delta_{n+2}^\gamma {\alpha}^{\gamma,1-\frac{n}{2}},\ \ \ {\alpha}\in
      {\rm O}(n+2).\label{SO33-}\end{eqnarray}\ees

Note that for even $n$ the numbers $\pm1-\frac{n}{2}$ are 
integers. Therefore, $\Lambda^{\pm1-\frac{n}{2}}(\cV)$
are well defined. In the above construction, we can then use  $\cY$ instead of its double cover $\tilde\cY$. We also do not have problems in the complex case.

For odd $n$  the numbers $\pm1-\frac{n}{2}$ are not integers, and so
$\Lambda^{\pm1-\frac{n}{2}}(\cV)$ are ill defined.
Therefore,
we have to use $\Lambda_+^{\pm1-\frac{n}{2}}(\cV)$ 
and  $\tilde\cY$.




\subsection{Conformal invariance of the flat Laplacian}
\label{subsec-conf}

 In this subsection we illustrate
 the somewhat abstract theory of the previous subsections
with the example of  the 
 {\em flat section} described 
 in (\ref{bije}). Recall that the flat section identifies an open subset of $\cY$ with $\rr^n$. Therefore we obtain
 an action of  $\so(n+2)$ and $\mathrm{O}(n+2)$  on $\rr^n$.
As a result we will obtain  the invariance of the  Laplacian on the flat pseudo-Euclidean space with respect to conformal
 transformations.
 The results of this subsection will be needed  for our discussion of symmetries of the heat equation.

 We will use the notation of
 (\ref{poiu}) and  (\ref{poiu2}),
where instead of $\gamma$ we 
   write
 ``$\fl$'', for the
   {\em flat section}.
   We will describe   conformal symmetries 
 on two levels:
 \ben\item[(a)] the ambient space
 $\rr^{n+2}$
\item[(b)]
  the space $\rr^n$.\een

 We will use the split coordinates,
 that is,  $ z\in\rr^{n+2}$ and $y\in \rr^n$
 have the square
\bes\begin{eqnarray}
 \langle z|z\rangle&=&\sum_{|j|\leq m+1}z_{-j}z_j,\\
\langle y|y\rangle &=&\sum_{|j|\leq m}y_{-j}y_j
.\end{eqnarray}
\ees

As a rule, if a given operator does not depend on $\eta$, we omit the subscript $\eta$.

Derivation of all the following identities will be sketched in  Subsect. \ref{Computations}.  

\medskip

\noindent
{\bf Cartan algebra of $\so(n+2)$}

\medskip

\noindent
Cartan operators  of $\so(n)$, $i=1,\dots,m$:
\bes\begin{eqnarray}
N_i&=&-z_{-i}\p_{z_{-i}}+z_{i}\p_{z_i},\\
N_i^{\fl}&=&-y_{-i}\p_{y_{-i}}+y_{i}\p_{y_i}.
\end{eqnarray}\ees
Generator of dilations:
\bes\begin{eqnarray}
N_{m+1}&=&-z_{-m-1}\p_{z_{-m-1}}+z_{m+1}\p_{z_{m+1}},\\
N_{m+1}^{\fl ,\eta}&=&\sum\limits_{|i|\leq m}y_i\p_{y_i}-\eta\ =\ A_n-\eta.
\end{eqnarray}\ees

\medskip

\noindent{\bf Root operators}

\medskip

\noindent
Roots of $\so(n)$,
$|i|<|j|\leq m$:
\bes\begin{eqnarray}
B_{i,j}&=&z_{-i}\p_{z_j}-z_{-j}\p_{z_i},\\
B_{i,j}^{\fl }&
=&y_{-i}\p_{y_j}-y_{-j}\p_{y_i}.
\end{eqnarray}\ees
Generators of translations, $|j|\leq m$:
\bes\begin{eqnarray}
B_{m+1,j}&=&z_{-m-1}\p_{z_j}-z_{-j}\p_{z_{m+1}},\\
B_{m+1,j}^{\fl }&=&\p_{y_j}.
\end{eqnarray}\ees
Generators of special conformal transformations,  $|j|\leq m$:
\bes\begin{eqnarray}
B_{-m-1,j}&=&z_{m+1}\p_{z_j}-z_{-j}\p_{z_{-m-1}},\\
B_{-m-1,j}^{\fl ,\eta}&=&-\12 \langle y|y\rangle \p_{y_j}
+y_{-j}\sum\limits_{|i|\leq m}y_i\p_{y_i}-\eta y_{-j}.
\end{eqnarray}\ees

\medskip

\noindent{\bf Weyl symmetries}

\noindent We will write $K$ for a function on
$\rr^{n+2}$ and $f$ for a function on
$\rr^n$. We only give some typical elements that generate the whole Weyl group.

\medskip

\noindent
Reflection in the $0$th coordinate (for odd $n$):
\bes\begin{eqnarray}
\tau_0
 K(z_0,\dots)&=&
 K(-z_0,
\dots),\\
\tau_0^{\fl }f(y_0,\dots)&&
=\ f(-y_0,\dots
).
\end{eqnarray}\ees
Flips, $j=1,\dots,m$:
\bes\begin{align}
\tau_j
 K(\dots,z_{-j},z_j,\dots,z_{-m-1},z_{m+1})&\notag\\
&\hspace{-20ex}=
 K(
\dots,
z_{j},z_{-j},\dots,z_{-m-1},z_{m+1}),\\
\tau_j^{\fl}f(\dots,y_{-j},y_j,\dots)
=& f(\dots
y_{j},y_{-j},\dots).
\end{align}
\ees
Inversion:
\bes\begin{align}
\tau_{m+1} K(\dots,z_{-m-1},z_{m+1})
&=  K(\dots,z_{m+1},z_{-m-1}),\\
\tau_{m+1}^{\fl,\eta} f(y)&=
\Bigl(-\frac{\langle y|y\rangle }{2}\Bigr)^\eta f\Bigl(-\frac{2y}{\langle y|y\rangle }\Bigr).
\end{align}\ees
Permutations, $\pi\in S_m$:
\bes\begin{align}
\sigma_\pi
K(\dots,z_{-j},z_j,\dots,z_{-m-1},z_{m+1})&\notag\\
&\hspace{-14ex}
= K(
\dots,
z_{-\pi_j},z_{\pi_j},\dots,z_{-m-1},z_{m+1}),\\
\sigma_\pi^{\fl }f(\dots,y_{-j},y_j,\dots)&
=\ f(\dots
y_{-\pi_j},y_{\pi_j},\dots).
\end{align}\ees
Special conformal transformations, $j=1,\dots,m$:
\bes\begin{align}
  \sigma_{(j,m+1)} K(z_{-1},z_1,\dots,z_{-j},z_j,\dots,z_{-m-1},z_{m+1})&\notag\\
&\hspace{-30ex} =
 K(z_{-1},z_1,\dots,z_{-m-1},z_{m+1},\dots,z_{-j},z_j),
\\
\sigma_{(j,m+1)}^{\fl ,\eta} f
(y_{-1},y_{1},\dots,y_{-j},y_j,\dots)&\notag\\
&\hspace{-30ex}=\quad y_{-j}^\eta
f\Bigl(
\frac{y_{-1}}{y_{-j}},\frac{y_{1}}{y_{-j}},\dots,\frac{1}{y_{-j}},-\frac{\langle y|y\rangle
}{2y_{-j}}
\dots\Bigr).\end{align}\ees

\medskip

\noindent{\bf Laplacian}

\medskip

\bes\begin{eqnarray}
\Delta_{n{+}2}&=&\sum\limits_{|i|\leq m+1}\p_{z_i}\p_{z_{-i}},\\
\Delta_{n+2}^\fl&=&\sum\limits_{|i|\leq m}\p_{y_i}\p_{y_{-i}}\ =\ \Delta_n.
\end{eqnarray}\ees

We have the representations
on functions on $\rr^n$:
\bes\begin{eqnarray}
\so(n+2)\ni B&\mapsto&B^{\fl
  ,\eta}
,\label{pos1}\\
{\rm O}(n+2)\ni
 {\alpha}&\mapsto &{\alpha}^{\fl
   ,\eta}.
\label{pos2}\end{eqnarray}
\ees

They yield generalized symmetries:

\bes
\begin{eqnarray}
B^{\fl ,\frac{-2-n}{2}}\Delta_n&=&\Delta_n B^{\fl ,\frac{2-n}{2}},\ \ \ B\in \so(n+2),\\
{\alpha}^{\fl ,\frac{-2-n}{2}}\Delta_n&=&\Delta_n
{\alpha}^{\fl ,\frac{2-n}{2}},\ \ \ {\alpha}\in {\rm O}(n+2).
\end{eqnarray}
\ees

\medskip

\subsection{Computations}
\label{Computations}

Below we sketch explicit computations that lead to
the formulas on from the previous subsection.
Consider $\rr^{n}\times\rr^\times\times\rr$ (defined by
$z_{-m-1}\neq0$),
 which is an open dense subset of
$\rr^{n+2}$. Clearly, $\cV_0$ is contained in
 $\rr^{n}\times\rr^\times\times\rr$.

We will write $\Lambda^\eta(\rr^{n}\times\rr^\times\times\rr)$ for the space
of functions {\em homogeneous} of degree $\eta$ on
$\rr^{n}\times\rr^\times\times\rr$.

Instead of using the maps $\phi^{\fl ,\eta}$ and $\psi^{\fl ,\eta}$,
as   in (\ref{phi}) and (\ref{psi}),
we will prefer
 $\Phi^{\fl ,\eta}:C^\infty(\rr^n)\to\Lambda^\eta(\rr^{n}\times\rr^\times\times\rr)$ and
 $\Psi^{\fl ,\eta}:\Lambda^\eta(\rr^{n}\times\rr^\times\times\rr)\to C^\infty(\rr^n)$
defined below.

For $ K\in
\Lambda^\eta\bigl(\rr^{n}\times\rr^\times\times\rr\bigr)$, we define $\Psi^{\fl,\eta} K\in C^\infty(\rr^n)$ by
\[\big(\Psi^{\fl,\eta} K\big)(y)= K\Big(y,1,-\frac{\langle y|y\rangle }{2}\Big),  \ \ y\in \rr^n.\]

Let $f\in C^\infty(\rr^n)$.
Then there exists a unique function in
$\Lambda^\eta\bigl(\rr^{n}\times\rr^\times\times\rr\bigr)$
that extends $f$ and does not depend on $z_{m+1}$. It is given by
\[\big(\Phi^{\fl ,\eta} f\big)
(z,z_{-m-1},z_{m+1})
:=z_{-m-1}^\eta
f\Big(\frac{z}{z_{-m-1}}\Big),\quad z\in\rr^n.\]

The map $\Psi^{\fl,\eta}$ is a left inverse of   $\Phi^{\fl ,\eta}$:
\begin{eqnarray*}\Psi^{\fl,\eta}\Phi^{\fl ,\eta}& =&\id,\end{eqnarray*}
where $\id$ denotes the identity. Clearly,
\begin{eqnarray*}\Phi^{\fl ,\eta}f\Big|_{\cV_0}&=&\phi^{\fl ,\eta}f,\\
\Psi^{\fl,\eta}K&=&\psi^{\fl,\eta}\Big(K\Big|_{\cV_0}\Big).\end{eqnarray*}
Moreover, functions in $\Lambda^\eta(\rr^{n}\times\rr^\times\times\rr)$ restricted to $\cV_0$ are in $\Lambda^\eta(\cV_0)$.
Therefore,
\begin{eqnarray*}
B^{\fl ,\eta}&=&\Psi^{\fl,\eta} B\Phi^{\fl ,\eta},\ \ \ \ B\in \so(\rr^{n+2}),\\
{\alpha}^{\fl ,\eta}&=&\Psi^{\fl,\eta} {\alpha}\Phi^{\fl ,\eta},\ \ \ {\alpha}\in
{\rm O}(\rr^{n+2}).\end{eqnarray*}
(Note that $\alpha,B$ preserve
 $\Lambda^\eta(\rr^{n}\times\rr^\times\times\rr)$).
Note also that
\[\Delta_{n+2}^\fl=\Psi^{\fl,\eta} \Delta_{n+2}\Phi^{\fl
  ,\eta}=\Delta_n
.\]

In practice, the above idea can be implemented by
 the following change
of coordinates on $\rr^{n+2}$:
\[\begin{array}{l}
y_i:=\frac{z_i}{z_{-m-1}},\ |i|\leq m,\\[3mm]
R:=\sum\limits_{|i|\leq m+1}z_{i}z_{-i},\\[3mm]
p:=z_{-m-1}.
\end{array}\]
The inverse transformation is
\[\begin{array}{l}
z_i=py_i,\ |i|\leq m,\\[3mm]
z_{m+1}=\12(\frac{R}{p}-p\sum\limits_{|i|\leq m}y_iy_{-i}),\\[3mm]
z_{-m-1}=p
\end{array}\]
The derivatives are equal to
\[\begin{array}{l}
\p_{z_i}=z_{-m-1}^{-1}\p_{y_i}+2z_{-i}\p_R,\ |i|\leq m,\\[3mm]
\p_{z_{m+1}}=2z_{-m-1}\p_R,\\[3mm]
\p_{z_{-m-1}}=\p_p-z_{-m-1}^{-2}
\sum\limits_{|i|\leq m}z_i\p_{y_i}+2 z_{m+1}\p_R.
\end{array}\]

Note that these coordinates are defined on
$\rr^{n}\times\rr^\times\times\rr$.
The set $\cV_0$ is given by
the condition $R=0$. The flat section is given by  $p=1$.

For a function $y\mapsto f(y)$ we have
\[\big(\Phi^{\fl ,\eta} f\big)(y,R,p)=p^\eta f(y).\]
For a function
$(y,R,p)\mapsto K(y,R,p)$ we have
\[\big(\Psi^{\fl,\eta} K\big)(y)= K(y,1,0).\]
Note also that on $\Lambda^\eta(\rr^{n}\times\rr^\times\times\rr)$ we have
\[p\p_p+2R\p_R=\eta.\]

\section{Laplacian in 4 dimensions and the hypergeometric equation}
\label{s6}
\init

The goal of this section is to derive the ${}_2\cF_1$  equation together with its symmetries from the Laplacian in  $4$ dimensions, or actually from the Laplacian in $6$ dimensions, if one takes into account the ambient space.  Let us describe the main steps of this derivation:
\ben
\item\label{it1}
  We start from the  $4+2=6$ dimensional ambient space, with the obvious representations of $\so(6)$ and $\mathrm{O}(6)$, and the Laplacian $\Delta_6$.
  \item\label{it2} As explained      in Subsect. \ref{s-diamond1},
we introduce the representations
    $\so(6)\ni B\mapsto B^{\diamond,\eta}$ and
$\mathrm{O}(6)\ni \alpha\mapsto \alpha^{\diamond,\eta}$. Besides,
 as explained
 in Subsect. \ref{s-diamond2},
we obtain the reduced Laplacian $\Delta_6^\diamond$.
The most relevant values of $\eta$ are $1-\frac{4}{2}=-1$ and  $-1-\frac{4}{2}=-3$, which yield  generalized symmetries of $\Delta_6^\diamond$.
\item\label{it3} We fix a section $\gamma$ of the null quadric.
It allows us to construct the representations $B^{\gamma,\eta}$, $\alpha^{\gamma,\eta}$ and  the operator $\Delta_6^\gamma$, acting on a  $4$ dimensional manifold whose pseudo-Riemannian structure depends on $\gamma$.
\item\label{it4} We choose coordinates $w,u_1,u_2,u_3$, so that the Cartan operators are expressed in terms of $u_1$, $u_2$, $u_3$. We compute
  $\Delta_6^\gamma$, $B^{\gamma,\eta}$, and  $\alpha^{\gamma,\eta}$ in the new coordinates.
\item\label{it5} We make an ansatz that diagonalizes the Cartan operators, whose eigenvalues,  denoted by $\alpha$, $\beta$, $\mu$, 
  become parameters. 
  $\Delta_6^\gamma$,  $B^{\gamma,\eta}$, and  $\alpha^{\gamma,\eta}$ involve now only the single variable $w$. 
  $\Delta_6^\gamma$ turns out to be the  ${}_2\cF_1$ hypergeometric operator.
  The generalized symmetries of
  $\Delta_6^\gamma$ yield transmutation relations and discrete symmetries of the
 ${}_2\cF_1$  operator.
\een

Step \ref{it1} is described in Subsect. \ref{so(6) in 6 dimensions}.

We have a considerable freedom in the choice of the section $\gamma$
of Step \ref{it3}. For instance, it can be the flat section, which we described in Subsects
\ref{Projective null quadric as a compactification of a (pseudo-)Euclidean space} and \ref{subsec-conf}. However, to simplify computations we prefer to choose a different section, which we call the spherical section. (Both approaches are described in   \cite{DeMaj}).

We perform Steps \ref{it2}, \ref{it3} and \ref{it4} at once.
They are  described jointly in Subsect.
\ref{so(6) on the  spherical section}. We choose coordinates $w,r,p,u_1,u_2,u_3$ in $6$ dimensions, so that the null quadric, the spherical section and the homogeneity of functions are expressed in a simple way. 
In these coordinates, after the reductions of Steps \ref{it2} and \ref{it3}, the variables $r,p$ disappear. We are left with the variables $w,u_1,u_2,u_3$, and we are ready for Step \ref{it5}.

 Step \ref{it5} is described in Subsects
\ref{Hypergeometric equation} and 
\ref{Transmutation relations and discrete symmetries}.

Subsects \ref{Factorizations of the Laplacian} and
\ref{Factorizations of the hypergeometric operator}
are devoted to factorizations of the ${}_2\cF_1$ operator. Again, we see that the additional dimensions make all the formulas  more symmetric.
The role of factorizations is explained in Subsect.
\ref{Factorization relations}.

Subsects \ref{Transmutation relations and discrete symmetries} and
\ref{Factorizations of the hypergeometric operator} contain long lists of identities for the hypergeometric operator. We think that it is easy to appreciate and understand them at a glance, without studying them  line by line.
Actually, the analogous lists of identities in the next sections, corresponding to other types of equations, are shorter but in a sense more complicated, because they correspond to ``less symmetric'' groups.

All the material so far has been devoted to the ${}_2\cF_1$ operator
and its multidimensional ``parents''. Starting with Subsect. \ref{The ${}_2F_1$ hypergeometric function} we discuss  the  ${}_2F_1$ function and, more generally, distinguished solutions of the ${}_2\cF_1$ equation. The symmetries of the 
${}_2\cF_1$ operator are helpful in deriving and organizing the identities concerning these solutions.

Subsects \ref{Wave packets in 6 dimensions},
\ref{Integral representations}, \ref{Integral representations of standard solutions} are devoted to integral representations of solutions of the 
${}_2\cF_1$ equation. In particular, Subsect.  \ref{Wave packets in 6 dimensions}
shows that these representation are disguised ``wave packets'' solving the Laplace equation and diagonalizing Cartan operators.

In Subsect. \ref{Connection formulas} we derive connection formulas, where we use the pairs of solutions with a simple behavior at $0$  and at $\infty$ as two bases of solutions. The connection formulas follow easily from integral representations. These identities look  symmetric when expressed in terms of the Lie-algebraic parameters.

\subsection{$\so(6)$ in 6 dimensions}
\label{so(6) in 6 dimensions}

We consider $\rr^6$ with the split coordinates
\beq
z_{-1},z_1,z_{-2},z_2,z_{-3},z_3\label{sq0}\eeq
and the scalar product given by
\beq
\langle z|z\rangle=2z_{-1}z_1+2z_{-2}z_2+2z_{-3}z_3.\label{sq1}\eeq

The Lie algebra  $\so(6)$ acts naturally on $\rr^6$. Below we describe its natural basis. Then we consider its Weyl group, $D_3$,  acting on functions on $\rr^6$. For brevity, we list only elements
from its subgroup  $D_3\cap \SO(6)$. Finally, we write down the Laplacian.

\medskip

\noindent{\bf Lie algebra $\so(6)$.}
Cartan algebra
\bes\begin{align}
  N_1&=-z_{-1}\partial_{z_{-1}}+z_{1}\partial_{z_{1}},\\
  N_2&=-z_{-2}\partial_{z_{-2}}+z_{2}\partial_{z_{2}},\\
   N_3&=-z_{-3}\partial_{z_{-3}}+z_{3}\partial_{z_{3}}.
  \end{align}\ees
Root operators
\begin{subequations}
\begin{align}
  B_{-2,-1}&=z_{2}\p_{z_{-1}}-z_{1}\p_{z_{-2}},\\
  B_{2,1}&=z_{-2}\p_{z_1}-z_{-1}\p_{z_2},\\
B_{2,-1}&=z_{-2}\p_{z_{-1}}-z_{1}\p_{z_2},\\
B_{-2,1}&=z_{2}\p_{z_1}-z_{-1}\p_{z_{-2}};
\end{align}
\begin{align}  B_{-3,-2}&=z_{3}\p_{z_{-2}}-z_{2}\p_{z_{-3}},\\
  B_{3,2}&=z_{-3}\p_{z_2}-z_{-2}\p_{z_3},\\
B_{3,-2}&=z_{-3}\p_{z_{-2}}-z_{2}\p_{z_3},\\
B_{-3,2}&=z_{3}\p_{z_2}-z_{-2}\p_{z_{-3}};
\end{align}
\begin{align}
  B_{-3,-1}&=z_{3}\p_{z_{-1}}-z_{1}\p_{z_{-3}},\\
  B_{3,1}&=z_{-3}\p_{z_1}-z_{-1}\p_{z_3},\\
B_{3,-1}&=z_{-3}\p_{z_{-1}}-z_{1}\p_{z_3},\\
B_{-3,1}&=z_{3}\p_{z_1}-z_{-1}\p_{z_{-3}}.
\end{align}
\end{subequations}

\noindent{\bf Weyl symmetries}

\bes
\begin{align}
\sigma_{123} K(z_{-1},z_1,z_{-2},z_2,z_{-3},z_{3})=&K(z_{-1},z_1,z_{-2},z_2,z_{-3},z_{3}),\\
\sigma_{-12-3}K(z_{-1},z_1,z_{-2},z_2,z_{-3},z_{3})=&K(z_{1},z_{-1},z_{-2},z_2,z_{3},z_{-3}),\\
\sigma_{1-2-3}K(z_{-1},z_1,z_{-2},z_2,z_{-3},z_{3})=
&K(z_{-1},z_{1},z_{2},z_{-2},z_{3},z_{-3}),\\
\sigma_{-1-23}K(z_{-1},z_1,z_{-2},z_2,z_{-3},z_{3})=&K(z_{1},z_{-1},z_{2},z_{-2},z_{-3},z_{3});\end{align}
\begin{align}
\sigma_{213}K(z_{-1},z_1,z_{-2},z_2,z_{-3},z_{3})=&K(z_{-2},z_2,z_{-1},z_1,z_{-3},z_3),\\
\sigma_{-21-3}K(z_{-1},z_1,z_{-2},z_2,z_{-3},z_{3})=&K(z_{2},z_{-2},z_{-1},z_1,z_{3},z_{-3}),\\
\sigma_{2-1-3}K(z_{-1},z_1,z_{-2},z_2,z_{-3},z_{3})=&K(z_{-2},z_2,z_{1},z_{-1},z_{3},z_{-3}),\\
\sigma_{-2-13}K(z_{-1},z_1,z_{-2},z_2,z_{-3},z_{3})=&K(z_{2},z_{-2},z_{1},z_{-1},z_{-3},z_3);\end{align}\begin{align}
\sigma_{321}K(z_{-1},z_1,z_{-2},z_2,z_{-3},z_{3})=&K(z_{-3},z_3,z_{-2},z_2,z_{-1},z_1),\\
\sigma_{-32-1}K(z_{-1},z_1,z_{-2},z_2,z_{-3},z_{3})=&K(z_{3},z_{-3},z_{-2},z_2,z_{1},z_{-1}),\\
\sigma_{3-2-1}K(z_{-1},z_1,z_{-2},z_2,z_{-3},z_{3})=&K(z_{-3},z_3,z_{2},z_{-2},z_{1},z_{-1}),\\
\sigma_{-3-21}K(z_{-1},z_1,z_{-2},z_2,z_{-3},z_{3})=&K(z_{3},z_{-3},z_{2},z_{-2},z_{-1},z_1);\end{align}\begin{align}
\sigma_{312}K(z_{-1},z_1,z_{-2},z_2,z_{-3},z_{3})=&K(z_{-3},z_3,z_{-1},z_1,z_{-2},z_2),\\
\sigma_{-31-2}K(z_{-1},z_1,z_{-2},z_2,z_{-3},z_{3})=&K(z_{3},z_{-3},z_{-1},z_1,z_{2},z_{-2}),\\
\sigma_{3-1-2}K(z_{-1},z_1,z_{-2},z_2,z_{-3},z_{3})=&K(z_{-3},z_3,z_{1},z_{-1},z_{2},z_{-2}),\\
\sigma_{-3-12}K(z_{-1},z_1,z_{-2},z_2,z_{-3},z_{3})=&K(z_{3},z_{-3},z_{1},z_{-1},z_{-2},z_2);\end{align}\begin{align}
\sigma_{231}K(z_{-1},z_1,z_{-2},z_2,z_{-3},z_{3})=
&K(z_{-2},z_2,z_{-3},z_3,z_{-1},z_1),\\
\sigma_{-23-1}K(z_{-1},z_1,z_{-2},z_2,z_{-3},z_{3})=
&K(z_{2},z_{-2},z_{-3},z_3,z_{1},z_{-1}),\\
\sigma_{2-3-1}K(z_{-1},z_1,z_{-2},z_2,z_{-3},z_{3})=
&K(z_{-2},z_2,z_{3},z_{-3},z_{1},z_{-1}),\\
\sigma_{-2-31}K(z_{-1},z_1,z_{-2},z_2,z_{-3},z_{3})=
&K(z_{2},z_{-2},z_{3},z_{-3},z_{-1},z_1);\end{align}\begin{align}
\sigma_{132}K(z_{-1},z_1,z_{-2},z_2,z_{-3},z_{3})=&K(z_{-1},z_1,z_{-3},z_3,z_{-2},z_2),\\
\sigma_{-13-2}K(z_{-1},z_1,z_{-2},z_2,z_{-3},z_{3})=&K(z_{1},z_{-1},z_{-3},z_3,z_{2},z_{-2}),\\
\sigma_{1-3-2}K(z_{-1},z_1,z_{-2},z_2,z_{-3},z_{3})=&K(z_{-1},z_1,z_{3},z_{-3},z_{2},z_{-2}),\\
\sigma_{-1-32}K(z_{-1},z_1,z_{-2},z_2,z_{-3},z_{3})=&K(z_{1},z_{-1},z_{3},z_{-3},z_{-2},z_2).
\end{align}\ees

\noindent{\bf Laplacian}
\beq
\Delta_6=2\partial_{z_{-1}}\partial_{z_{1}}+
2\partial_{z_{-2}}\partial_{z_{2}}+
2\partial_{z_{-3}}\partial_{z_{3}}
.\label{sq2}\eeq



\subsection{$\so(6)$ on the  spherical section}
\label{so(6) on the  spherical section}

In this subsection
we perform Steps \ref{it2}, \ref{it3} and \ref{it4} described in the introduction to this section. Recall that in Step \ref{it2}
we use the null quadric
\[\cV^5:=\{z\in\rr^6\backslash\{0\}\ :\ 2z_{-1}z_1+2z_{-2}z_2+2z_{-3}z_3=0\}.\]
Then, in Step \ref{it3}, we   fix a section of the null quadric. We choose the section  given by the equations
\[4 =2\left(z_{-1}z_1+z_{-2}z_2\right)=-2z_3z_{-3}.\]
We will call it the {\em spherical section}, because it coincides
with $\ss^3(4)\times\ss^1(-4)$. The superscript used for this section
will be ``$\sph$'' for spherical.

In Step \ref{it4} we
introduce the coordinates
\begin{subequations}\label{coo}
\begin{align}
r& =\sqrt{2\left(z_{-1}z_1+z_{-2}z_2\right)}\;,&
w& =\frac{z_{-1}z_1}{z_{-1}z_1+z_{-2}z_{2}}\;,&
\\
u_1& =\frac{z_{1}}{\sqrt{z_{-1}z_{1}+z_{-2}z_2}} \;,&
u_2& =\frac{z_{2}}{\sqrt{z_{-1}z_1+z_{-2}z_{2}}}\;,\\
p& =\sqrt{-2z_3z_{-3}}\;,&
u_3& =\sqrt{-\frac{z_{3}}{z_{-3}}}.
\end{align}\label{subeq}
\end{subequations}
with the inverse transformation
\bes\begin{align}
  z_{-1}&=\frac{rw}{\sqrt2u_1},&z_1&=\frac{u_1r}{\sqrt2},\\
    z_{-2}&=\frac{r(1-w)}{\sqrt2u_2},&z_2&=\frac{u_2r}{\sqrt2},\\
z_{-3}&=-\frac{p}{\sqrt2u_3},&z_3&=\frac{pu_3}{\sqrt2}.  \end{align}\ees

The null quadric in these coordinates is given by
 $r^2=p^2$. We will restrict ourselves to the sheet $r=p$.
The generator of dilations is
\[A_6=r\ddr+p\ddp.\]
The spherical section is given by the condition
$r^2=4$.

All the objects of the previous subsection will be now presented in the above coordinates after the  reduction to the spherical section.
This reduction allows us to eliminate the variables $r,p$.
We omit the superscript $\eta$, whenever there is no dependence on this parameter.
\medskip

\noindent{\bf Lie algebra $\so(6)$.}
\noindent  Cartan operators:
\begin{align*}
N_{1}^\sph& =u_1\ddu{1}\,,\\
N_{2}^\sph& =u_2\ddu{2}\,,\\
N_{3}^\sph& =u_3\ddu{3}\,.
\end{align*}
Roots:

\begin{align*}
  B_{-2,-1}^\sph=u_1 u_2&\p_w,\\
  B_{2,1}^\sph=\frac{1}{u_1u_2}&\left((1-w)w\p_w+(1-w)u_1\p_{u_1}-wu_2\p_{u_2}\right)\, ,\\
B_{2,-1}^\sph=\frac{u_1}{u_2}&\left((1-w)\p_w-u_2\p_{u_2}\right)\, ,\\
B_{-2,1}^\sph=\frac{u_2}{u_1}&\left(w\p_w+u_1\p_{u_1}\right)\, ;
\end{align*}
\begin{align*}  B_{-3,-2}^{\sph,{\eta}}{=}{-}u_2u_3&\left(w\p_w+\frac12\big(u_1\p_{u_1}+u_2\p_{u_2}+u_3\p_{u_3}-\eta\big)\right)\, ,\\ B_{3,2}^{\sph,{\eta}}{=}{-}\frac{1}{u_2u_3}&\left(w(w{-}1)\p_w{+}\frac{(w{-}1)}{2}(u_1\p_{u_1}{+}u_2\p_{u_2}{-}u_3\p_{u_3}{-}\eta){+}u_2\p_{u_2}\right) ,\\
B_{3,-2}^{\sph,{\eta}}=\frac{u_2}{u_3}&\left(w\p_w+\frac12\big(u_1\p_{u_1}+u_2\p_{u_2}-u_3\p_{u_3}-\eta\big)\right)\, ,\\
B_{-3,2}^{\sph,{\eta}}{=}\frac{ u_3}{u_2}&\left(w(w{-}1)\p_w{+}\frac{(w{-}1)}{2}(u_1\p_{u_1}{+}u_2\p_{u_2}{+}u_3\p_{u_3}{-}\eta){+}u_2\p_{u_2}\right);
\end{align*}
\begin{align*}
  B_{-3,-1}^{\sph,{\eta}}{=} {-}u_1u_3&\left((w-1)\p_w+\frac12\big(u_1\p_{u_1}+u_2\p_{u_2}+u_3\p_{u_3}-\eta\big)\right)\, ,\\ B_{3,1}^{\sph,{\eta}}{=}\frac{1}{u_1u_3}&\left(w(w{-}1)\p_w{+}\frac{w}{2}\big(u_1\p_{u_1}{+}u_2\p_{u_2}{-}u_3\p_{u_3}{-}\eta\big){-}u_1\p_{u_1}\right)\, ,\\
B_{3,-1}^{\sph,{\eta}}=\frac{ u_1}{u_3}&\left((w-1)\p_w+\frac12\big(u_1\p_{u_1}+u_2\p_{u_2}-u_3\p_{u_3}-\eta\big)\right)\, ;\\
B_{-3,1}^{\sph,{\eta}}{=}{-}\frac{ u_3}{u_1}&\left(w(w{-}1)\p_w{+}\frac{w}{2}\big(u_1\p_{u_1}{+}u_2\p_{u_2}{+}u_3\p_{u_3}{-}\eta\big)
{-}u_1\p_{u_1}\right).
\end{align*}

\noindent{\bf Weyl symmetries}
  \begin{alignat*}{4}
\sigma_{123}^{\sph,{\eta}} f(w,u_1,u_2,u_3)& =&&f\left(w,u_1,u_2,u_3\right),\\
\sigma_{-12-3}^{\sph,{\eta}} f(w,u_1,u_2,u_3)&
=&&f\left(w,\frac{w}{u_1},u_2,\frac{1}{u_3}\right),\\
\sigma_{1-2-3}^{\sph,{\eta}} f(w,u_1,u_2,u_3)&
=&&f\left(w,u_1,\frac{1-w}{u_2},\frac{1}{u_3}\right),\\
\sigma_{-1-23}^{\sph,{\eta}} f(w,u_1,u_2,u_3)&
=&&f\left(w,\frac{w}{u_1},\frac{1-w}{u_2},u_3\right);\end{alignat*}\begin{alignat*}{4}
    \sigma_{213}^{\sph,{\eta}} f(w,u_1,u_2,u_3)& =&&f\left(1-w,u_2,u_1,u_3\right),\\
    \sigma_{-21-3}^{\sph,{\eta}} f(w,u_1,u_2,u_3)& =&&f\left(1-w,\frac{1-w}{u_2},u_1,\frac1{u_3}\right),\\
    \sigma_{2-1-3}^{\sph,{\eta}} f(w,u_1,u_2,u_3)& =&&f\left(1-w,u_2,\frac{w}{u_1},\frac1{u_3}\right),\\
    \sigma_{-2-13}^{\sph,{\eta}} f(w,u_1,u_2,u_3)& =&&f\left(1-w,\frac{1-w}{u_2},\frac{w}{u_1},u_3\right);\end{alignat*}\begin{alignat*}{4}
    \sigma_{321}^{\sph,{\eta}} f(w,u_1,u_2,u_3)& =&
    \left(\sqrt{-w}\right)^\eta&f\left(\frac{1}{w},\frac{u_3}{\sqrt{-w}}
        ,\frac{u_2}{\sqrt{-w}},\frac{u_1}{\sqrt{-w}}\right),\\
    \sigma_{-32-1}^{\sph,{\eta}} f(w,u_1,u_2,u_3)& =&\left(\sqrt{-w}\right)^\eta&f\left(\frac{1}{w},\frac{1}{\sqrt{-w}u_3}
        ,\frac{u_2}{\sqrt{-w}},\frac{\sqrt{-w}}{u_1}\right),\\
                    \sigma_{3-2-1}^{\sph,{\eta}} f(w,u_1,u_2,u_3)& =&\left(\sqrt{-w}\right)^\eta&f\left(\frac{1}{w},\frac{u_3}{\sqrt{-w}}
            ,\frac{(w-1)}{\sqrt{-w} u_2},\frac{\sqrt{-w}}{u_1}\right),\\
                \sigma_{-3-21}^{\sph,{\eta}} f(w,u_1,u_2,u_3)& =&\left(\sqrt{-w}\right)^\eta&f\left(\frac{1}{w},\frac{1}{\sqrt{-w}u_3}
    ,\frac{(w-1)}{\sqrt{-w}u_2},\frac{u_1}{\sqrt{-w}}\right);\end{alignat*}\begin{alignat*}{4}
    \sigma_{312}^{\sph,{\eta}} f(w,u_1,u_2,u_3)& {=}&\left(\sqrt{w{-}1}\right)^\eta&f\left(\frac{1}{1{-}w},\frac{u_3}{\sqrt{w{-}1}},\frac{u_1}{\sqrt{w{-}1}},\frac{u_2}{\sqrt{w{-}1}}\right),\\
    \sigma_{-31-2}^{\sph,{\eta}} f(w,u_1,u_2,u_3)& {=}&\left(\sqrt{w{-}1}\right)^\eta&f\left(\frac{1}{1{-}w},\frac{1}{\sqrt{w{-}1}u_3},\frac{u_1}{\sqrt{w{-}1}},\frac{\sqrt{w{-}1}}{u_2}\right),\\
    \sigma_{3-1-2}^{\sph,{\eta}} f(w,u_1,u_2,u_3)& {=}&\left(\sqrt{w{-}1}\right)^\eta&f\left(\frac{1}{1{-}w},\frac{u_3}{\sqrt{w{-}1}},\frac{w}{\sqrt{w{-}1}u_1},\frac{\sqrt{w{-}1}}{u_2}\right),\\
    \sigma_{-3-12}^{\sph,{\eta}} f(w,u_1,u_2,u_3)& {=}&\left(\sqrt{w{-}1}\right)^\eta&f\left(\frac{1}{1{-}w},\frac{1}{\sqrt{w{-}1}u_3},\frac{w}{\sqrt{w{-}1}u_1},\frac{u_2}{\sqrt{w{-}1}}\right);\end{alignat*}\begin{alignat*}{4}
    \sigma_{231}^{\sph,{\eta}} f(w,u_1,u_2,u_3)&
=&\left(\sqrt{-w}\right)^\eta&f\left(\frac{w-1}{w},\frac{u_2}{\sqrt{-w}},\frac{u_3}{\sqrt{-w}},\frac{u_1}{\sqrt{-w}}\right),\\
    \sigma_{-23-1}^{\sph,{\eta}} f(w,u_1,u_2,u_3)&
=&\left(\sqrt{-w}\right)^\eta&f\left(\frac{w-1}{w},\frac{(w-1)}{\sqrt{-w}u_2},\frac{u_3}{\sqrt{-w}},\frac{\sqrt{-w}}{u_1}\right),\\
    \sigma_{2-3-1}^{\sph,{\eta}} f(w,u_1,u_2,u_3)&
=&\left(\sqrt{-w}\right)^\eta&f\left(\frac{w-1}{w},\frac{u_2}{\sqrt{-w}},\frac{1}{\sqrt{-w}u_3},\frac{\sqrt{-w}}{u_1}\right),\\
    \sigma_{-2-31}^{\sph,{\eta}} f(w,u_1,u_2,u_3)&
=&\left(\sqrt{-w}\right)^\eta&f\left(\frac{w-1}{w},\frac{(w-1)}{\sqrt{-w}u_2},\frac{1}{\sqrt{-w}u_3},\frac{u_1}{\sqrt{-w}}\right);\end{alignat*}\begin{alignat*}{4}
\sigma_{132}^{\sph,{\eta}} f(w,u_1,u_2,u_3)& =&\left(\sqrt{w{-}1}\right)^\eta&f\left(\frac{w}{w{-}1},\frac{u_1}{\sqrt{w{-}1}},\frac{u_3}{\sqrt{w{-}1}},\frac{u_2}{\sqrt{w{-}1}}\right),\\
\sigma_{-13-2}^{\sph,{\eta}} f(w,u_1,u_2,u_3)& =&\left(\sqrt{w{-}1}\right)^\eta&f\left(\frac{w}{w{-}1},\frac{w}{\sqrt{w{-}1}u_1},\frac{u_3}{\sqrt{w{-}1}},\frac{\sqrt{w{-}1}}{u_2}\right),\\
\sigma_{1{-}3{-}2}^{\sph,{\eta}} f(w,u_1,u_2,u_3)& {=}&\left(\sqrt{w{-}1}\right)^\eta&f\left(\frac{w}{w{-}1},\frac{u_1}{\sqrt{w{-}1}},\frac{1}{\sqrt{w{-}1}u_3},\frac{\sqrt{w{-}1}}{u_2}\right),\\
\sigma_{{-}1{-}32}^{\sph,{\eta}} f(w,u_1,u_2,u_3)& {=}&\left(\sqrt{w{-}1}\right)^\eta&f\left(\frac{w}{w{-}1},\frac{w}{\sqrt{w{-}1}u_1},\frac{1}{\sqrt{w{-}1}u_3},\frac{u_2}{\sqrt{w{-}1}}\right).
\end{alignat*}

\noindent{\bf Laplacian}

\begin{align}
  \Delta_6^\sph=&\quad w(1-w)\p_w^2-\big((1+u_1\p_{u_1})(w-1)+
  (1+u_2\p_{u_2})w\big)\p_w\notag\\
  &-\frac14(u_1\p_{u_1}+u_2\p_{u_2}+1)^2
  +\frac14(u_3\p_{u_3})^2.\label{lapiq2}
\end{align}


Let us give the computations that yield (\ref{lapiq2}). Using
\begin{align*} \p_{z_{-1}}&=\frac{u_1}{\sqrt2r}\Big(-u_1\p_{u_1}-u_2\p_{u_2}+r\p_r+2(1-w)\p_w\Big),\\  \p_{z_1}&=\frac{\sqrt2}{ru_1}\Big(\big(1-\frac{w}{2}\big)u_1\p_{u_1}-\frac{w}{2}u_2\p_{u_2}+\frac{w}{2}r\p_r+w(1-w)\p_w\Big),\\
  \p_{z_{-2}}&=\frac{u_2}{\sqrt2r}\Big(-u_1\p_{u_1}-u_2\p_{u_2}+r\p_r-2w\p_w\Big),\\  \p_{z_2}&=\frac{\sqrt2}{ru_2}\Big(\frac{(w-1)}{2}u_1\p_{u_1}+\frac{(w+1)}{2}u_2\p_{u_2}+\frac{(1-w)}{2}r\p_r+w(w-1)\p_w\Big),\\
  \p_{z_{-3}}&=\frac{u_3}{\sqrt2 p}\Big(u_3\p_{u_3}-p\p_p\Big),\\
    \p_{z_3}&=\frac{1}{\sqrt2 pu_3}\Big(u_3\p_{u_3}+p\p_p\Big),
\end{align*}
we compute the  Laplacian in  coordinates (\ref{coo}):
\begin{align}
\Delta_6&
=\quad \frac1{r^2}\Big(4w(1-w)\p_w^2-4\big((1+u_1\p_{u_1})(w-1)+
  (1+u_2\p_{u_2})w\big)\p_w\notag\\
  &-(u_1\p_{u_1}+u_2\p_{u_2}+1)^2+(r\p_r)^2+2r\p_r+1\Big)\notag\\
  &+\frac{1}{p^2}\Big((u_3\p_{u_3})^2-(p\p_p)^2\Big).
\end{align}
Next we note that
\beq\frac1{r^2}\big((r\p_r)^2+2r\p_r-(p\p_p)^2+1\big)=
\frac1{r^2}\big(r\p_r-p\p_p+1\big)\big(r\p_r+p\p_p+1\big).\label{lapiq3}\eeq
Using $p^2=r^2$ and $r\p_r+p\p_p=-1$, we see that (\ref{lapiq3}) is zero
 on functions of degree $-1$.
Thus we obtain
\begin{align}
\Delta_6^\diamond ={}&\frac{4}{r^2}\Big(
 w(1-w)\p_w^2-\big((1+u_1\p_{u_1})(w-1)+
  (1+u_2\p_{u_2})w\big)\p_w\notag\\
  &-\frac14(u_1\p_{u_1}+u_2\p_{u_2}+1)^2
  +\frac14(u_3\p_{u_3})^2\Big).
\end{align}
To convert
$\Delta_6^\diamond$ into
the   $\Delta_6^\sph$, we simply remove
the prefactor $\frac{4}{r^2}$.

\subsection{Hypergeometric equation}
\label{Hypergeometric equation}

Let us make the ansatz
\beq
f(u_1,u_2,u_3,w)=u_1^\alpha  u_2^\beta u_3^\mu F(w).\label{form0}\eeq
Clearly,
\bes\begin{eqnarray}
N_1^\sph f&=&\alpha f,\\
N_2^\sph f&=&\beta f,\\
N_3^\sph f&=&\mu f,\\
 u_1^{-\alpha}  u_2^{-\beta} u_3^{-\mu}\Delta_6^\sph f&
=& \cF_{\alpha,\beta,\mu}(w,\p_w)F(w),
\end{eqnarray}\ees
where
\begin{align}\notag
 \cF_{\alpha,\beta,\mu}(w,\p_w):=
&w(1-w)\p_w^2-\big((1+\alpha)(w-1)+
(1+\beta)w\big)\p_w\\
&  -\frac14(\alpha+\beta+1)^2
  +\frac14\mu^2,\label{liea}\end{align}
  which is the  ${}_2\cF_1$ {\em hypergeometric operator} in the Lie-algebraic parameters.

  Traditionally, the {\em hypergeometric equation} is given by the operator
\beq
\mathcal{F}(a,b;c;w,\p_w):=w(1-w)\p_w^2+\big(c-(a+b+1)w\big)\p_w-ab,\label{hy1-tra}\eeq
where $a,b,c\in\cc$ will be called the 
{
\em classical parameters}.
Here is the relationship between the Lie-algebraic and classical parameters:
\bes\begin{alignat}{3}
\alpha:=c-1,&\ \ \beta: =a+b-c,&\ \ \mu:=a-b;\\[2ex]
\label{newnot}
a=\frac{1+\alpha+\beta +\mu}{2},&\ \ b=\frac{1+\alpha+\beta -\mu}{2},&\ \ c=1+\alpha.
\end{alignat}
\ees

Note that the Lie-algebraic parameters $\alpha,\beta,\mu$ are differences of the indices of the singular points $0,1,\infty$. For many purposes, they are more convenient than the traditional parameters $a,b,c$. They are used e.g. in Subsect. 2.7.2 of \cite{BE}, where they are called $\lambda,\nu,\mu$. In the standard notation for Jacobi Polynomials $P_n^{\alpha,\beta}$, the parameters $\alpha,\beta$ correspond to our $\alpha, \beta$ (where the singular points have been moved from $0,1$ to $-1,1$).

  \subsection{Transmutation relations and discrete symmetries}
\label{Transmutation relations and discrete symmetries}

By \eqref{dada1}, we have the following generalized symmetries
\begin{subequations}
\begin{align}
  B^{\sph,-3}\Delta_6^\sph&=\Delta_6^\sph B^{\sph,-1},\quad B\in \so(6),\label{siy1}\\
  \alpha^{\sph,-3}\Delta_6^\sph&=\Delta_6^\sph \alpha^{\sph,-1},\quad \alpha\in \mathrm{O}(6).\label{siy2}
\end{align}
\end{subequations}
Applying (\ref{siy1}) to the roots of $\so(6)$ we obtain the {\em transmutation relations for the hypergeometric operator}:

\noindent
 \begin{align*}
 \p_w&\cF_{\alpha,\beta ,\mu }\notag\\
= \cF_{\alpha+1,\beta +1,\mu } &\p_w,\\
\big(w(1-w)\p_w+(1-w)\alpha-w\beta \big)&\cF_{\alpha,\beta ,\mu }\notag\\
= \cF_{\alpha-1,\beta -1,\mu }&\big(w(1-w)\p_w+(1-w)\alpha-w\beta \big),\\
\big((1-w)\p_w -\beta \big)&\cF_{\alpha,\beta ,\mu }\notag\\
= \cF_{\alpha+1,\beta -1,\mu }& \big((1-w)\p_w -\beta \big),\\
   (w\p_w+\alpha)&\cF_{\alpha,\beta ,\mu }\notag
\\
= \cF_{\alpha-1,\beta +1,\mu }&(w\p_w+\alpha);
 \end{align*}\begin{align*}
\Big(w\p_w+\12 (\alpha+ \beta +\mu +1)\Big)&w\cF_{\alpha,\beta ,\mu }\notag
\\
= w\cF_{\alpha,\beta +1,\mu +1}&\Big(w\p_w+\12 (\alpha+ \beta +\mu +1)\Big),\\
\Big(w(w{-}1)\p_w{+}\frac12(w{-}1)(\alpha{+}\beta {-}\mu {+}1){-}\beta \Big)&w\cF_{\alpha,\beta ,\mu }\notag\\
= w\cF_{\alpha,\beta {-}1,\mu {-}1}\qquad&\hspace{-4ex}
\Big(w(w{-}1)\p_w{+}\frac12(w{-}1)(\alpha{+}\beta {-}\mu {+}1){-}\beta \Big),\\
\Big( w\p_w{+}\12(\alpha{+}\beta {-}\mu {+}1)\Big)&w\cF_{\alpha,\beta ,\mu }\notag\\
= w\cF_{\alpha,\beta {+}1,\mu {-}1}&\Big( w\p_w{+}\12(\alpha{+}\beta {-}\mu {+}1\Big),
\\
\Big(w(w{-}1)\p_w{-}\12(1{-}w)(\alpha{+}\beta {+}\mu {+}1)
{+}\beta \Big)&w\cF_{\alpha,\beta ,\mu }\notag\\
= w\cF_{\alpha,\beta -1,\mu +1}\qquad&\hspace{-4ex}\Big(w(w{-}1)\p_w{-}\12(1{-}w)(\alpha{+}\beta {+}\mu {+}1)
{+}\beta \Big);\end{align*}\begin{align*}
\Big((w-1)\p_w+\12(\alpha+\beta +\mu +1)\Big)&(1-w)\cF_{\alpha,\beta ,\mu }\notag\\
= (1-w)\cF_{\alpha+1,\beta ,\mu +1}& \Big((w-1)\p_w+\12(\alpha+\beta +\mu +1\Big)
,
\\
\Big(w(w{-}1)\p_w{+}\12w(\alpha{+}\beta {-}\mu {+}1)
{+}\alpha\Big)&(1-w)\cF_{\alpha,\beta ,\mu }\notag
\\
= (1-w)\cF_{\alpha-1,\beta ,\mu -1}&\Big(w(w{-}1)\p_w{+}\12w(\alpha{+}\beta {-}\mu {+}1)
{+}\alpha\Big),
\\
\Big((w-1)\p_w+\12(\alpha+\beta -\mu +1)\Big)&
(1-w)\cF_{\alpha,\beta ,\mu }\notag
\\= 
(1-w)\cF_{\alpha+1,\beta ,\mu -1}&\Big((w-1)\p_w+\12(\alpha+\beta -\mu +1)\Big),
\\
\Big(w(w{-}1)\p_w{+}\12w(\alpha{+}\beta {+}\mu {+}1)
-\alpha\Big)&(1-w)\cF_{\alpha,\beta ,\mu }\notag
\\
= (1-w)\cF_{\alpha-1,\beta ,\mu +1}&\Big(w(w{-}1)\p_w{+}\12w(\alpha{+}\beta {+}\mu {+}1)
{-}\alpha\Big).
\end{align*}

Applying (\ref{siy2}) to the Weyl group $D_3$ we obtain the {\em discrete symmetries of the hypergeometric operator}. We describe them below, restricting ourselves to $D_3\cap \SO(6)$.
  
All the operators below equal $\cF_{\alpha,\beta ,\mu}(w,\p_w)$ for the
corresponding $w$:
\noindent\[\begin{array}{rrcl}
w=v:&&
\cF_{\alpha,\beta ,\mu}(v,\p_v),&
\\[0ex]
&(-v)^{-\alpha}(v-1)^{-\beta }&\cF_{-\alpha,-\beta ,\mu}(v,\p_v)&(-v)^{\alpha}(v-1)^{\beta }
\\[0ex]
&(v-1)^{-\beta }&\cF_{\alpha,-\beta ,-\mu}(v,\p_v)&(v-1)^{\beta },\\[0ex]
&(-v)^{-\alpha}&\cF_{-\alpha,\beta ,-\mu}(v,\p_v)&
(-v)^{\alpha};\end{array}\]\[\begin{array}{rrcl}
w=1-v:&&\cF_{\beta ,\alpha,\mu}(v,\p_v),&
\\[0ex]
&(v-1)^{-\alpha}(-v)^{-\beta }&\cF_{-\beta ,-\alpha,\mu}(v,\p_v)&(v-1)^{\alpha}(-v)^{\beta },
\\[0ex]
&(v-1)^{-\alpha}&\cF_{\beta ,-\alpha,-\mu}(v,\p_v)&
(v-1)^{\alpha},\\[0ex]
&(-v)^{-\beta }&\cF_{-\beta ,\alpha,-\mu}(v,\p_v)
&(-v)^{\beta };\end{array}\]\[\begin{array}{rrcl}
&w=\frac{1}{v}:\quad(-v)^{\frac{\alpha+\beta +\mu +1}{2}}&(-v)\cF_{\mu ,\beta ,\alpha}(v,\p_v)&
(-v)^{\frac{-\alpha-\beta -\mu -1}{2}},\\[0ex]
&(-v)^{\frac{\alpha+\beta -\mu +1}{2}}(v-1)^{-\beta }&
(-v)\cF_{-\mu ,-\beta ,\alpha}(v,\p_v) &(-v)^{\frac{-\alpha-\beta +\mu -1}{2}}(v-1)^{\beta },
\\[0ex]
&(-v)^{\frac{\alpha+\beta +\mu +1}{2}}(v-1)^{-\beta }&
(-v)
\cF_{\mu ,-\beta ,-\alpha}(v,\p_v)& (-v)^{\frac{-\alpha-\beta -\mu -1}{2}}(v-1)^{\beta }
,\\[0ex]
&(-v)^{\frac{\alpha+\beta -\mu +1}{2}}&(-v)\cF_{-\mu ,\beta ,-\alpha}(v,\p_v)
& (-v)^{\frac{-\alpha-\beta +\mu -1}{2}}
;\end{array}\]\[\begin{array}{rrcl}
&w=\frac{v-1}{v}:\quad(-v)^{\frac{\alpha+\beta +\mu +1}{2}}&
(-v)\cF_{\mu ,\alpha,\beta }(v,\p_v)&(-v)^{\frac{-\alpha-\beta -\mu -1}{2}},\\[0ex]
&(-v)^{\frac{\alpha+\beta -\mu +1}{2}}(v-1)^{-\alpha}&
(-v)\cF_{-\mu ,-\alpha,\beta }(v,\p_v) &(-v)^{\frac{-\alpha-\beta +\mu -1}{2}}(v-1)^{\alpha},\\[0ex]
&(-v)^{\frac{\alpha+\beta +\mu +1}{2}}(v-1)^{-\alpha}&
(-v) \cF_{\mu ,-\alpha,-\beta }(v,\p_v)& (-v)^{\frac{-\alpha-\beta -\mu -1}{2}}(v-1)^{\alpha},
\\[0ex]
 &(-v)^{\frac{\alpha+\beta -\mu +1}{2}}&(-v)\cF_{-\mu ,\alpha,-\beta }(v,\p_v)&
  (-v)^{\frac{-\alpha-\beta +\mu -1}{2}};\end{array}\]\[\begin{array}{rrcl}
&w{=}\frac{1}{1{-}v}:\,
 (v{-}1)^{\frac{\alpha+\beta +\mu +1}{2}}&
(v{-}1)\cF_{\beta ,\mu ,\alpha}(v,\p_v)& (v{-}1)^{\frac{-\alpha-\beta -\mu -1}{2}},\\[0ex]
&  (-v)^{-\beta }(v{-}1)^{\frac{\alpha+\beta -\mu +1}{2}}
&(v{-}1)\cF_{-\beta ,-\mu ,\alpha}(v,\p_v) & (-v)^{\beta }(v{-}1)^{\frac{-\alpha-\beta +\mu -1}{2}},\\[0ex]
& (v{-}1)^{\frac{\alpha+\beta -\mu +1}{2}}
&(v{-}1)\cF_{\beta ,-\mu ,-\alpha}(v,\p_v) & (v{-}1)^{\frac{-\alpha-\beta +\mu -1}{2}},\\[0ex]
&(-v)^{-\beta }(v{-}1)^{\frac{\alpha+\beta +\mu +1}{2}}&(v{-}1)
\cF_{-\beta ,\mu ,-\alpha}(v,\p_v)&
 (-v)^{\beta }(v{-}1)^{\frac{-\alpha-\beta -\mu -1}{2}};
\end{array}\]\[\begin{array}{rrcl}
&w{=}\frac{v}{v{-}1}:\,
(v{-}1)^{\frac{\alpha+\beta +\mu+1}{2}}&(v{-}1)\cF_{\alpha,\mu,\beta }(v,\p_v)&
(v{-}1)^{\frac{-\alpha-\beta -\mu-1}{2}}
,\\[0ex]
&  (-v)^{-\alpha}(v{-}1)^{\frac{\alpha+\beta -\mu+1}{2}}&(v{-}1)\cF_{-\alpha,-\mu,\beta }(v,\p_v) &
 (-v)^{\alpha}(v{-}1)^{\frac{-\alpha-\beta +\mu-1}{2}}
,\\[0ex]
& (v{-}1)^{\frac{\alpha+\beta -\mu+1}{2}}&
(v{-}1)\cF_{\alpha,-\mu,-\beta }(v,\p_v) & (v{-}1)^{\frac{-\alpha-\beta +\mu-1}{2}}
,\\[0ex]
& (-v)^{-\alpha}(v{-}1)^{\frac{\alpha+\beta +\mu +1}{2}}
&(v{-}1)\cF_{-\alpha,\mu ,-\beta }(v,\p_v)&
  (-v)^{\alpha}(v{-}1)^{\frac{-\alpha-\beta -\mu -1}{2}}.
\end{array}\]

\subsection{Factorizations of the Laplacian}
\label{Factorizations of the Laplacian}

In the Lie algebra $\so(6)$ represented on $\rr^6$
we have 3 distinguished Lie subalgebras isomorphic to $\so(4)$:
\beq\so_{12}(4),\ \ \so_{23}(4),\ \ \so_{13}(4),\eeq
where we use a hopefully obvious notation.
 By (\ref{casimir2}), the corresponding Casimir operators are
\begin{subequations}
\begin{eqnarray}
\cC_{12}&=&
4B_{1,2}B_{-1,-2}-(N_1+N_2+1)^2+1\\
&=&4B_{-1,-2}B_{1,2}-(N_1+N_2-1)^2+1\\
&=&4B_{1,-2}B_{-1,2}-(N_1-N_2+1)^2+1\\
&=&4B_{-1,2}B_{1,-2}-(N_1-N_2-1)^2+1;\\
\cC_{23}&=&
4B_{2,3}B_{-2,-3}-(N_2+N_3+1)^2+1\\
&=&4B_{-2,-3}B_{2,3}-(N_2+N_3-1)^2+1\\
&=&4B_{2,-3}B_{-2,3}-(N_2-N_3+1)^2+1\\
&=&4B_{-2,3}B_{2,-3}-(N_2-N_3-1)^2+1;\\
\cC_{13}&=&
4B_{1,3}B_{-1,-3}-(N_1+N_3+1)^2+1\\
&=&4B_{-1,-3}B_{1,3}-(N_1+N_3-1)^2+1\\
&=&4B_{1,-3}B_{-1,3}-(N_1-N_3+1)^2+1\\
&=&4B_{-1,3}B_{1,-3}-(N_1-N_3-1)^2+1.
\end{eqnarray}\label{subu-}
\end{subequations}
Of course, for any $\eta$ we can append the superscript ${}^{\diamond,\eta}$ to all the operators in (\ref{subu-}).

After the reduction described in (\ref{deq1a}), we obtain the
identities
\begin{subequations}\label{facto}
\begin{eqnarray}\label{facto1--}
(2z_{-1}z_1+2z_{-2}z_2)\Delta_6^\diamond=-1+\cC_{12}^{\diamond,-1}+(N_3^{\diamond,-1})^2,\\
\label{facto3--}
(2z_{-2}z_2+2z_{-3}z_3)\Delta_6^\diamond=-1+\cC_{23}^{\diamond,-1}+(N_1^{\diamond,-1})^2,\\
\label{facto2--}(2z_{-1}z_1+2z_{-3}z_3)\Delta_6^\diamond=-1+\cC_{13}^{\diamond,-1}+(N_2^{\diamond,-1})^2.
\end{eqnarray}
\end{subequations}

We insert (\ref{subu-}) with superscript ${}^{\diamond,-1}$ to 
(\ref{facto}), obtaining
\begin{subequations}
\begin{align}\notag
&(2z_{-1}z_1+2z_{-2}z_2)\Delta_6^\diamond\\=&
4B_{1,2}B_{-1,-2}-(N_1+N_2+N_3+1)(N_1+N_2-N_3+1)\\
=&4B_{-1,-2}B_{1,2}-(N_1+N_2+N_3-1)(N_1+N_2-N_3-1)\\
=&4B_{1,-2}B_{-1,2}-(N_1-N_2+N_3+1)(N_1-N_2-N_3+1)\\
=&4B_{-1,2}B_{1,-2}-(N_1-N_2+N_3-1)(N_1-N_2-N_3-1);\\
&(2z_{-2}z_2+2z_{-3}z_3)\Delta_6^\diamond\notag\\
=&
4B_{2,3}B_{-2,-3}-(N_1+N_2+N_3+1)(-N_1+N_2+N_3+1)\\
=&4B_{-2,-3}B_{2,3}-(N_1+N_2+N_3-1)(-N_1+N_2+N_3-1)\\
=&4B_{2,-3}B_{-2,3}-(N_1+N_2-N_3+1)(-N_1+N_2-N_3+1)\\
=&4B_{-2,3}B_{2,-3}-(N_1+N_2-N_3-1)(-N_1+N_2-N_3-1);\\\notag
&(2z_{-1}z_1+2z_{-3}z_3)\Delta_6^\diamond\\=&
4B_{1,3}B_{-1,-3}-(N_1+N_2+N_3+1)(N_1-N_2+N_3+1)\\
=&4B_{-1,-3}B_{1,3}-(N_1+N_2+N_3-1)(N_1-N_2+N_3-1)\\
=&4B_{1,-3}B_{-1,3}-(N_1+N_2-N_3+1)(N_1-N_2-N_3+1)\\
=&4B_{-1,3}B_{1,-3}-(N_1+N_2-N_3-1)(N_1-N_2-N_3-1);
\end{align}\label{subu-.}
\end{subequations}
where  for typographical reasons we omitted
 the superscript ${}^{\diamond,{-}1}$ at all the operators $B$ and $N$.


If we use the coordinates \eqref{coo} and the spherical section, then we have to rewrite (\ref{subu-.})
by making the replacements
\begin{subequations}
\begin{align}\label{facto1b}
2z_{-1}z_1+2z_{-2}z_2&\to\quad1,\\
  \label{facto3b}
2z_{-2}z_2+2z_{-3}z_3&\to\quad-w,\\
\label{facto2b}
2z_{-1}z_1+2z_{-3}z_3&\to\quad w-1,
\end{align}
\end{subequations}
as well as replacing the superscript ${}^{\diamond}$ with ${}^{\sph}$.

\subsection{Factorizations of the hypergeometric operator}
\label{Factorizations of the hypergeometric operator}

The factorizations of $\Delta_6^\sph$ described in Subsect.
\ref{Factorizations of the Laplacian} yield
the following factorizations of the hypergeometric operator:
\begin{align*}
&\cF_{\alpha,\beta,\mu}\\=&
\Big(w(1-w)\partial_w+\big((1+\alpha)(1-w)-(1+\beta)  w\big)\Big)\partial_w\\&-\frac14(\alpha+\beta+\mu+1)(\alpha+\beta-\mu+1)\\
=&
\partial_w\Big(w(1-w)\partial_w+\big(\alpha(1-w)-\beta  w\big)\Big)\\&-\frac14(\alpha+\beta+\mu-1)(\alpha+\beta-\mu-1)\\
=&
\Big(w\partial_w+\alpha+1\Big)\Big((1-w)\partial_w-\beta\Big)\\
&-\frac14(\alpha-\beta+\mu+1)(\alpha-\beta-\mu+1)\\
=&
\Big((1-w)\partial_w-\beta-1\Big)\Big(w\partial_w+\alpha\Big)\\
&-\frac14(\alpha-\beta+\mu-1)(\alpha-\beta-\mu-1);
\end{align*}\begin{align*}
&w\cF_{\alpha,\beta,\mu}\\=&
\Big(w\partial_w{+}\frac12(\alpha{+}\beta{+}\mu{-}1)\Big)
\Big(w(1{-}w)\partial_w{+}\frac12(1{-}w)(\alpha{+}\beta{-}\mu{+}1){-}\beta\Big)\\
&{-}\frac14(\alpha{+}\beta{+}\mu{-}1)(\alpha{-}\beta{-}\mu{+}1)\\
=&
\Big(w(1{-}w)\partial_w{+}\frac12(1{-}w)(\alpha{+}\beta{-}\mu{+}1){-}\beta{-}1\Big)
\Big(w\partial_w{+}\frac12(\alpha{+}\beta{+}\mu{+}1)\Big)
\\
&{-}\frac14(\alpha{+}\beta{+}\mu{+}1)(\alpha{-}\beta{-}\mu{-}1)\\
=&
\Big(w\partial_w{+}\frac12(\alpha{+}\beta{-}\mu{-}1)\Big)
\Big(w(1{-}w)\partial_w{+}\frac12(1{-}w)(\alpha{+}\beta{+}\mu{+}1){-}\beta\Big)\\
&{-}\frac14(\alpha{+}\beta{-}\mu{-}1)(\alpha{-}\beta{+}\mu{+}1)\\
=&
\Big(w(1{-}w)\partial_w{+}\frac12(1{-}w)(\alpha{+}\beta{+}\mu{+}1){-}\beta{-}1\Big)
\Big(w\partial_w{+}\frac12(\alpha{+}\beta{-}\mu{+}1)\Big)
\\
&{-}\frac14(\alpha{+}\beta{-}\mu{+}1)(\alpha{-}\beta{+}\mu{-}1);
\end{align*}\begin{align*}
&(w{-}1)\cF_{\alpha,\beta,\mu}\\
=&
\Big(w(w{-}1)\partial_w{+}\frac12w(\alpha{+}\beta{-}\mu{+}1){-}\alpha{-}1\Big)
\Big((w{-}1)\partial_w{+}\frac12(\alpha{+}\beta{+}\mu{+}1)\Big)
\\
&{-}\frac14(\alpha{+}\beta{+}\mu{+}1)(\alpha{-}\beta{+}\mu{+}1)\\
=&
\Big((w{-}1)\partial_w{+}\frac12(\alpha{+}\beta{+}\mu{-}1)\Big)
\Big(w(w{-}1)\partial_w{+}\frac12w(\alpha{+}\beta{-}\mu{+}1){-}\alpha\Big)\\
&{-}\frac14(\alpha{+}\beta{+}\mu{-}1)(\alpha{-}\beta{+}\mu{-}1)\\
=&
\Big(w(w{-}1)\partial_w{+}\frac12w(\alpha{+}\beta{+}\mu{+}1){-}\alpha{-}1\Big)
\Big((w{-}1)\partial_w{+}\frac12(\alpha{+}\beta{-}\mu{+}1)\Big)
\\
&{-}\frac14(\alpha{+}\beta{-}\mu{+}1)(\alpha{-}\beta{-}\mu{+}1)\\
=&
\Big((w{-}1)\partial_w{+}\frac12(\alpha{+}\beta{-}\mu{-}1)\Big)
\Big(w(w{-}1)\partial_w{+}\frac12w(\alpha{+}\beta{+}\mu{+}1){-}\alpha\Big)\\
&{-}\frac14(\alpha{+}\beta{-}\mu{-}1)(\alpha{-}\beta{-}\mu{-}1).
\end{align*}

\subsection{The ${}_2F_1$ hypergeometric function}
  \label{The ${}_2F_1$ hypergeometric function}

$0$ is a regular singular point of 
the ${}_2\cF_1$ hypergeometric equation.  Its indices are 
 $0$ and $1-c$. 
For $c\neq 0,-1,-2,\dots$ the Frobenius method yields  the unique solution of
the hypergeometric equation equal to $1$
 at $0$, given by the series 
\[F(a,b;c;w)=\sum_{j=0}^\infty
\frac{(a)_j(b)_j}{
(c)_j}\frac{w^j}{j!}\]
convergent for $|w|<1$. The function extends to the whole complex plane cut at $[1,\infty[$ and is 
called the {\em hypergeometric function}.
Sometimes it is more convenient to consider the function
\[ {\bf F}  (a,b;c;w):=\frac{F(a,b,c,w)}{\Gamma(c)}
=\sum_{j=0}^\infty
\frac{(a)_j(b)_j}{
\Gamma(c+j)}\frac{w^j}{j!}\]
 defined for all $a,b,c\in\cc$.
Another useful function proportional to $F$ is
\[ {\bf F}^\I  (a,b;c;w):=\frac{\Gamma(b)\Gamma(c-b)}{\Gamma(c)}
F(a,b;c;w)
=\sum_{j=0}^\infty
\frac{\Gamma(b+j)\Gamma(c-b)(a)_j}{
\Gamma(c+j)}\frac{w^j}{j!}.
\]

We will usually prefer to parametrize all varieties of the
hypergeometric function with the Lie-algebraic parameters:
\begin{eqnarray*}
 F_{\alpha,\beta ,\mu }(w)&=&F\Bigl(
\frac{1+\alpha+\beta +\mu}{2},\frac{1+\alpha+\beta -\mu}{2};1+\alpha;w\Bigr),\\
 {\bf F}_{\alpha,\beta ,\mu }(w)&=&{\bf F} \Bigl(
\frac{1+\alpha+\beta +\mu}{2},\frac{1+\alpha+\beta -\mu}{2};1+\alpha;w\Bigr)\\
&=&
\frac{1}{\Gamma(\alpha+1)} F_{\alpha,\beta ,\mu }(w),\\
 {\bf F}_{\alpha,\beta ,\mu }^{\I}(w)&=&{\bf F}^{\I}\Bigl(
\frac{1+\alpha+\beta +\mu}{2},\frac{1+\alpha+\beta -\mu}{2};1+\alpha;w\Bigr)\\
&=&
\frac{\Gamma\big(\frac{1+\alpha+\beta-\mu}{2}\big)\Gamma\big(\frac{1+\alpha-\beta+\mu}{2}\big)}{\Gamma(\alpha+1)}
F_{\alpha,\beta ,\mu }(w).
\end{eqnarray*}

\subsection{Standard solutions}

The hypergeometric equation has 3 singular points. With each of them we can associate two solutions with a simple behavior. Therefore, we obtain 6 standard solutions.

Applying the discrete symmetries from $D_3\cap \SO(6)$ to  the hypergeometric function, we obtain 24  expressions for solutions of the hypergeometric equation, which go under the name of {\em Kummer's table}. Some of them coincide as functions, so that we obtain 6 standard solutions, each expressed in 4 ways:

\begin{align*}
\text{Solution $\sim1$ at $0$:}\quad& F_{\alpha,\beta ,\mu }(w)\\
=&(1-w)^{-\beta } F_{\alpha,-\beta ,-\mu }(w)\\
=&(1-w)^{\frac{-1-\alpha-\beta +\mu }{2}}
 F_{\alpha,-\mu ,-\beta }(\frac{w}{w-1})\\
 =&(1-w)^{\frac{-1-\alpha-\beta -\mu }{2}} F_{\alpha,\mu ,\beta }(\frac{w}{w-1});
 \end{align*}\begin{align*}
\text{Solution $\sim w^{-\alpha}$ at $0$:}\quad&w^{-\alpha} F_{-\alpha,\beta ,-\mu }(w)\\
=&w^{-\alpha}(1-w)^{-\beta } F_{-\alpha,-\beta ,\mu }(w)\\
=&w^{-\alpha}(1-w)^{\frac{-1+\alpha-\beta +\mu }{2}}
 F_{-\alpha,-\mu ,\beta }(\frac{w}{w-1})\\
 =&w^{-\alpha}(1-w)^{\frac{-1+\alpha-\beta -\mu }{2}} F_{-\alpha,\mu ,-\beta }(\frac{w}{w-1});
  \end{align*}\begin{align*}
  \text{Solution $\sim1$  at $1$:}\quad
& F_{\beta ,\alpha,\mu }(1-w)\\
=&w^{-\alpha} F_{\beta ,-\alpha,-\mu }(1-w)\\
=&w^{\frac{-1-\alpha-\beta +\mu }{2}}
 F_{\beta ,-\mu ,-\alpha}(1-w^{-1})\\
=&w^{\frac{-1-\alpha-\beta -\mu }{2}}
F_{\beta ,\mu ,\alpha}(1-w^{-1});
 \end{align*}\begin{align*}
\text{Solution $\sim(1-w)^{-\beta }$ at $1$:}\quad&(1-w)^{-\beta } F_{-\beta ,\alpha,-\mu }(1-w)\\
=&w^{-\alpha}(1-w)^{-\beta } F_{-\beta ,-\alpha,\mu }(1-w)\\
=&w^{\frac{-1-\alpha+\beta -\mu }{2}}(1-w)^{-\beta }
 F_{-\beta ,\mu ,-\alpha}(1-w^{-1})\\
=&w^{\frac{-1-\alpha+\beta +\mu }{2}}(1-w)^{-\beta }
F_{-\beta ,-\mu ,\alpha}(1-w^{-1});
 \end{align*}\begin{align*}
\text{Solution $\sim w^{-a}$ at $\infty$:}\quad&(-w)^{\frac{-1-\alpha-\beta -\mu }{2}} F_{\mu ,\beta ,\alpha}(w^{-1})\\
=&(-w)^{\frac{-1-\alpha+\beta -\mu }{2}}(1-w)^{-\beta } F_{\mu ,-\beta ,-\alpha}(w^{-1})\\
=&(1-w)^{\frac{-1-\alpha-\beta -\mu }{2}} F_{\mu ,\alpha,\beta }((1-w)^{-1})\\
=&(-w)^{-\alpha}(1-w)^{\frac{-1+\alpha-\beta -\mu }{2}} F_{\mu ,-\alpha,-\beta }((1-w)^{-1});
 \end{align*}\begin{align*}
  \text{Solution $\sim w^{-b}$ at $\infty$:}\quad
  &(-w)^{\frac{-1-\alpha-\beta +\mu }{2}} F_{-\mu ,\beta ,-\alpha}(w^{-1})\\
=&(-w)^{\frac{-1-\alpha+\beta +\mu }{2}}(1-w)^{-\beta } F_{-\mu ,-\beta ,\alpha}(w^{-1})\\
=&(1-w)^{\frac{-1-\alpha-\beta +\mu }{2}} F_{-\mu ,\alpha,-\beta }((1-w)^{-1})\\
=&(-w)^{-\alpha}(1-w)^{\frac{-1+\alpha-\beta +\mu }{2}} F_{-\mu ,-\alpha,\beta }((1-w)^{-1}).
\end{align*}

\subsection{Recurrence relations}
\label{s3.15}
To each root of $\so(6)$ there corresponds a recurrence relation:
\begin{align*}
\p_w {\bf F}^\I  _{\alpha,\beta ,\mu }(w)&=
\frac{1{+}\alpha{+}\beta {+}\mu }{2} {\bf F}^\I  _{\alpha+1,\beta +1,\mu }(w),\\[0ex]
-(w(1{-}w)\p_w{+}\alpha(1{-}w){-}\beta w) {\bf F}^\I  _{\alpha,\beta ,\mu }(w)&=
 \frac{1{-}\alpha{-}\beta{+}\mu}{2}{\bf F}^\I  _{\alpha-1,\beta -1,\mu }(w),\\[0ex]
((1-w)\p_w-\beta ) {\bf F}^\I  _{\alpha,\beta ,\mu }(w)&=\frac{1{+}\alpha{-}\beta {-}\mu }{2}{\bf F}^\I  _{\alpha{+}1,\beta {-}1,\mu }(w),\\[0ex]
- (w\p_w+\alpha) {\bf F}^\I  _{\alpha,\beta ,\mu }(w)&= \frac{1{-}\alpha{+}\beta{-}\mu}{2}{\bf F}^\I  _{\alpha-1,\beta +1,\mu }(w);
\end{align*}
\begin{align*}
 \left(w\p_w+\frac{1+\alpha+\beta +\mu }{2}\right) {\bf F}^\I  _{\alpha,\beta ,\mu }(w)&=
\frac{1{+}\alpha{+}\beta {+}\mu }{2} {\bf F}^\I  _{\alpha,\beta +1,\mu +1}(w)
, \\[0ex]
{-}\! \left(\!w(w{-}1)\p_w{+}\beta {+}\frac{1{+}\alpha{+}\beta {-}\mu }{2}(w{-}1)\!
\right)\! {\bf F}^\I  _{\alpha,\beta ,\mu }(w)&=\frac{1{+}\alpha{-}\beta {-}\mu }{2} {\bf F}^\I  _{\alpha,\beta -1,\mu -1}(w),\\[0ex]
- \left(w\p_w+\frac{1+\alpha+\beta -\mu }{2}\right) {\bf F}^\I  _{\alpha,\beta ,\mu }(w)&=
\frac{1{-}\alpha{+}\beta{-}\mu }{2} {\bf F}^\I  _{\alpha,\beta +1,\mu -1}(w)
, \\[0ex]
 \left(w(w{-}1)\p_w{+}\beta {+}\frac{1{+}\alpha{+}\beta {+}\mu }{2}(w{-}1)\right) {\bf F}^\I  _{\alpha,\beta ,\mu }(w)&=
\frac{1{-}\alpha{-}\beta{+}\mu }{2}
 {\bf F}^\I  _{\alpha,\beta -1,\mu +1}(w)
 ;\end{align*}
\begin{align*}
\left((w-1)\p_w+\frac{1+\alpha+\beta +\mu }{2}\right) {\bf F}^\I  _{\alpha,\beta ,\mu }(w)&=
\frac{1{+}\alpha{+}\beta {+}\mu }{2}
{\bf F}^\I  _{\alpha{+}1,\beta ,\mu {+}1}(w)
, \\[0ex]
 \left(w(w-1)\p_w{-}\alpha{+}\frac{1{+}\alpha{+}\beta {-}\mu }{2}w\right)
 {\bf F}^\I  _{\alpha,\beta ,\mu }(w)&=\frac{1{-}\alpha{+}\beta{-}\mu}{2} {\bf F}^\I  _{\alpha-1,\beta ,\mu -1}(w),
\\[0ex]
 \left((w-1)\p_w+\frac{1+\alpha+\beta -\mu }{2}\right)
 {\bf F}^\I  _{\alpha,\beta ,\mu }(w)&=\frac{1{+}\alpha{-}\beta {-}\mu }{2}
 {\bf F}^\I  _{\alpha{+}1,\beta ,\mu {-}1}(w)
, \\[0ex]
 \left(w(w-1)\p_w{-}\alpha{+}\frac{1{+}\alpha{+}\beta {+}\mu }{2}w\right) {\bf F}^\I  _{\alpha,\beta ,\mu }(w)&=\frac{1{-}\alpha{-}\beta{+}\mu}{2} {\bf F}^\I  _{\alpha-1,\beta ,\mu +1}(w)
 .
 \end{align*}
The recurrence relations are essentially fixed by the transmutation relations. The only missing piece of information is the coefficient on the right hand side, which can be derived by analyzing the behavior of both sides around zero. Another way to obtain these coefficients is to use the integral representations described in the following subsections.

\subsection{Wave packets in 6 dimensions}
\label{Wave packets in 6 dimensions}

We start with the following easy fact:

\bel For any $\tau$, the following function is harmonic on $\rr^6$:
\beq
(z_1-\tau^{-1}z_{-2})^{\alpha +\nu}(z_2+\tau^{-1}z_{-1})^{\beta +\nu}z_3^\mu
\label{cacy}\eeq
\label{cacko}\eel

\proof
Set
$e_1:=(1,0,0,-\tau^{-1})$, $e_2:=(0,\tau^{-1},1,0)$. 
Then
\[\langle e_1|e_1\rangle=\langle e_2|e_2\rangle=\langle e_2|e_1\rangle=0.\]
Hence, (\ref{cacy}) is harmonic by Prop. \ref{isotro}. \qed

Let us make a wave packet out of
(\ref{cacy}), which is an eigenfunction of the Cartan operators:
\begin{align}\notag
&K_{\alpha,\beta ,\mu,\nu}(z_{-1},z_1,z_{-2},z_2,z_{-3},z_{3})\\
:=&\int_\gamma(z_1-\tau^{-1}z_{-2})^{\alpha+\nu}(z_2+\tau^{-1}z_{-1})^{\beta +\nu}
z_3^{\mu}\tau^{\nu-1}\frac{\d\tau}{2\pi\i}.\label{wave}
\end{align}

\bep \label{caca7}
Let the contour $]0,1[\ni s\overset{\gamma}\mapsto\tau(s)$
    satisfy
\begin{align}    (z_1-\tau^{-1}z_{-2})^{\alpha+\nu}(z_2+\tau^{-1}z_{-1})^{\beta +\nu}
\tau^{\nu-1}\Big|_{\tau(0)}^{\tau(1)}&=0.\label{condi2}
    \end{align}
Then  $ K_{\alpha,\beta ,\mu,\nu}$ is harmonic and
\bes\begin{align}
 \label{cac1.}
      N_1K_{\alpha,\beta ,\mu,\nu}&=\alpha K_{\alpha,\beta ,\mu,\nu},\\
      \label{cac2.}      N_2K_{\alpha,\beta ,\mu,\nu}&=\beta  K_{\alpha,\beta ,\mu,\nu},\\      
\label{cac00.}      N_3K_{\alpha,\beta ,\mu,\nu}&=\mu K_{\alpha,\beta ,\mu,\nu}.
\end{align}\ees
\eep

\proof $ K_{\alpha,\beta ,\mu,\nu}$ is harmonic by Lemma \ref{cacko}.
    Writing
\bes    \begin{align}
K_{\alpha,\beta,\mu ,\nu}(z)&=
\int_\gamma(\tau z_1-z_{-2})^{\alpha+\nu}(z_2+\tau^{-1}z_{-1})^{\beta +\nu}z_3^\mu
\tau^{-\alpha-1}\frac{\d\tau}{2\pi\i}\\
&=
\int_\gamma( z_1-\tau^{-1}z_{-2})^{\alpha+\nu}(\tau z_2+z_{-1})^{\beta +\nu}z_3^\mu
\tau^{-\beta -1}\frac{\d\tau}{2\pi\i},
      \end{align}\ees
    we see that (\ref{cac1.}) and (\ref{cac2.}) follow from
    assumption (\ref{condi2}) by
Prop. \ref{cac0}. (\ref{cac00.}) is obvious. \qed

\bep If in addition to (\ref{condi2}) we assume that    \begin{align}
    (z_1-\tau^{-1}z_{-2})^{\alpha+\nu}(z_2+\tau^{-1}z_{-1})^{\beta +\nu}
\tau^{\nu}\Big|_{\tau(0)}^{\tau(1)}&=0,\label{condi1}\end{align}
and that we are allowed to differentiate under the integral sign, we obtain the recurrence relations
\begin{subequations}\begin{align}
 \label{cac3.} B_{-12}K_{\alpha,\beta ,\mu,\nu}&=(\beta +\nu) K_{\alpha+1,\beta -1,\mu,\nu},\\
       \label{cac4.} B_{1-2}K_{\alpha,\beta ,\mu,\nu}&=-(\alpha+\nu) K_{\alpha-1,\beta +1,\mu,\nu},\\
      \label{cac5.}  B_{12}K_{\alpha,\beta ,\mu,\nu}&=(\nu+1) K_{\alpha-1,\beta -1,\mu,\nu+1},\\
\label{cac6.}          B_{-1-2}K_{\alpha,\beta ,\mu,\nu}&=-(\alpha+\beta +\nu+1) K_{\alpha+1,\beta +1,\mu,\nu-1},\\\label{cac7.}
B_{1-3}K_{\alpha,\beta ,\mu,\nu}&=-(\alpha+\nu)K_{\alpha-1,\beta ,\mu+1,\nu},\\
\label{cac8.}B_{-1-3}K_{\alpha,\beta ,\mu,\nu}&=-(\beta +\nu)K_{\alpha+1,\beta ,\mu+1,\nu-1},\\\label{cac9.}
B_{2-3}K_{\alpha,\beta ,\mu,\nu}&=-(\beta+\nu)K_{\alpha,\beta-1 ,\mu+1,\nu},\\
\label{cac10.}B_{-2-3}K_{\alpha,\beta ,\mu,\nu}&=(\alpha +\nu)K_{\alpha,\beta+1 ,\mu+1,\nu-1}.\end{align}
\label{caco}  \end{subequations}
\eep

\proof Relations
(\ref{cac3.}),  (\ref{cac4.}),
(\ref{cac7.}),  (\ref{cac8.}), (\ref{cac9.}) and  (\ref{cac10.})
are elementary. They follow by simple differentiation under the integral sign  and do not need assumptions (\ref{condi1}) and (\ref{condi2}).

Relations (\ref{cac5.}) and (\ref{cac6.}) require assumption (\ref{condi2}) and follow by the following computations:
\begin{align}\label{caca1}
  &B_{12}
(z_1-\tau^{-1}z_{-2})^{\alpha+\nu}(z_2+\tau^{-1}z_{-1})^{\beta +\nu}
\tau^{\nu+1}\\\notag
=\,&\partial_{\tau^{-1}}
(z_1-\tau^{-1}z_{-2})^{\alpha+\nu}(z_2+\tau^{-1}z_{-1})^{\beta +\nu}
\tau^{\nu+1}\\\notag
&+(\nu+1)(z_1-\tau^{-1}z_{-2})^{\alpha+\nu}(z_2+\tau^{-1}z_{-1})^{\beta +\nu}
\tau^{\nu},\\\label{caca2}
&B_{-1-2}
(\tau z_1-z_{-2})^{\alpha+\nu}(\tau z_2+z_{-1})^{\beta +\nu}
\tau^{-\alpha-\beta -\nu-1}\\\notag
=&-\partial_\tau
(\tau z_1-z_{-2})^{\alpha+\nu}(\tau z_2+z_{-1})^{\beta +\nu}
\tau^{-\alpha-\beta -\nu-1}\\\notag
&-(\alpha+\beta +\nu+1)
(\tau z_1-z_{-2})^{\alpha+\nu}(\tau z_2+z_{-1})^{\beta +\nu}
\tau^{-\alpha-\beta -\nu-2},  
  \end{align}
where in (\ref{caca2}) we used yet another representation:
\beq
K_{\alpha,\beta,\mu ,\nu}(z):=
\int_\gamma(\tau z_1-z_{-2})^{\alpha+\nu}(\tau z_2+z_{-1})^{\beta +\nu}z_3^\mu
\tau^{-\alpha-\beta -\nu-1}\frac{\d\tau}{2\pi\i}.
\eeq
\qed


If in addition
\[\nu=\frac{-\alpha-\beta-\mu-1}{2},\]
then (\ref{wave}) is homogeneous of degree $-1$, so that we can reduce
it  to 4 dimensions.
Let us substitute the coordinates (\ref{subeq}), and then set $\tau=\frac{s}{u_1u_2}$, $s=t-w$:
\beq
  K_{\alpha,\beta,\mu,\nu}(u_1,u_2,u_3,r,p,w)=
2^{\frac12}  u_1^\alpha
u_2^\beta
u_3^\mu
r^{-\mu-1}p^\mu F(w),\end{equation}
\begin{align}\notag
F(w)&=\int_\gamma(s-1+w)^{\frac{\alpha-\beta-\mu-1}{2}}
(s+w)^{\frac{-\alpha+\beta-\mu-1}{2}}s^{\frac{-\alpha-\beta+\mu-1}{2}}\d s\\
&=\int_\gamma(t-1)^{\frac{\alpha-\beta-\mu-1}{2}}
t^{\frac{-\alpha+\beta-\mu-1}{2}}(t-w)^{\frac{-\alpha-\beta+\mu-1}{2}}\d t.
\label{cacy1}\end{align}
On the spherical section we can remove $r$ and $p$.
Therefore, the function $F$ given by (\ref{cacy1})
satisfies the hypergeometric equation: \beq
\cF_{\alpha,\beta,\mu}(w,\p_w)F(w)=0.\label{pwe}\eeq

From (\ref{caco}) we can also easily  obtain the recurrence relations for $F$. Note that in this list the recurrence relations corresponding to
$B_{1,3}$, $B_{-1,3}$, $B_{2,3}$ and $B_{-2,3}$ are missing.
However, they can be obtained after the reduction to $4$ dimensions by an application of the factorization formulas.

\subsection{Integral representations}
\label{Integral representations}

Below we independently prove (\ref{pwe}),
without going through the additional variables. We will use the classical parameters.
\bet
Let $[0,1]\ni \tau \overset{\gamma}\mapsto t(\tau)$ satisfy
\[t^{a-c+1}(1-t)^{c-b}(t-w)^{-a-1}\Big|_{t(0)}^{t(1)}=0.\]
Then
\beq 
\cF(a,b;c;w,\p_w)
\int_\gamma t^{a-c}(1-t)^{c-b-1}(t-w)^{-a}\d t=0.
\label{f4}\eeq
\label{intr}\eet

\proof 
We check that for any contour $\gamma$
\[\text{lhs of  (\ref{f4})}=-a \int_\gamma\d t\,
\p_t t^{a-c+1}(1-t)^{c-b}(t-w)^{-a-1}
.\]
\qed

Analogous (and nonequivalent) 
integral representations can be obtained by interchanging $a$ and $b$
in Theorem \ref{intr}.

The hypergeometric function with the type I normalization has
the integral representation
\begin{eqnarray}\label{eqa1}
&&\int_1^\infty t^{a-c}(t-1)^{c-b-1}(t-w)^{-a}\d t\\
&=&
 {\bf F}^\I  (a,b;c;w),\ \ \ \ \Re(c-b)>0,\ \Re b>0,\ \ \ w\not\in[1,\infty[.
\nonumber\end{eqnarray}
Indeed, by Theorem \ref{intr} the left hand side of  (\ref{eqa1})
 is annihilated by
the hypergeometric operator
(\ref{hy1-tra}). Besides, by  Euler's identity it equals
 $\frac{\Gamma(b)\Gamma(c-b)}{\Gamma(c)}$ at $0$. So does
the right hand side. Therefore,  (\ref{eqa1}) follows by the uniqueness of the solution by the Frobenius method.

\subsection{Integral representations of standard solutions}
\label{Integral representations of standard solutions}

 The integrand of \eqref{f4}
has four singularities: $\{0,1,\infty,w\}$. It is natural to chose
$\gamma$  as the interval joining a pair of  singularities. This choice leads to $6$ standard solutions with the $\rm I$-type normalization:
\begin{align*}
\text{  $\sim1$ at $0$:}\quad&[1,\infty];\\
\text{  $\sim w^{-\alpha}$ at $0$:}\quad&[0,w];\\
  \text{  $\sim1$  at $1$:}\quad
&[0,\infty];\\
\text{  $\sim(1-w)^{-\beta }$ at $1$:}\quad&[1,w];\\
\text{  $\sim w^{-a}$ at $\infty$:}\quad&[w,\infty];\\
  \text{  $\sim w^{-b}$ at $\infty$:}\quad
  &[0,1].
\end{align*}
Below we give explicit formulas.
To highlight their symmetry, we  use Lie-algebraic parameters.

\begin{align}
  \Re(1+\alpha)> |\Re(\beta -\mu )|:&\\
\int_1^\infty t^{\frac{-1-\alpha+\beta +\mu }{2}}
(t-1)^{\frac{-1+\alpha-\beta +\mu }{2}}(t-w)^{\frac{-1-\alpha-\beta -\mu }{2}}\d t&\notag\\
= {\bf F}^\I  _{\alpha,\beta ,\mu }(w),&\quad
w\not\in[1,\infty[;\nonumber\end{align}\begin{align}
\Re(1-\alpha)> |\Re(\beta -\mu )|:&\\
\int_0^w t^{\frac{-1-\alpha+\beta +\mu }{2}}
(1-t)^{\frac{-1+\alpha-\beta +\mu }{2}}(w-t)^{\frac{-1-\alpha-\beta -\mu }{2}}\d t&\notag\\
=
w^{-\alpha} {\bf F}^\I  _{-\alpha,\beta ,-\mu }(w),&\quad w\not\in]{-}\infty,0]{\cup}[1,\infty[
,\nonumber\notag\\
\int_w^0 (-t)^{\frac{-1-\alpha+\beta +\mu }{2}}
(1-t)^{\frac{-1+\alpha-\beta +\mu }{2}}(t-w)^{\frac{-1-\alpha-\beta -\mu }{2}}\d t&\notag\\
=
(-w)^{-\alpha} {\bf F}^\I  _{-\alpha,\beta ,-\mu }(w),&\quad w\not\in[0,\infty[
    ;\nonumber\end{align}\begin{align}
    \Re(1+\beta )> |\Re(\alpha-\mu )|:&\\
\int_{-\infty}^0 (-t)^{\frac{-1-\alpha+\beta +\mu }{2}}
(1-t)^{\frac{-1+\alpha-\beta +\mu }{2}}(w-t)^{\frac{-1-\alpha-\beta -\mu }{2}}\d t
&\notag\\=
 {\bf F}^\I  _{\beta ,\alpha ,\mu }(1-w),\quad
&w\not\in]-\infty,0];\notag\end{align}\begin{align}
\Re(1-\beta )> |\Re(\alpha+\mu )|:&\\
\int_w^1 t^{\frac{-1-\alpha+\beta +\mu }{2}}
(1-t)^{\frac{-1+\alpha-\beta +\mu }{2}}(t-w)^{\frac{-1-\alpha-\beta -\mu }{2}}\d t
&\notag\\=(1-w)^{-\beta } {\bf F}^\I  _{-\beta ,\alpha,-\mu }(1-w),
&\quad w\not\in]-\infty,0]\cup[1,\infty[,\notag\end{align}\begin{align}
\int_1^w t^{\frac{-1-\alpha+\beta +\mu }{2}}
(t-1)^{\frac{-1+\alpha-\beta +\mu }{2}}(w-t)^{\frac{-1-\alpha-\beta -\mu }{2}}\d t\notag\\
=(w-1)^{-\beta } {\bf F}^\I  _{-\beta ,\alpha,-\mu }(1-w),&\quad
w\not\in]-\infty,1];\notag\end{align}\begin{align}
\Re(1-\mu )> |\Re(\alpha+\beta )|:&\\
\int_w^\infty t^{\frac{-1-\alpha+\beta +\mu }{2}}
(t-1)^{\frac{-1+\alpha-\beta +\mu }{2}}(t-w)^{\frac{-1-\alpha-\beta -\mu }{2}}\d t&\notag\\
=
w^{\frac{-1-\alpha-\beta {+}\mu }{2}} {\bf F}^\I   _{-\mu ,\beta ,-\alpha}(w^{-1})
,&\quad
&w\not\in]-\infty,1],\notag\end{align}\begin{align}
\int_{-\infty}^w (-t)^{\frac{-1-\alpha+\beta +\mu }{2}}
(1-t)^{\frac{-1+\alpha-\beta +\mu }{2}}(w-t)^{\frac{-1-\alpha-\beta -\mu }{2}}\d t&\notag\\
=
(-w)^{\frac{-1-\alpha-\beta {+}\mu }{2}} {\bf F}^\I   _{-\mu ,\beta ,-\alpha}(w^{-1})
,&\quad
w\not\in]0,\infty];\notag\end{align}\begin{align}
\Re(1+\mu )> |\Re(\alpha-\beta )|:&\\
\int_0^1 t^{\frac{-1-\alpha+\beta -\mu }{2}}
(1-t)^{\frac{-1+\alpha-\beta +\mu }{2}}(t-w)^{\frac{-1-\alpha-\beta -\mu }{2}}\d t&\notag\\
=
(-w)^{\frac{-1-\alpha-\beta {-}\mu }{2}} {\bf F}^\I  _{\mu ,\beta ,\alpha}(w^{-1}),
&\quad w\not\in[0,\infty[,\notag\\
\int_0^1 t^{\frac{-1-\alpha+\beta -\mu }{2}}
(1-t)^{\frac{-1+\alpha-\beta +\mu }{2}}(w-t)^{\frac{-1-\alpha-\beta -\mu }{2}}\d t&\notag\\=
w^{\frac{-1-\alpha-\beta {-}\mu }{2}} {\bf F}^\I  _{\mu ,\beta ,\alpha}(w^{-1}),&
\quad w\not\in[-\infty,1[.\notag
\end{align}

    \subsection{Connection formulas}
\label{Connection formulas}

    Generically, each pair of standard solution is a basis of solutions to the hypergeometric equation.
For instance, we can use the pair of solutions $\sim 1$ and $\sim w^{-\alpha}$ at $0$ 
as one basis, and the pair $\sim w^{-a}$ and $\sim w^{-b}$ as another basis. We also assume that
$ w\not\in[0,\infty[$.

    Introduce the matrix
    \[
A_{\alpha,\beta,\mu}:=      \frac{\pi}{\sin(\pi\mu)}\left[\begin{array}{cc}\frac{-1}
{\Gamma\left(\frac{1+\alpha+\beta -\mu }{2}\right)
\Gamma\left(\frac{1+\alpha-\beta -\mu }{2}\right)}&
\frac{1}
{\Gamma\left(\frac{1+\alpha+\beta +\mu }{2}\right)
\Gamma\left(\frac{1+\alpha-\beta +\mu }{2}\right)}\\[2.5ex]
\frac{-1}
{\Gamma\left(\frac{1-\alpha-\beta -\mu }{2}\right)
\Gamma\left(\frac{1-\alpha+\beta -\mu }{2}\right)}&
\frac{1}{\Gamma\left(\frac{1-\alpha-\beta +\mu }{2}\right)
  \Gamma\left(\frac{1-\alpha+\beta +\mu }{2}\right)}
  \end{array}\right].
\]
Then
\begin{align}\label{matri2}&\left[\begin{array}{c} {\bf F}  _{\alpha,\beta ,\mu }(w)\\[2.5ex]
    (-w)^{-\alpha} {\bf F}  _{-\alpha,\beta ,-\mu }(w)
    \end{array}\right]\\=\quad &A_{\alpha,\beta,\mu}
  \left[\begin{array}{c} 
(-w)^{\frac{-1-\alpha-\beta -\mu }{2}}  {\bf F}  _{\mu ,\beta ,\alpha}(w^{-1})\\[2.5ex]
    (-w)^{\frac{-1-\alpha-\beta +\mu }{2}}  {\bf F}  _{-\mu ,\beta ,-\alpha}(w^{-1})
  \end{array}\right]
  .\notag
\end{align}

Note that in the Lie-algebraic parameters the matrix $A_{\alpha,\beta,\mu}$ has a very symmetric form. Here are some of its properties:
\begin{align}
A_{\alpha,\beta,\mu}=A_{\alpha,-\beta,\mu}&=-
\begin{bmatrix}0&1\\1&0
\end{bmatrix}A_{-\alpha,\beta,-\mu}
\begin{bmatrix}0&1\\1&0
\end{bmatrix}=A_{\mu,\beta,\alpha}^{-1},\\
\det A_{\alpha,\beta,\mu}&=-\frac{\sin(\pi\alpha)}{\sin(\pi\mu)}.
\end{align}

Relation (\ref{matri2}) can be derived from the integral representations. Indeed,
consider $\Im w<0$. Take the branches of the powers of $-t$ and $1-t$ and $w-t$   continued from the left clockwise onto the upper halfplane.
Then (under some conditions on $\alpha,\beta,\mu$) we can write
\begin{align*}
  \Bigg(\int_{-\infty}^0+\int_0^1+\int_1^{+\infty}\Bigg)
(- t)^{\frac{-1-\alpha+\beta \pm\mu }{2}}
  (1-t)^{\frac{-1+\alpha-\beta \pm\mu }{2}}(w-t)^{\frac{-1-\alpha-\beta \mp\mu }{2}}\d t&=0.\label{alberti1}
\end{align*}
We obtain
\begin{align*}
   {\bf F}^\I_{\beta ,\alpha ,\pm\mu }(1{-}w)
   -\e^{\i\pi\alpha}
 ( - w)^{\frac{-1-\alpha-\beta \mp\mu }{2}} {\bf F}^\I  _{\pm\mu ,\beta ,\alpha}(w^{-1})
   -\i \e^{\i\pi\frac{\alpha+\beta\mp\mu}{2}}{\bf F}^\I_{\alpha,\beta ,\pm\mu }(w)&=0.
\end{align*}
Using
\[{\bf F}_{\alpha,\beta,\mu}^\I(w)=
\Gamma\Big(\frac{1+\alpha+\beta-\mu}{2}\Big)\Gamma\Big(\frac{1+\alpha-\beta+\mu}{2}\Big)
{\bf F}_{\alpha,\beta ,\mu }(w),
\] we express everything in terms of ${\bf F}$. We eliminate $   {\bf F}_{\beta ,\alpha ,\mu }(1{-}w)={\bf F}_{\beta ,\alpha ,-\mu }(1{-}w)$.
We find
\begin{align*}
  {\bf F}  _{\alpha,\beta ,\mu }(w)\,=\,&
  -  \frac{\pi
  (-w)^{\frac{-1-\alpha-\beta -\mu }{2}}  {\bf F}  _{\mu ,\beta ,\alpha}(w^{-1})
}{\sin(\pi\mu)
\Gamma\left(\frac{1+\alpha+\beta -\mu }{2}\right)
\Gamma\left(\frac{1+\alpha-\beta -\mu }{2}\right)}
\\
&
+\frac{\pi(-w)^{\frac{-1-\alpha-\beta +\mu }{2}}  {\bf F}  _{-\mu ,\beta ,-\alpha}(w^{-1})}
{\sin(\pi\mu)\Gamma\left(\frac{1+\alpha+\beta +\mu }{2}\right)
  \Gamma\left(\frac{1+\alpha-\beta +\mu }{2}\right)}    
 ,\end{align*}
  which is the first line of  
(\ref{matri2}).
A similar argument, starting with the integral  $ \int_{-\infty}^0+\int_0^w+\int_w^{+\infty}$, yields the second line of (\ref{matri2}).

\section{Laplacian in 3 dimensions and the Gegenbauer equation}
\label{s7}
\init
The Gegenbauer equation
is equivalent to a subclass of the ${}_2\cF_1$ equation.
Nevertheless, not all its symmetries are directly inherited from the symmetries of the ${}_2\cF_1$ equation. Therefore it deserves a separate treatment, which is given in this section.
We start from the Laplacian in 5 dimensions,  pass through 3 dimensions, and eventually we derive the
Gegenbauer equation.

This section is to a large extent parallel to the previous one, devoted to the ${}_2\cF_1$ equation.  The number of symmetries, parameters, etc. is now smaller than in the previous section, since we are in  lower dimensions. Nevertheless, some things are here more complicated and less symmetric.
This is related to the fact that the number of dimensions is odd, which corresponds to a less symmetric orthogonal group and Lie algebra.

Let us describe the main  steps of our  derivation of the Gegenbauer equation, even though they are almost the same as  for the ${}_2\cF_1$ equation.
\ben
\item
  We start from the $3+2=5$ dimensional ambient space, with the obvious representation of $\so(5)$ and $\mathrm{O}(5)$, and the Laplacian $\Delta_5$.
  \item\label{it2.} We go to the representations
    $\so(5)\ni B\mapsto B^{\diamond,\eta}$ and
    $\mathrm{O}(5)\ni \alpha\mapsto \alpha^{\diamond,\eta}$
and to the reduced Laplacian $\Delta_5^\diamond$. The most relevant values of $\eta$ are $1-\frac{3}{2}=-\frac12$ and  $-1-\frac{3}{2}=-\frac52$.
\item\label{it3.} We fix a section $\gamma$ of the null quadric, obtaining the representations $B^{\gamma,\eta}$ and $\alpha^{\gamma,\eta}$, as well as  the operator $\Delta_5^\gamma$, acting on an appropriate pseudo-Riemannian  $3$ dimensional manifold.
\item\label{it4.} We choose coordinates $w,u_2,u_3$, so that the Cartan elements can be expressed in terms of  $u_2$, $u_3$. We compute
  $B^{\gamma,\eta}$,   $\alpha^{\gamma,\eta}$ and  $\Delta_5^\gamma$ in the new coordinates.
\item\label{it5.} We make an ansatz that diagonalizes the Cartan elements. The eigenvalues, denoted by $\alpha$, $\lambda$,
  become parameters.
  $B^{\gamma,\eta}$,   $\alpha^{\gamma,\eta}$ and   $\Delta_5^\gamma$ involve now only the single variable $w$. 
  $\Delta_5^\gamma$ turns out to be the  Gegenbauer operator. We obtain its transmutation relations and discrete symmetries.
\een

Again, we choose a special section which makes computations relatively easy.
We perform Steps \ref{it2.}, \ref{it3.} and \ref{it4.} at once, by choosing convenient coordinates  $w,r,p,u_2,u_3$ in $5$ dimensions.
After the reductions of Steps \ref{it2.} and \ref{it3.}, we are left with the variables $w,u_2,u_3$, and we can perform Step \ref{it5.}.
  

The remaining material of this section is  parallel to the analogous material of the previous section except for Subsect. \ref{Quadratic transformation},
which describes a quadratic relation  reducing the Gegenbauer equation to the ${}_2\cF_1$ equation. We describe a derivation of this relation starting from the level of the ambient space.

\subsection{$\so(5)$ in 5 dimensions}
\label{sub-geg1}

We consider $\rr^5$ with the coordinates
\beq z_0,z_{-2},z_2,z_{-3},z_3\label{geg1}\eeq
and the scalar product given by
\beq \langle z|z\rangle=z_0^2+2z_{-2}z_2+2z_{-3}z_3.\label{geg5}\eeq
Note that we omit the indices $-1,1$; this  makes it easier to compare
$\rr^5$ with $\rr^6$.

The Lie algebra  $\so(5)$ acts naturally on $\rr^5$. Below we describe its natural basis. Then we consider the Weyl group $B_2$ acting on functions on $\rr^5$. For brevity, we list only elements
from its subgroup  $B_2\cap \SO(5)$. Finally, we write down the Laplacian.

\medskip

\noindent{\bf Lie algebra $\so(5)$.}
Cartan algebra
\bes\begin{align}
  N_2&=-z_{-2}\partial_{z_{-2}}+z_{2}\partial_{z_{2}},\\
   N_3&=-z_{-3}\partial_{z_{-3}}+z_{3}\partial_{z_{3}}.
  \end{align}\ees
Root operators
\begin{subequations}\label{coor5}
\begin{align}
  B_{0,-2}&=z_{0}\p_{z_{-2}}-z_{2}\p_{z_0},\\
  B_{0,2}&=z_{0}\p_{z_2}-z_{-2}\p_{z_0},\\
B_{0,-3}&=z_{0}\p_{z_{-3}}-z_{3}\p_{z_0},\\
B_{0,3}&=z_{0}\p_{z_3}-z_{-3}\p_{z_0};
\end{align}
\begin{align}  B_{-3,-2}&=z_{3}\p_{z_{-2}}-z_{2}\p_{z_{-3}},\\
  B_{3,2}&=z_{-3}\p_{z_2}-z_{-2}\p_{z_3},\\
B_{3,-2}&=z_{-3}\p_{z_{-2}}-z_{2}\p_{z_3},\\
B_{-3,2}&=z_{3}\p_{z_2}-z_{-2}\p_{z_{-3}}.
\end{align}
\end{subequations}

\noindent{\bf Weyl symmetries} 
\begin{subequations}
  \begin{align}
  \sigma_{23} K(z_0,z_{-2},z_2,z_{-3},z_{3})=&K(z_0,z_{-2},z_2,z_{-3},z_{3}),\\
  \tau_{2-3} K(z_0,z_{-2},z_2,z_{-3},z_{3})=&K(-z_0,z_{-2},z_2,z_{3},z_{-3}),\\
  \sigma_{-2-3}K(z_0,z_{-2},z_2,z_{-3},z_{3})=&K(z_0,z_{2},z_{-2},z_{3},z_{-3}),\\
  \tau_{-23}K(z_0,z_{-2},z_2,z_{-3},z_{3})=&K(-z_0,z_{2},z_{-2},z_{-3},z_{3});
  \end{align}\begin{align}
    \sigma_{32} K(z_0,z_{-2},z_2,z_{-3},z_{3})=&K(z_0,z_{-3},z_3,z_{-2},z_{2}),\\
  \tau_{3-2} K(z_0,z_{-2},z_2,z_{-3},z_{3})=&K(-z_0,z_{-3},z_3,z_2,z_{-2}),\\
  \sigma_{-3-2}K(z_0,z_{-2},z_2,z_{-3},z_{3})=&K(z_0,z_{3},z_{-3},z_{2},z_{-2}),\\
    \tau_{-32}K(z_0,z_{-2},z_2,z_{-3},z_{3})=&K(-z_0,z_{3},z_{-3},z_{-2},z_{2}).
\end{align}\label{weylge1}
\end{subequations}

\noindent{\bf Laplacian}
\beq\Delta_{5}=\partial_{z_0}^2+2\partial_{z_{-2}}\partial_{z_{2}}+
2\partial_{z_{-3}}\partial_{z_{3}}
.\eeq

\subsection{$\so(5)$ on the spherical section}

In this subsection
we perform Steps \ref{it2.}, \ref{it3.} and \ref{it4.}, as described in the introduction to this section. Recall that Step \ref{it2.} involves
restricting to the null quadric
\[\cV^4:=\{z\in\rr^5\ :\ z_0^2+2z_{-2}z_2+2z_{-3}z_3=0\}.\]
To perform Step \ref{it3.} we need to fix a section of this quadric.
We choose the  section given by the equations
\[1 =z_0^2+2z_{-2}z_2=-2z_3z_{-3}.\]
We will call it the {\em spherical section}, because it is
$\ss^2(1)\times\ss^1(-1)$.  The superscript used for this section
will be ``$\sph$'' for spherical.

We introduce the coordinates $w,r,p,u_2,u_3$:
\begin{subequations}\begin{align}
r& =\sqrt{z_0^2+2z_{-2}z_2}\;,&&\\
w& =\frac{z_0}{\sqrt{2z_{-2}z_2+z_0^2}}\;,\ \ \ 
    &u_2& =\frac{\sqrt2
      z_2}{\sqrt{z_0^2+2z_{-2}z_{2}}}\;,\phantom{\sqrt{\frac{1}{2}}}\\
    p& =\sqrt{-2z_3z_{-3}}\;,
    &u_3& =\sqrt{-\frac{z_{3}}{z_{-3}}}\;.
  \end{align}\label{cou}
\end{subequations}
Here is the inverse transformation:
\begin{subequations}\begin{align}
  z_0&=wr,&\quad z_{-2}&=\frac{r(1-w^2)}{\sqrt2u_2},&\quad
  z_2&=\frac{u_2r}{\sqrt2},\\
  &&\quad z_{-3}&=-\frac{p}{\sqrt2 u_3},&\quad z_3&=\frac{pu_3}{\sqrt{2}}.
  \end{align}\end{subequations}

Similarly as in the previous section,
the null quadric in these coordinates is given by
 $r^2=p^2$. We choose the sheet $r=p$.
The generator of dilations is
\[A_{5}=r\ddr+p\ddp.\]
The spherical section is given by the condition
$r^2=1$.

\medskip

\noindent{\bf Lie algebra $\so(5)$.}
Cartan operators
\begin{align*}
N_{2}^\sph& =u_2\ddu{2}\,,\\
N_{3}^\sph& =u_3\ddu{3}\,.
\end{align*}
Roots
\begin{align*}
  B_{0,-2}^\sph&=-\frac{u_2}{\sqrt2}\p_w,\\
  B_{0,2}^\sph&=\frac{1}{\sqrt2 u_2}\Big((w^2-1)\p_w+2wu_2\p_{u_2}\Big),\\ B_{0,-3}^{\sph,\eta}&=\frac{u_3}{\sqrt2}\Big((w^2-1)\p_w+wu_2\p_{u_2}+wu_3\p_{u_3}-w\eta\Big),\\ B_{0,3}^{\sph,\eta}&=\frac{1}{\sqrt2u_3}\Big((1-w^2)\p_w-wu_2\p_{u_2}+wu_3\p_{u_3}+w\eta\Big);
\end{align*}
\begin{align*}  B_{-3,-2}^{\sph,{\eta}}&=\frac{u_2u_3}{2}
  \left(-w\p_w-u_2\p_{u_2}-u_3\p_{u_3}+\eta\right),\\
B_{3,2}^{\sph,{\eta}}&=\frac{1}{2u_2u_3}\left(w(1-w^2)\p_w-(1+w^2)u_2\p_{u_2}+(w^2-1)u_3\p_{u_3}+(w^2-1)\eta\right),\\
  B_{3,-2}^{\sph,{\eta}}&=\frac{u_2}{2u_3}\left(w\p_w+
u_2\p_{u_2}-u_3\p_{u_3}-\eta\right),\\
B_{-3,2}^{\sph,{\eta}}&=\frac{ u_3}{2u_2}\left(w(w^2-1)\p_w+(1+w^2)u_2\p_{u_2}+(w^2-1)u_3\p_{u_3}+(1-w^2)\eta
\right). 
\end{align*}

\noindent{\bf Weyl symmetries} 
  \begin{align*}
  \sigma_{23}^{\sph,\eta} f(w,u_2,u_3)=& f(w,u_2,u_3),\\
  \tau_{2-3}^{\sph,\eta}  f(w,u_2,u_3)=& f\Big(-w,u_2,\frac{1}{u_3}\Big),\\
  \sigma_{-2-3}^{\sph,\eta} f(w,u_2,u_3)=& f\Big(w,\frac{1-w^2}{u_2},\frac{1}{u_3}\Big),\\
  \tau_{-23}^{\sph,\eta} f(w,u_2,u_3)=& f\Big(-w,\frac{1-w^2}{u_2},u_3\Big);
  \end{align*}\begin{align*}
    \sigma_{32}^{\sph,\eta}  f(w,u_2,u_3)=&(w^2-1)^{\frac{\eta}{2}}f\Big(\frac{w}{\sqrt{w^2-1}},\frac{u_3}{\sqrt{w^2-1}},\frac{u_2}{
      \sqrt{w^2-1}}\Big),\\
  \tau_{3-2}^{\sph,\eta}  f(w,u_2,u_3)=&(w^2-1)^{\frac{\eta}{2}}f\Big(\frac{-w}{\sqrt{w^2-1}},\frac{u_3}{\sqrt{w^2-1}},\frac{
      \sqrt{w^2-1}}{u_2}\Big),\\
  \sigma_{-3-2}^{\sph,\eta} f(w,u_2,u_3)=&(w^2-1)^{\frac{\eta}{2}}f\Big(\frac{w}{\sqrt{w^2-1}},\frac{-1}{u_3\sqrt{w^2-1}},\frac{
      \sqrt{w^2-1}}{u_2}\Big),\\
    \tau_{-32}^{\sph,\eta} f(w,u_2,u_3)=&(w^2-1)^{\frac{\eta}{2}}f\Big(\frac{-w}{\sqrt{w^2-1}},\frac{-1}{u_3\sqrt{w^2-1}},\frac{u_2}{
      \sqrt{w^2-1}}\Big).
\end{align*}

\noindent
{\bf Laplacian}

\begin{align}\notag
\Delta_5^\sph
=&(1-w^2)\p_w^2
-2(1+u_2\p_{u_2})w\p_w
  -\Big(u_2\p_{u_2}+\frac12\Big)^2+(u_3\p_{u_3})^2.
\label{lapiq}  \end{align}

 Let us sketch the computations that lead to (\ref{lapiq}). Using
\begin{align*}
  \p_{z_0}&=\frac1r\Big(wr\p_r-wu_2\p_{u_2}+(1-w^2)\p_w\Big),\\
  \p_{z_{-2}}&=\frac{u_2}{\sqrt{2}r}\Big(r\p_r-u_2\p_{u_2}-w\p_w\Big),\\
  \p_{z_2}&=\frac{1}{\sqrt{2}ru_2}\Big((1-w^2)r\p_r+(1+w^2)u_2\p_{u_2}+(w^2-1)w\p_w\Big)
  ,\\
  \p_{z_{-3}}&=\frac{u_3}{\sqrt2 p}\Big(u_3\p_{u_3}-p\p_p\Big),\\
    \p_{z_3}&=\frac{1}{\sqrt2 pu_3}\Big(u_3\p_{u_3}+p\p_p\Big).
  \end{align*}
we change the variables in the Laplacian:
\begin{align}
  \Delta_5=&\frac{1}{r^2}\Big((1-w^2)\p_w^2-2(1+u_2\p_{u_2})w\p_w
  -(u_2\p_{u_2})^2-u_2\p_{u_2}\notag\\  &+(r\p_r)^2+r\p_r\Big)+\frac{1}{p^2}\Big(-(p\p_p)^2+(u_3\p_{u_3})^2\Big).\end{align}
Now,
\begin{align*}
(r\ddr)^2+r\ddr-\frac{r^2}{p^2}(p\ddp)^2& =\biggl(r\ddr-p\p_p+\frac{1}{2}\biggr)\biggl(r\ddr+p\p_p+\frac{1}{2}\biggr)\\
&\quad+\biggl(1-\frac{r^2}{p^2}\biggr)(p\p_p)^2-\frac{1}{4}.
\end{align*}
Therefore, using $r^2=p^2$, $r\p_r+p\p_p=-\frac12$, we obtain
\begin{align}\notag
\Delta_5^\diamond&
=\frac{1}{r^2}\Big((1-w^2)^2\p_w^2
-2(1+u_2\p_{u_2})w\p_w\\
&  -\big(u_2\p_{u_2}+\frac12\big)^2+(u_3\p_{u_3})^2\Big).
  \end{align}
 To obtain the Laplacian at the spherical section we drop $\frac1{r^2}$.

\subsection{The Gegenbauer operator}

Let us make the ansatz
\beq
f(u_2,u_3,w)= u_2^\alpha u_3^\lambda S(w).\label{form0-}\eeq
Clearly,
\bes\begin{eqnarray}
N_2^\sph f&=&\alpha f,\\
N_3^\sph f&=&\lambda f,\\
u_2^{-\alpha} u_3^{-\lambda} \Delta_5^\sph  f&
=&  \cS_{\alpha,\lambda}(w,\p_w)S(w),\end{eqnarray}\ees
where
\begin{eqnarray}
\cS_{\alpha ,\lambda }(w,\p_w)
&:=&(1-w^2)\p_w^2-2(1+\alpha )w\p_w
+\lambda ^2-\Big(\alpha +\12\Big)^2\end{eqnarray}
is the  {\em Gegenbauer operator}.
Here is another parametrization of the Gegenbauer operator, which we call {\em classical}:
\beq
\cS({a},{b};w,\p_w):=(1-w^2)\p_w^2-({a}+{b}+1)w\p_w-{a}{b}.
\label{gego}\eeq
Here is the relationship between the classical and Lie-algebraic  parameters:
\bes\begin{align}
&\alpha :=\frac{{a}+{b}-1}{2},& \lambda: =\frac{{b}-{a}}{2},\\
&{a}=\frac12+\alpha -\lambda ,&{b}=\frac12+\alpha +\lambda .
\end{align}
\ees

The Gegenbauer operator is  the   ${}_2\cF_1$ operator with its finite singular points moved to $-1$ and $1$, which in addition is reflection invariant.
Because of the reflection invariance, the third classical parameter can be obtained from the first two: $c=\frac{a+b+1}{2}$. Therefore, we  use only  ${a},{b}\in\cc$
 as the (classical) parameters of the Gegenbauer equation.

We can  reduce the Gegenbauer equation to the ${}_2\cF_1$ equation by
two affine transformations. They move the singular points from $-1$, $1$ to
$0$, $1$ or $1$, $0$:
\beq\begin{array}{l}
\cS({a},{b};w,\p_w)=\cF\big({a},{b};\frac{{a}+{b}+1}{2}
;v,\p_v\big),\end{array}\label{ha2}\eeq
where
\bes\begin{align}&v=\frac{1-w}{2},\ \ \ w=1-2v,\\
\hbox{or}\qquad&v=\frac{1+w}{2},\ \ \ w=-1+2v.
\end{align}\ees
In the Lie-algebraic parameters
\beq\cS_{\alpha,\lambda}(w,\p_w)=\cF_{\alpha,\alpha,2\lambda}(v,\p_v).\eeq

\subsection{Quadratic transformation}
\label{Quadratic transformation}

Let us go back to 6 dimensions and the Laplacian
\beq\Delta_6=2\partial_{z_{-1}}\partial_{z_{1}}+
2\partial_{z_{-2}}\partial_{z_{2}}+
2\partial_{z_{-3}}\partial_{z_{3}}
.\label{sq20}\eeq
Let us use the reduction described in Subsect.
\ref{Dimensional reduction}.
Introduce new variables
\beq z_0:=\sqrt{2z_{-1}z_{1}},\quad
u:=\sqrt{\frac{z_{1}}{z_{-1}}}.\label{coo1}\eeq
  In the new variables,
  \begin{align}
        N_1=&u\p_u,\\\notag
    \Delta_6
    =&\Big(\partial_{z_0}+\frac1{2z_0}\Big)^2
    -\frac1{z_0^2}\Big(u\partial_u-\frac12\Big)
    \Big(u\partial_u+\frac12\Big)\\&+
   2\partial_{z_{-2}}\partial_{z_{2}}+
2\partial_{z_{-3}}\partial_{z_{3}}.
  \end{align}
  Therefore,
\bes  \begin{align}
  (uz_0)^{\frac12}\Delta_6  (uz_0)^{-\frac12}
   & =
    -\frac1{z_0^2}N_1
    \Big(N_1- 1\Big)+
    \Delta_5,\\
      (u^{-1}z_0)^{\frac12}\Delta_6  (u^{-1}z_0)^{-\frac12}
   & =
    -\frac1{z_0^2}N_1
    \Big(N_1+1\Big)+
    \Delta_5.
    \end{align}\ees
  Compare the coordinates
  the coordinates (\ref{coo}) for 6 dimensions and
   (\ref{cou}) for 5 dimensions.
  The coordinates $p$, $u_3$ are the same.
  Taking into account $z_0:=\sqrt{2z_{-1}z_{1}}$,  the coordinates $r,u_2$ also coincide. This is not the case of
  $w$, so let us rename  $w$ from (\ref{cou})
  as $v$. We then have $w=v^2$.
  We also have
  \[uz_0=\sqrt2z_1=u_1r,\quad u^{-1}z_0=\sqrt 2z_{-1}=rw u_1^{-1}.\]
  Hence on functions that do not depend on $u$ we obtain
\bes  \begin{align}
    r^{\frac12}u_1^{\frac12}\Delta_6 r^{-\frac12}u_1^{-\frac12}&=
    \Delta_5 ,\\
        r^{\frac12}u_1^{-\frac12}v\Delta_6 r^{-\frac12}u_1^{\frac12}v^{-1}&=
        \Delta_5 .
        \end{align}
\ees

This implies that a quadratic substitution
 transforms the ${}_2\cF_1$ operator with $\alpha=\pm\frac12$ into the Gegenbauer operator.
Explicitly, if \[w=v^2,\ \ \  \ 
v=\sqrt w;\]
then in the classical parameters
\bes\begin{align}\cS({a},{b};v,\p_v)&=
4\cF\Big(\frac{{a}}{2},\frac{{b}}{2};\12;w,\p_w\Big),\\[2ex]
v^{-1}\cS({a},{b};v,\p_v)v&=
4\cF\Big(\frac{{a}+1}{2},\frac{{b}+1}{2};\frac32;w,\p_w\Big),
\label{ha1}\end{align}\ees
and in the Lie-algebraic parameters
\bes\begin{eqnarray}
  \label{1/2}\cS_{\alpha,\lambda}(v,\p_v)&=&4\cF_{-\12,\alpha,\lambda}(w,\p_w),
  \label{pol1}\\
\label{-1/2}v^{-1}\cS_{\alpha,\lambda}(v,\p_v)v&=&4\cF_{\12,\alpha,\lambda}(w,\p_w).
\label{pol2}\end{eqnarray}
\ees

\subsection{Transmutation relations and discrete symmetries}
\label{symcom}

We have the following generalized symmetries:
\bes\begin{eqnarray}
B^{\sph,-\frac52}\Delta_{5}^\sph&=&\Delta_{5}^\sph B^{\sph,-\frac12},\ \ \ \ B\in \so(5);\label{vv1}\\
\alpha^{\sph,-\frac52}\Delta_{5}^\sph&=&\Delta_{5}^\sph\alpha^{\sph,-\frac12},\ \ \ \alpha\in \mathrm{O}(5).\label{vv2}\end{eqnarray}\ees
Equality (\ref{vv1}) applied to the roots of $\so(5)$ yield the following transmutation relations:

\[\hspace{-3ex}\begin{array}{rrl}
&\p_w&
\cS_{\alpha ,\lambda } \\[0.4ex]
&=\ \ \ \cS_{\alpha +1,\lambda }&\p_w,\\[0.5ex]
&((1-w^2)\p_w-2\alpha w)&
\cS_{\alpha ,\lambda }\\[0.4ex]
&=\ \ \ \cS_{\alpha -1,\lambda }&((1-w^2)\p_w-2\alpha w),\\[0.5ex]
&((1-w^2)\p_w-
(\alpha +\lambda +\frac12 )w)&(1-w^2)\cS_{\alpha ,\lambda }\\[0.4ex]
&=\ \ \ (1-w^2)\cS_{\alpha ,\lambda +1}&
((1-w^2)\p_w-
(\alpha +\lambda +\12)w),\\[0.5ex]
&((1-w^2)\p_w-(\alpha -\lambda +\frac12)w)&
(1-w^2)\cS_{\alpha ,\lambda }\\[0.4ex]
&=\ \ \ (1-w^2)\cS_{\alpha ,\lambda -1}&((1-w^2)\p_w-(\alpha -\lambda +\12)w);
\\[1ex]
&(w\p_w+
\alpha -\lambda +\frac12)&w^2\cS_{\alpha ,\lambda }\\[0.4ex]
&=\ \ \ w^2\cS_{\alpha +1,\lambda -1}&(w\p_w+
\alpha -\lambda +\12),\\[0.5ex]
&(w(1{-}w^2)\p_w{-}\alpha{-}\lambda{+}\frac12{-}(\alpha{-}\lambda{+}\frac12)w^2)
&w^2\cS_{\alpha ,\lambda }\\[0.4ex]
&=\ \ \ w^2\cS_{\alpha -1,\lambda +1}
&(w(1{-}w^2)\p_w
{-}\alpha{-}\lambda{+}\frac12{-}(\alpha{-}\lambda{+}\frac12)w^2),
\\[0.5ex]
&(w\p_w+
\alpha -\lambda +\frac12)&w^2\cS_{\alpha ,\lambda }\\[0.4ex]
&=\ \ \ w^2\cS_{\alpha +1,\lambda +1}&(w\p_w+
\alpha -\lambda +\12),\\[0.5ex]
&(w(1{-}w^2)\p_w
{-}\alpha{+}\lambda{+}\frac12{-}(\alpha{+}\lambda{+}\frac12)w^2
)&
w^2\cS_{\alpha ,\lambda }\\[0.4ex]
&=\ \ \ w^2\cS_{\alpha -1,\lambda -1}&(w(1{-}w^2)\p_w
{-}\alpha{+}\lambda{+}\frac12{-}(\alpha{+}\lambda{+}\frac12)w^2).
\end{array}\]

Next we describe  discrete symmetries of the Gegenbauer operator, which follow from
Relation (\ref{vv2}) applied to Weyl symmetries.
All the operators below equal $\cS_{\alpha ,\lambda }(w,\p_w)$ for the
appropriate $w$:
\[\begin{array}{rrcl}
w=\pm v:
&&\cS_{\alpha ,\pm\lambda }(v,\p_v),&\\[1ex]
w=\pm v:
&(v^2-1)^{-\alpha }&\cS_{-\alpha ,\mp\lambda }(v,\p_v)&
(v^2-1)^{\alpha },\\[1ex]
w=\frac{\pm v}{(v^2-1)^{\12}}:
& (v^2-1)^{\12(\alpha +\lambda +\frac52)}
&\cS_{\lambda ,\pm\alpha }(v,\p_v)
& (v^2-1)^{\12(-\alpha -\lambda -\12)},\\[1ex]
w=\frac{\pm v}{(v^2-1)^{\12}}:
& (v^2-1)^{\12(\alpha -\lambda +\frac52)}&
\cS_{-\lambda ,\mp\alpha }(v,\p_v)&
 (v^2-1)^{\12(-\alpha +\lambda -\12)}.
\end{array}\]

Note that we use  $\pm$ to describe two symmetries at once. Therefore, the above list has all $2\times 4=8$ symmetries corresponding to the lists of Weyl symmetries (\ref{weylge1}).

\subsection{Factorizations of the Laplacian}
\label{Factorizations of the Laplacian2}
In the Lie algebra $\so(5)$ represented on $\rr^5$
we have 3 distinguished Lie subalgebras:
two isomorphic to $\so(3)$ and
 one isomorphic to $\so(4)$:
 \beq \so_{02}(3),\; \; \so_{03}(3),\; \; \so_{23}(4),\eeq
 where we use an obvious notation.
 By  (\ref{casimir1}) and
 (\ref{casimir2}), the corresponding Casimir operators are
\begin{subequations}
\begin{align}
\cC_{02}
&=2B_{0,-2}B_{0,2}-\Big(N_2-\frac12\Big)^2+\frac14\\
&=2B_{0,2}B_{0,-2}-\Big(N_2+\frac12\Big)^2+\frac14,\\
\cC_{03}
&=2B_{0,-3}B_{0,3}-\Big(N_3-\frac12\Big)^2+\frac14\\
&=2B_{0,3}B_{0,-3}-\Big(N_3+\frac12\Big)^2+\frac14,\\
\cC_{23}
&=4B_{2,3}B_{-2,-3}-(N_2+N_3+1)^2+1\\
&=4B_{-2,-3}B_{2,3}-(N_2+N_3-1)^2+1\\
&=4B_{2,-3}B_{-2,3}-(N_2-N_3+1)^2+1\\
&=4B_{-2,3}B_{2,-3}-(N_2-N_3-1)^2+1.
\end{align}\label{facto1}
\end{subequations}

After the reduction described in 
  (\ref{deq3a}) and (\ref{deq1a}), we obtain the
identities
\begin{subequations}
\begin{align}
\label{factor2}
(z_0^2+2z_{-2}z_2)\Delta_{5}^\diamond&=-\frac14+\cC_{02}^{\diamond,-\frac12}
+(N_3^{\diamond,-\frac12})^2,
\\\label{factor3}
(z_0^2+2z_{-3}z_3)\Delta_{5}^\diamond&=-\frac14+\cC_{03}^{\diamond,-\frac12}
+(N_2^{\diamond,-\frac12})^2,\\
\label{factor1}(2z_{-2}z_2+2z_{-3}z_3)\Delta_{5}^\diamond&=-\frac34+\cC_{23}^{\diamond,-\frac12}.
\end{align}\label{facto2}
\end{subequations}

Inserting (\ref{facto1}) into (\ref{facto2}), we obtain
\bes
\begin{align}
&(z_0^2+2z_{-2}z_2)\Delta_{5}^\diamond\notag\\
=&2B_{0,-2}B_{0,2}-\Big(N_2+N_3-\frac12\Big)\Big(N_2-N_3-\frac12\Big)\\
=&2B_{0,2}B_{0,-2}-\Big(N_2+N_3+\frac12\Big)\Big(N_2-N_3+\frac12\Big),
  \\
  &(z_0^2+2z_{-3}z_3)\Delta_{5}^\diamond\notag\\
=&2B_{0,-3}B_{0,3}-\Big(N_2+N_3-\frac12\Big)\Big(-N_2+N_3-\frac12\Big)
\\
=&2B_{0,3}B_{0,-3}-\Big(N_2+N_3+\frac12\Big)\Big(-N_2+N_3+\frac12\Big),\\
&(2z_{-2}z_2+2z_{-3}z_3)\Delta_{5}^\diamond\notag\\
=&4B_{2,3}B_{-2,-3}-\Big(N_2+N_3+\frac32\Big)\Big(N_2+N_3+\frac12\Big)
\\
=&4B_{-2,-3}B_{2,3}-\Big(N_2+N_3-\frac32\Big)\Big(N_2+N_3-\frac12\Big)\\
=&4B_{2,-3}B_{-2,3}-\Big(N_2-N_3+\frac32\Big)\Big(N_2-N_3+\frac12\Big)\\
=&4B_{-2,3}B_{2,-3}-\Big(N_2-N_3-\frac32\Big)\Big(N_2-N_3-\frac12\Big),
\end{align}\ees
where all the $B$ and $N$ operators need to have the superscript ${}^{\diamond,-\frac12}$.

If we use the spherical section,
we need to make the replacements
\begin{subequations}
\begin{align}
\label{factor2-}
z_0^2+2z_{-2}z_2&\quad\to\quad 1,
\\\label{factor3-}
z_0^2+2z_{-3}z_3&\quad\to\quad  w^2-1,\\
\label{factor1-}
2z_{-2}z_2+2z_{-3}z_3&\quad\to\quad -w^2,
\end{align}\label{facto2-}
\end{subequations}
and replace the superscript ${}^\diamond$ with ${}^\sph$.

\subsection{Factorizations of the Gegenbauer equation}
\label{symcoma}
The factorizations of $\Delta_5^\sph$ of Subsect.
\ref{Factorizations of the Laplacian2} yield the following factorizations of the Gegenbauer operator:

\begin{align*}\cS_{\alpha,\lambda}
=&\partial_w\Big((1-w^2)\partial_w
  -2\alpha w\Big)\\
&+\Big(\alpha+\lambda-\frac12\Big)\Big(-\alpha+\lambda+\frac12\Big)\\
  =&\Big((1-w^2)\partial_w
  -2(1+\alpha)w\Big)\partial_w\\
  &+\Big(\alpha+\lambda+\frac12\Big)\Big(-\alpha+\lambda-\frac12\Big),\\
(1{-}w^2)\cS_{\alpha,\lambda}=&\Big((1-w^2)\partial_w
-\big(\alpha+\lambda-\frac12\big)w\Big)\Big((1-w^2)\partial_w
-\big(\alpha-\lambda+\frac12\big)w\Big)\\
&+\Big(\alpha+\lambda-\frac12\Big)\Big(\alpha-\lambda+\frac12\Big)\\
=&\Big((1-w^2)\partial_w
-\big(\alpha-\lambda-\frac12\big)w\Big)\Big((1-w^2)\partial_w
-\big(\alpha+\lambda+\frac12\big)w\Big)\\
&+\Big(\alpha+\lambda+\frac12\Big)\Big(\alpha-\lambda-\frac12\Big),
\end{align*}
\begin{align*}w^2\cS_{\alpha,\lambda}=&\Big(w(1-w^2)\partial_w
-\alpha-\lambda-\frac32+\big(-\alpha+\lambda-\frac12\big)w^2\Big)
\Big(w\partial_w+\alpha+\lambda+\frac12\Big)\\
&+\Big(\alpha+\lambda+\frac12\Big)\Big(\alpha+\lambda+\frac32\Big)\\
=&
\Big(w\partial_w+\alpha+\lambda-\frac32\Big)
\Big(w(1-w^2)\partial_w
-\alpha-\lambda+\frac12+\big(-\alpha+\lambda-\frac12\big)w^2\Big)\\
&+\Big(\alpha+\lambda-\frac12\Big)\Big(\alpha+\lambda-\frac32\Big)\\
=&\Big(w(1-w^2)\partial_w
-\alpha+\lambda-\frac32+\big(-\alpha-\lambda-\frac12\big)w^2\Big)
\Big(w\partial_w+\alpha-\lambda+\frac12\Big)\\
&+\Big(\alpha-\lambda+\frac12\Big)\Big(\alpha-\lambda+\frac32\Big)\\
=&
\Big(w\partial_w+\alpha-\lambda-\frac32\Big)
\Big(w(1-w^2)\partial_w
-\alpha+\lambda+\frac12+\big(-\alpha-\lambda-\frac12\big)w^2\Big)\\
&+\Big(\alpha-\lambda-\frac12\Big)\Big(\alpha-\lambda-\frac32\Big).
\end{align*}

\subsection{Standard solutions}

As usual, by standard solutions we mean solutions with a simple behavior around singular points. The singular points of the Gegenbauer equation are $\{1,-1,\infty\}$. The discussion of the point $-1$ can be easily reduced to that of $1$. Therefore, it is enough to discuss $2\times2=4$ solutions corresponding to two indices at $1$ and $\infty$.

The standard solutions can be expressed in terms of  the function
\begin{align}\notag
S_{\alpha,\lambda}(w)=S(a,b;w)&:=F\Big({a},{b};\frac{{a}+{b}+1}{2};\frac{1-w}{2}\Big)\label{qqe}\\
&=
F\Big(\frac{{a}}{2},\frac{{b}}{2};\frac{{a}+{b}+1}{2};1-w^2\Big).
\end{align}

Here are the 4 standard solutions. We consistently use the Lie-algebraic parameters.

\begin{alignat*}{2}
    \text{  $\sim1$  at $1$:}&\qquad S_{\alpha ,\lambda }(w)&\\
  =&F_{\alpha ,\alpha ,2\lambda }
\Big(\frac{1-w}{2}\Big)
=&
F_{\alpha ,-\frac12,\lambda }(1-w^2),
\end{alignat*}
\begin{alignat*}{2}
&\text{  $\sim \frac{1}{2^{\alpha}(1-w)^{\alpha }}$ at $1$:}\qquad
(1-w^2)^{-\alpha }S_{-\alpha ,-\lambda }(w)\\=&2^{-\alpha}(1-w)^{-\alpha }F_{-\alpha ,\alpha ,-2\lambda }
\Big(\frac{1-w}{2}\Big)=
(1-w^2)^{-\alpha }F_{-\alpha ,-\frac12,-\lambda }(1-w^2),\end{alignat*}
\begin{alignat*}{2}
&  \text{  $\sim  w^{-{a}}$ at $\infty$:}\qquad
  (w^2-1)^{\frac{-1-2\alpha +2\lambda }{4}}S_{-\lambda ,-\alpha }\Big(\frac{w}{\sqrt{w^2-1}}\Big)\\
=&(1+w)^{-\frac12-\alpha +\lambda }F_{-2\lambda ,\alpha ,-\alpha }
\Big(\frac{2}{1+w}\Big)=
w^{-\frac12-\alpha +\lambda }F_{-\lambda ,\alpha ,\frac12}(w^{-2}),
\end{alignat*}
\begin{alignat*}{2}
&\text{  $\sim  w^{-{b}}$ at $\infty$:}\qquad(w^2-1)^{\frac{-1-2\alpha -2\lambda }{4}}S_{\lambda ,\alpha }\Big(\frac{w}{\sqrt{w^2-1}}\Big)\\
=&(1+w)^{-\frac12-\alpha -\lambda }F_{2\lambda ,\alpha ,\alpha }
\Big(\frac{2}{1+w}\Big)=
w^{-\frac12-\alpha -\lambda }F_{\lambda ,\alpha ,\frac12}(w^{-2}).
\end{alignat*}

\subsection{Recurrence relations}

We will use the following normalization to express recurrence relations:
\begin{align}
  {\bf S}_{\alpha,\lambda}(w)&:=\frac{1}{\Gamma(\alpha+1)}S_{\alpha,\lambda}(w)\notag\\
&=\frac{1}{\Gamma(\frac{a+b+1}{2})}
F\Big({a},{b};\frac{{a}+{b}+1}{2};\frac{1-w}{2}\Big)\notag\\&=
{\bf F}_{\alpha,\alpha,2\lambda}\Big(\frac{1-w}{2}\Big).
\end{align}

To each root of $\so(5)$ there corresponds a recurrence relation:
\begin{align*}
&\p_w {\bf S}_{\alpha ,\lambda }(w)=-\12\Big(\frac12+\alpha -\lambda \Big)\Big(\frac12+\alpha +\lambda \Big)
 {\bf S}_{\alpha +1,\lambda }(w),\\
&\left((1-w^2)\p_w -2\alpha w\right)
{\bf S}_{\alpha ,\lambda }(w)
=-2{\bf S}_{\alpha -1,\lambda }(w) ,\\[2ex]
&\left((1-w^2)\p_w -\Big(\frac12+\alpha +\lambda \Big) w\right){\bf S}_{\alpha ,\lambda }(w)
=-\Big(\frac12+\alpha +\lambda \Big) {\bf S}_{\alpha ,\lambda +1}(w),
\\
&\left((1-w^2)\p_w -\Big(\frac12+\alpha -\lambda \Big) w\right){\bf S}_{\alpha ,\lambda }(w)
=-\Big(\frac12+\alpha -\lambda \Big){\bf S}_{\alpha ,\lambda -1}(w);
\\[2ex]
&\left(w\p_w+\frac12+\alpha -\lambda \right){\bf S}_{\alpha ,\lambda }(w)
=\frac12\Big(\frac12+\alpha -\lambda \Big)\Big(\frac32+\alpha-\lambda\Big) {\bf S}_{\alpha +1,\lambda -1}(w),\\
&\left(w(1{-}w^2)\p_w{+}\Big(\frac12{-}\alpha {+}\lambda \Big)(1{-}w^2){-}2\alpha w^2
\right){\bf S}_{\alpha ,\lambda }(w)
=-2{\bf S}_{\alpha -1,\lambda +1}(w),\\[2ex]
&\left(w\p_w+\frac12+\alpha +\lambda \right){\bf S}_{\alpha ,\lambda }(w)
=\frac12\Big(\frac12+\alpha +\lambda \Big)\Big(\frac32+\alpha+\lambda\Big) {\bf S}_{\alpha +1,\lambda +1}(w),\\
&\left(w(1{-}w^2)\p_w{+}\Big(\frac12{-}\alpha {-}\lambda \Big)(1{-}w^2){-}2\alpha  w^2
\right){\bf S}_{\alpha ,\lambda }(w)
=-2{\bf S}_{\alpha -1,\lambda -1}(w).
\end{align*}


\subsection{Wave packets in 5 dimensions}

We easily check the following lemma:
\bel For any $\tau$, the function
$z_2^\alpha\big(\sqrt2 z_0-\tau^{-1} z_{-3}+\tau z_3\big)^\nu$ is harmonic.
\eel

Let us make a wave packet from the above functions.
\bep Let
the contour $]0,1[\ni s\overset{\gamma}\mapsto\tau(s)$
    satisfy
\beq
\big(\sqrt2 z_0-\tau^{-1} z_{-3}+\tau z_3\big)^\nu\tau^{-\lambda}\Big|_{\tau(0)}^{\tau(1)}=0.\eeq
Then the function
  \[
  K_{\alpha,\nu,\lambda}(z_0,z_{-2},z_{2},z_{-3},z_3):=\int_\gamma z_2^\alpha\big(\sqrt2 z_0-\tau^{-1} z_{-3}+\tau z_3\big)^\nu\tau^{-\lambda-1}\d\tau
\]
is harmonic and
\bes  \begin{align}
    N_2K_{\alpha,\nu,\lambda}&=\alpha K_{\alpha,\nu,\lambda},\label{caz1}\\
    N_3K_{\alpha,\nu,\lambda}&=\lambda K_{\alpha,\nu,\lambda}.\label{caz2}
  \end{align}\ees
\eep
\proof (\ref{caz1}) is obvious.
To obtain (\ref{caz2}) we use Prop. \ref{cac0}. \qed

If in addition
\[\nu=-\alpha-\frac12,\]
then $K_{\alpha,\nu,\lambda}$ is homogeneous of degree $-\frac12$. Therefore, we can reduce it to dimension $3$.
Let us express it in the coordinates $w,r,p,u_2,u_3$:
\begin{align*}
    K
    (w,r,p,u_2,u_3)&=\int u_2^\alpha r^\alpha
    \Big(wr\sqrt2 +\frac{p}{\tau u_3\sqrt 2}+\frac{\tau pu_3}{\sqrt2}\Big)^{-\alpha-\frac12}\tau^{-\lambda-1}\d\tau\\
    &=(\sqrt2)^{\alpha+\frac12}u_2^\alpha u_3^\lambda r^{-\frac12}
   \int
    \Big(2w\sigma+(1+\sigma^2)\frac{p}{r}\Big)^{-\alpha-\frac12}
    \sigma^{\alpha-\lambda-\frac12}\d\sigma,
    \end{align*}
 where we set $\sigma:=u_3\tau$.
 Noting that on the spherical section $p=r$, we see that
 \beq S(w):=\int
    \big(2w\sigma+1+\sigma^2\big)^{-\alpha-\frac12}
    \sigma^{\alpha-\lambda-\frac12}\d\sigma
    \eeq
    satisfies the Gegenbauer equation.

\subsection{Integral representations}

In this subsection we describe two kinds of integral representations for solutions to the Gegenbauer equation. The first is essentially inherited
from the ${}_2\cF_1$ equation. The second was derived using additional variables in the previous subsection. Here we give  independent derivations.
We will use classical parameters.
\bet\ben\item[a)]
Let $[0,1]\ni \tau\overset{\gamma}\mapsto t(\tau)$ satisfy
\[(t^2-1)^{\frac{{b}-{a}+1}{2}}(t-w)^{-{b}-1}\Big|_{t(0)}^{t(1)}=0.\]
Then
\beq
\cS({a},{b};w,\p_w)
\int_\gamma (t^2-1)^{\frac{{b}-{a}-1}{2}}(t-w)^{-{b}}\d t=0
.\label{dad9}\eeq
\item[b)]
Let  $[0,1]\ni \tau\overset{\gamma}\mapsto t(\tau)$ satisfy
\begin{eqnarray*}&&
(t^2+2tw+1)^{\frac{-{b}-{a}}{2}+1}t^{b-2}
\Big|_{t(0)}^{t(1)}=0.\end{eqnarray*}
Then
\beq
\cS({a},{b};w,\p_w)
\int_\gamma(t^2+2tw+1)^{\frac{-{b}-{a}}{2}}t^{{b}-1}\d t=0
.\label{dad8}\eeq
\een
\label{gqw}\eet

\proof
For any contour $\gamma$ we have
\begin{eqnarray*}
\text{lhs of (\ref{dad9})}&=&{a}
\int_\gamma \d t\,\p_t(t^2-1)^{\frac{{b}-{a}+1}{2}}(t-w)^{-{b}-1},\\
\text{lhs of (\ref{dad8})}&=&\int\limits_\gamma\d t\,
\p_t(t^2+2tw+1)^{\frac{-{b}-{a}}{2}+1}t^{b-2}.\end{eqnarray*}
 \qed




Note that in the above theorem we can interchange $a$ and $b$. Thus we obtain four kinds of integral representations.

\subsection{Integral representations of the standard solutions}

As described in Thm \ref{gqw}, we have two types of integral representations of
solutions of Gegenbauer equations: a) and b).
It is natural to use singular points of the integrands as the endpoints of the contours of integration.
For the representations of type a) we have singular points at $\infty, -1,1, w$. For representations of type b) singular points are at $\infty,0$ and the two roots of $t^2+2tw+1=0$. Choosing an appropriate contour we obtain all standard solutions with both types of representations with some special normalizations. It is convenient to introduce special notation for these normalizations:

\begin{eqnarray}
{\bf S}_{\alpha,\lambda}^\I(w)&:=&2^{-\frac12-\alpha-\lambda}\frac{\Gamma(\frac{1+2\alpha+2\lambda}{2})\Gamma(\frac{1-2\lambda}{2})}{\Gamma(\alpha+1)}S_{\alpha,\lambda}
(w)
\\\notag
&=&
2^{-b}\frac{\Gamma(b)\Gamma(\frac{a-b+1}{2})}{\Gamma(\frac{a+b+1}{2})}
F\Big({a},{b};\frac{{a}+{b}+1}{2};\frac{1-w}{2}\Big)\\
\notag&=&2^{-\frac12-\alpha-\lambda}
{\bf F}_{\alpha,\alpha,2\lambda}^\I\Big(\frac{1-w}{2}\Big),\\[4ex]
{\bf S}_{\alpha,\lambda}^{\II}(w)&:=&\frac{\Gamma(\frac{1+2\alpha-2\lambda}{2})\Gamma(\frac{1+2\alpha+2\lambda}{2})}{\Gamma(2\alpha+1)}S_{\alpha,\lambda}(w)\\
&=&\notag
\frac{\Gamma(a)\Gamma(b)}{\Gamma(a+b)}
F\Big({a},{b};\frac{{a}+{b}+1}{2};\frac{1-w}{2}\Big),\\[4ex]
{\bf S}_{\alpha,\lambda}^{\0}(w)&:=&
\sqrt{\pi}
\frac{\Gamma(\frac{1+2\alpha}{2})}{\Gamma(\alpha+1)}S_{\alpha,\lambda}(w)
\\\notag
&=&
\sqrt{\pi}
\frac{\Gamma(\frac{1+2\alpha}{2})}{\Gamma(\alpha+1)}
F\Big({a},{b};\frac{{a}+{b}+1}{2};\frac{1-w}{2}\Big).
\end{eqnarray}

In the following table we list all standard solutions together with the contours of integration and the corresponding normalizations.

\medskip

\noindent
\begin{tabular}{cllc}&a)&b)& \\[3ex]
  $\sim 1$ at $1$:
  &$\begin{array}{l}]-\infty,-1],\\ \mathrm{I};
\end{array}$&$\begin{array}{l}[0,\infty[,\\ \mathrm{II};\end{array}$& \\[2ex]
\hspace{-2ex}$\sim \frac{1}{2^{\alpha}(1{-}w)^{\alpha }}$ at $1$:
&$\begin{array}{l}]{-}1,w],\\ \mathrm{I};\end{array}$&$\begin{array}{l}[{-}\i\sqrt{1{-}w^2}-w,\i\sqrt{1{-}w^2}{-}w],\\ 0;\end{array}$&\\[2ex]
      $\sim w^{-a}$ at $\infty$:
    &$\begin{array}{l}]-1,1],\\ 0;\end{array}$&$\begin{array}{l}[\sqrt{w^2-1}-w,0[,\\ \mathrm{I};\end{array}$& \\[2ex]
      $\sim w^{-b}$ at $\infty$:
      &$\begin{array}{l}]w,\infty],\\
\mathrm{II};\end{array}$&$\begin{array}{l}]-\infty,-\sqrt{w^2-1}-w],\\
\mathrm{I}.\end{array}$&\\
\end{tabular}

\medskip

Here are representations of type a):
\begin{align}\frac12>\Re\lambda>-\frac12-\Re\alpha:\quad&\\
\int\limits_{-\infty}^{-1}(t^2-1)^{-\frac12-\lambda }(w-t)^{-\frac12-\alpha +\lambda }\d
  t\nonumber \\
=
{\bf S}^\I_{\alpha ,\lambda }(w),&\quad w\not\in]-\infty,-1];\notag\end{align}
\begin{align}
\frac12>\Re\lambda>-\frac12+\Re\alpha:\quad&\\  
\int\limits_{w}^1(1-t^2)^{-\frac12-\lambda }(w-t)^{-\frac12-\alpha +\lambda }\d
  t\nonumber \\
=
  (1-w^2)^{-\alpha } {\bf S}^\I_{-\alpha ,-\lambda }(w),&\quad w\not\in]-\infty,-1]\cup[1,\infty[;\notag\end{align}
\begin{align}  \frac12>\Re\lambda:\quad\\
\int\limits_{-1}^1(1-t^2)^{-\frac12-\lambda }(w-t)^{-\frac12-\alpha +\lambda }\d
  t\nonumber \\
=
(w^2-1)^{\frac{-1-2\alpha +2\lambda }{4}}{\bf S}_{-\lambda ,\alpha }^{\0}\Big(\frac{w}{\sqrt{w^2-1}}\Big)
  ,&\quad w\not\in]-\infty,1];\notag\end{align}\begin{align}
  \Re\lambda+\frac12>|\Re\alpha|:\quad&\\
  \int\limits_w^\infty(t^2-1)^{-\frac12-\lambda }(t-w)^{-\frac12-\alpha +\lambda }\d
  t\nonumber &\\
=(w^2-1)^{\frac{-1-2\alpha -2\lambda }{4}}{\bf S}_{\lambda ,\alpha }^{\II}\Big(\frac{w}{\sqrt{w^2-1}}\Big),&\quad w\not\in]-\infty,1].\notag
  \end{align}

Next we list representations of type b):

  \begin{align}
\Re\alpha+\12>|\Re\lambda|:\quad&\\
\int\limits_0^{\infty}(t^2+2tw+1)^{-\alpha-\frac12}t^{-\frac12+\alpha +\lambda }\d t&\notag\\
\nonumber=
         {\bf S}_{\alpha ,\lambda }^{\II}(w)&\quad w\not\in]-\infty,-1];\notag\end{align}\begin{align}
         \12>\Re\alpha:\quad&
 \\
\int\limits_{-\i\sqrt{1-w^2}-w}^{\i\sqrt{1-w^2}-w}
(t^2+2tw+1)^{-\alpha -\frac12}(-t)^{-\frac12+\alpha +\lambda }\d t&
\nonumber\\
=
\i(1-w^2)^{-\alpha } 
    {\bf S}_{-\alpha ,-\lambda }^{\0}(w),&\quad w\not\in]-\infty,-1]\cup[1,\infty[;\notag\end{align}\begin{align}
    -\Re\lambda+\12>-\Re\alpha>-\12:\quad&\\
\int\limits_{\sqrt{w^2-1}-w}^{0}
(t^2+2tw+1)^{-\alpha-\frac12}(-t)^{-\frac12+\alpha -\lambda }\d t&\notag\\
=(w^2-1)^{\frac{-1-2\alpha +2\lambda }{4}}{\bf S}_{-\lambda ,\alpha }^\I\Big(\frac{w}{\sqrt{w^2-1}}\Big)  ,&\quad w\not\in]-\infty,1];\notag
        \end{align}\begin{align}
\Re\lambda+\12>-\Re\alpha>-\12:\quad&\\
\int\limits_{-\infty}^{-\sqrt{w^2-1}-w}
(t^2+2tw+1)^{-\alpha -\frac12}{(-t)}^{-\frac12+\alpha -\lambda }\d t&
\notag\\
=(w^2-1)^{-\frac14-\frac{\alpha}{2} -\frac{\lambda}{2}}{\bf S}_{\lambda ,\alpha }^\I\Big(\frac{w}{\sqrt{w^2-1}}\Big),&\quad w\not\in]-\infty,1].\notag
  \end{align}

\section{The Schr\"odinger Lie algebra and the heat equation}
\label{s8}
\init
By the {\em heat equation} on $\rr^n\oplus\rr$ we mean the equation given by the {\em heat operator}
  \beq \cL_n:=\Delta_n+2\partial_t.\eeq
 This operator has a large
family of generalized symmetries, the so-called {\em Schr\"odinger Lie algebra} and {\em group}. They can be derived from conformal symmetries of the Laplace equation. In this section we describe this derivation.

 In order to be consistent with Sect. \ref{s5}, it is convenient to consider
 $\cL_{n{-}2}$ instead of $\cL_n$.
Then the starting point, just as in Sect. \ref{s5}, is the $n+2$-dimensional ambient space. The Schr\"odinger Lie algebra and group are naturally contained in the pseudo-orhogonal Lie algebra and group for $n+2$ 
dimensions. Then, as described in Sect. \ref{subsec-conf}, we descend to the (flat) $n$ dimensional space and the corresponding Laplacian $\Delta_n$.  We assume that our functions  depend on $y_m$ only through the factor  $\e^{y_m}$. The variable $y_{-m}$ is renamed to $t$ (the ``time''). The Schr\"odinger Lie algebra and group respects functions of that form. The Laplacian $\Delta_n$ on such fuctions becomes 
the heat operator $\cL_{n{-}2}$. From the generalized symmetries of $\Delta_n$ we obtain generalized symmetries of $\cL_{n{-}2}$.

\subsection{$\sch(n{-}2)$ as a subalgebra of $\so(n{+}2)$}

We consider again the space $\rr^{n+2}$ with the split scalar product.
A special role will be played by the operator
\[B_{m+1,m}=z_{-m-1}\p_{z_{m}}-z_{-m}\p_{z_{m+1}}\in\so(n+2).\]
We define the {\em Schr\"odinger Lie algebra} and the  {\em
  Schr\"odinger group} as the commutants (centralizers) of this element:
\begin{subequations}
  \begin{align}\sch
(n-2)&:=\{B\in \so(n+2)\ :\ [B,B_{m+1,m}]=0\},\\
\Sch(n-2)&:=\{{\alpha}\in {\rm O}(n+2)\ :\
    {\alpha}B_{m+1,m}=B_{m+1,m}\alpha\}.
\end{align}\end{subequations}


\subsection{Structure of  $\sch(n{-}2)$}

Let us describe the structure of
 $\sch(n{-}2)$.

We will use our usual notation for elements of $\so(n{+}2)$ and ${\rm O}(n{+}2)$.
In particular,
\[N_{m}=-z_{-m}\p_{z_{-m}}+z_{m}\p_{z_{m}},\ \ \ \ \
N_{m+1}=-z_{-m-1}\p_{z_{-m-1}}+z_{m+1}\p_{z_{m+1}}.\]

Define
\beq
M:=-N_{m}+N_{m+1}.\eeq
Note that $M$ belongs to $\sch(n{-}2)$ and commutes with $\so(n{-}2)$, which is naturally embedded in  $\sch(n{-}2)$.

The Lie algebra $\sch(n{-}2)$ is spanned by the following operators:
\ben
\item
$B_{m+1,m}$, which spans the center of
 $\sch(n{-}2)$.
\item
$B_{m,j}$, $B_{m+1,j}$, $|j|=1,\dots,m-1$, which have the following nonzero commutator:
\beq[B_{m,j},B_{m+1,-j}]=B_{m+1,m}.\label{cent}\eeq
\item $B_{m+1,-m}$, $B_{-m-1,m}$, $M$,
which have the usual commutation relations of ${\rm sl}(2)\simeq \so(3)$:
\bes \begin{align}
  [B_{m+1,-m},B_{-m-1,m}]&=M,\\
   [M,B_{m+1,-m}]&
=-2B_{m+1, -m},\\
  [M,B_{-m-1,m}]&
=2B_{-m-1, m}
.\end{align}\ees
\item
$B_{i,j}$, $|i|<|j|\leq m-1$, $N_i$,  $i=1,\dots,m-1,$ with the usual commutation
  relations of $\so(n{-}2)$.
\een

The span of (2) can be identified with $\rr^{n-2}\oplus \rr^{n-2}\simeq
\rr^2\otimes \rr^{n-2}$, which has a natural structure of a symplectic space.
The span of (1) and (2) is  the central extension of
the abelian algebra $\rr^2\otimes \rr^{n-2}$ by (\ref{cent}). Such a Lie algebra is
usually called the {\em Heisenberg Lie algebra over }
$\rr^2\otimes \rr^{n-2}$ and can  be denoted by
\beq \heis(2(n{-}2))=\rr\rtimes(\rr^2\otimes \rr^{n-2}).\eeq

Lie algebras ${\rm sl}(2)$ and  $\so(n{-}2)$ act
in the obvious way on $\rr^2$, resp.
$\rr^{n-2}$. Thus ${\rm sl}(2)\oplus \so(n-2)$ acts on  $\rr^2\otimes \rr^{n-2}$. Thus
\begin{eqnarray}
\sch(n{-}2)&\simeq&
\rr\rrtimes( \rr^2\otimes\rr^{n-2})\rtimes \left({\rm sl}(2)\oplus \so(n{-}2)\right).
\end{eqnarray}

Note, in particular, that
 $\sch(n{-}2)$ is not semisimple.

The
subalgebra
spanned by the usual Cartan algebra of $\so(n{-}2)$, $M$ and $B_{-m-1,m}$ is a maximal commutative subalgebra of $\sch(n{-}2)$. It  will be called
the {\em  Cartan  algebra} of $\sch(n{-}2)$.

Let us introduce  $\kappa\in \SO(n{+}2)$:
\beq \kappa(\dots,z_{-m},z_{m},z_{-m-1},z_{m+1}):
=(\dots,z_{-m-1},z_{m+1},-z_{-m},-z_{m}).\eeq
Note that $\kappa^4=\id$ and $\kappa\in\Sch(n-2)$.
On the level of functions
\beq\kappa  K(\dots,z_{-m},z_{m},z_{-m-1},z_{m+1}):
= K(\dots,-z_{-m-1},-z_{m+1},z_{-m},z_{m}).\eeq

The subgroup of $\Sch(n{-}2)$ generated by
the Weyl group of $ {\rm O}(n{-}2)$ and $\kappa$
 will  be called
the  {\em Weyl group} of
$\sch(n{-}2)$.

\subsection{$\sch(n{+}2)$ in $n$ dimensions}

Recall from Subsect. \ref{subsec-conf}
that using the decomposition $\rr^{n+2}=\rr^n\oplus\rr^2$
we obtain the representations
\bes\begin{eqnarray}
\so(n{+}2)\ni B&\mapsto &
B^{\fl,\eta},
\label{repi1}\\
{\rm O}(n{+}2)\ni {\alpha}&\mapsto&
{\alpha}^{\fl,\eta}
\label{repi2}\end{eqnarray}\ees
acting on functions on $\rr^n$.
The Laplacian $\Delta_{n+2}$ becomes the Laplacian  $\Delta_n$ and it satisfies
 the generalized symmetry
\begin{subequations}
 \begin{eqnarray}B^{\fl,\frac{-2-n}{2}}
\Delta_n&=&
\Delta_nB^{\fl,\frac{2-n}{2}},\ \ B\in
\so(n{+}2),\label{syme1}\\
{\alpha}^{\fl,\frac{-2-n}{2}}
\Delta_n&=&
\Delta_n{\alpha}^{\fl,\frac{2-n}{2}},\ \ {\alpha}\in {\rm O}
(n{+}2).\label{syme2}\end{eqnarray}
\label{subu} \end{subequations}
The operator $B_{m+1,m}$ becomes
\beq B_{m+1,m}^{\fl ,\eta}=\partial_{y_{m}}.\eeq
Therefore,  all elements of $\sch(n{-}2)$
in the representation (\ref{repi1}) and all elements of $\Sch(n{-}2)$
in the representation (\ref{repi2})
have the form
\bes\begin{align}
B^{\fl ,\eta}&=C+D\partial_{y_{m}},\label{invo1}\\
 {\alpha}^{\fl ,\eta}f(\dots,y_{-m},y_m)&=\beta f\big(\dots,y_{-m},y_{m}+d(\dots,y_{-m})\big),\label{invo2}\end{align}
\label{invoo}
\ees
 where $C$, $D$,  $\beta$, $d$, do not involve the variable $y_{m}$.

\subsection{$\sch(n{-}2)$ in $(n-2)+1$ dimensions}
\label{subsec-sch}

We consider now the space $\rr^{n-2}\oplus\rr$ with the generic
variables $(y,t)=(\dots,y_{m-1},t)$. Note that $t$ should be
understood as the new name for $y_{-m}$, and we keep the old names for
the first $n{-}2$ coordinates.

We define the map
$\theta: C^\infty(\rr^{n-2}\oplus\rr)\to  C^\infty (\rr^n)$
by setting
\beq
(\theta h)(\dots,y_{m-1},y_{-m},y_m):=h(\dots,y_{m-1},y_{-m})\e^{y_{m}}.
\label{the1}\eeq
We also define $\zeta: C^\infty (\rr^n)\to  C^\infty (\rr^{n-2}\oplus\rr)$
\beq(\zeta f)(\dots,y_{m-1},t):=f(\dots,y_{m-1},t,0).\label{the2}\eeq
Clearly,
$\zeta$ is a left inverse of $\theta$:
\beq\zeta\circ\theta=\id.\eeq
Therefore, $\theta\circ\zeta=\id$ is true on the range of $\theta$.

The heat operator  in $n-2$ spatial
 dimensions
can be
obtained from the Laplacian in $n$ dimension:
\beq
\cL_{n-2}:=\Delta_{n{-}2}+2\p_{t}=\zeta
\Delta_n\theta.\label{conta}\eeq

For $B\in \sch(n-2)\subset \so(n+2)$
and ${\alpha}\in \Sch(n-2)\subset {\rm O}(n+2)$ we define
\bes\begin{align}
  B^{\sch,\eta}&:=\zeta B^{{\fl},\eta}\theta,\\
{\alpha}^{\sch,\eta}&:=\zeta {\alpha}^{\fl,\eta}\theta.
\end{align}\ees
It is easy to see, using
(\ref{invoo}), that
$\sch(n{-}2)$, $\Sch(n{-}2)$ and $\Delta_n$
 preserve the range of $\theta$. 
Therefore, for any $\eta$ we obtain representations
\bes\begin{eqnarray}\sch(n-2)\ni B&\mapsto& B^{\sch,\eta},\\
  \Sch(n-2)\ni {\alpha}&\mapsto &{\alpha}^{\sch,\eta}\end{eqnarray}\ees
acting on functions on $\rr^{n-2}\oplus\rr$.
By (\ref{subu}), we also have  generalized symmetries:
\bes\begin{eqnarray}
B^{\sch,\frac{-2-n}{2}}\cL_{n-2}&=&\cL_{n-2} B^{\sch,\frac{2-n}{2}} ,\ \ \ \ B\in \sch(n{-}2),\\
{\alpha}^{\sch,\frac{-2-n}{2}}\cL_{n-2}&=&\cL_{n-2} {\alpha}^{\sch,\frac{2-n}{2}} ,\ \ \ \ {\alpha}\in \Sch(n{-}2).\end{eqnarray}\ees
  
\subsection{Schr\"odinger symmetries in coordinates}

In this subsection we sum up information about Schr\"odinger
symmetries on
3 levels described in the previous subsections.

We start with generic
names of the variables and the corresponding squares:
\bes\begin{eqnarray}
z\in\rr^{n+2},& \langle z|z\rangle_{n+2}=&\sum_{|j|\leq m+1}z_{-j}z_j,\\
y\in\rr^n,&\langle y|y\rangle_n=&\sum_{|j|\leq m}y_{-j}y_j,\\
(y,t)\in\rr^{n-2}\oplus\rr,&\langle y|y\rangle_{n-2}=&\sum_{|j|\leq m-1}y_{-j}y_j.
\end{eqnarray}\ees

\noindent{\bf Cartan algebra of $\sch(n{-}2)$.}
 Central element

\bes\begin{eqnarray}
B_{m+1,m}&=&z_{-m-1}\p_{z_{m}}-z_{-m}\p_{z_{m+1}},\\
B_{m+1,m}^{{\fl}}&=&\p_{y_{m}},\\
B_{m+1,m}^{\sch}&=&1.
\end{eqnarray}\ees

\noindent Cartan algebra of $\so(n{-}2)$, $j=1,\dots,m-1$,
\bes\begin{eqnarray}
N_j&=&-z_{-j}\p_{z_{-j}}+z_j\p_{z_j},\\
N_j^\fl&=&-y_{-j}\p_{y_{-j}}+y_j\p_{y_j},\\
N_j^{\sch}&=&-y_{-j}\p_{y_{-j}}+y_j\p_{y_j}.
\end{eqnarray}\ees

\noindent Generator of scaling
\bes\begin{eqnarray}
M&=&z_{-m}\p_{z_{-m}}{-}z_{m}\p_{z_{m}}{-}
z_{-m-1}\p_{z_{-m-1}}{+}z_{m+1}\p_{z_{m+1}},\\
M^{\fl,\eta}&=&\sum\limits_{|j|\leq m-1}y_j\p_{y_j}+2y_{-m}\p_{y_{-m}}-\eta,\\
M^{\sch,\eta}&=&\sum\limits_{|j|\leq m-1}y_j\p_{y_j}+2t\p_t-\eta.
\end{eqnarray}\ees

\medskip

\noindent{\bf Root operators of $\sch(n{-}2)$.}
 Roots of $\so(n{-}2)$,  $|i|<|j|\leq m-1$,
\bes\begin{eqnarray}
B_{i,j}&=&z_{-i}\p_{z_j}-z_{-j}\p_{z_i},\\
B_{i,j}^{\fl}&=&y_{-i}\p_{y_j}-y_{-j}\p_{y_i},\\
B_{i,j}^{\sch}&=&y_{-i}\p_{y_j}-y_{-j}\p_{y_i}.\end{eqnarray}\ees

\noindent Space translations, $ |j|\leq m-1$,
\bes\begin{eqnarray}
B_{m+1,j}&=&z_{-m-1}\p_{z_j}-z_{-j}\p_{z_{m+1}},\\
B_{m+1,j}^{\fl}&=&\p_{y_j},\\
B_{m+1,j}^{\sch}&=&\p_{y_j}.\end{eqnarray}\ees

\noindent Time translation
\bes\begin{eqnarray}
B_{m+1,-m}&=&z_{-m-1}\p_{z_{-m}}-z_{m}\p_{z_{m+1}},\\
B_{m+1,-m}^{\fl}&=&\p_{y_{-m}},\\
B_{m+1,-m}^{\sch}&=&\p_t.
\end{eqnarray}
\ees

\noindent Additional roots, $|j|\leq m-1$,
\bes\begin{eqnarray}
B_{m,j}&=&z_{-m}\p_{z_j}-z_{-j}\p_{z_{m}},\\
B_{m,j}^{\fl}&=&y_{-m}\p_{y_j}-y_{-j}\p_{y_{m}},\\
B_{m,j}^{\sch}&=&t\p_{y_j}-y_{-j};\end{eqnarray}\ees

\bes\begin{eqnarray}
B_{-m-1,m}&=&z_{m+1}\p_{z_{m}}-z_{-m}\p_{z_{-m-1}},\\
B_{-m-1,m}^{\fl,\eta}&=&y_{-m}\big(\sum\limits_{|j|\leq m-1}y_j\p_{y_j}+y_{-m}\p_{y_{-m}}
-\eta \big)\notag
\\&&-\12 \sum_{|j|\leq m-1}y_{-j}y_j\p_{y_{m}},\\{}
B_{-m-1,m}^{\sch,\eta}&=&t\big(\sum\limits_{|j|\leq m-1}y_j\p_{y_j}+t\p_t
-\eta \big)\notag\\&&-\12 \sum_{|j|\leq m-1}y_{-j}y_j.
\end{eqnarray}\ees

\medskip

\noindent{\bf Weyl symmetries.}
We present a representative selection of elements of the Weyl group of $\Sch(n{-}2)$.
We will write $K$ for a function on $\rr^{n+2}$, $f$ for a function on
$\rr^n$, $h$ for a function on
$\rr^{n-2}\oplus \rr$ in the  coordinates $\big(\dots,y_{m-1},t\big)$.

\medskip

\noindent Reflection (for odd $n$)
\bes \begin{align}
  \tau_0 K(z_0,\dots,\dots,z_{-m},z_m,z_{-m-1},z_{m+1})\notag&\\
 & \hspace{-20ex}=
 K(-z_0,\dots,z_{-m},z_m,z_{-m-1},z_{m+1}),\\
 \tau_0^\fl f(y_0,\dots,y_{-m},y_m)\notag\\
 &\hspace{-10ex}=
 f(-y_0,\dots,y_{-m},y_m),\\
\tau_0^\sch h(y_0,\dots,t)&=
h(-y_0,\dots,t).
\end{align}\ees

\noindent Flips, $j=1,\dots,m-1$,
\bes\begin{align}
\tau_j K(\dots,z_{-j},z_j,\dots,z_{-m},z_m,z_{-m-1},z_{m+1})\notag&\\
&\hspace{-25ex}=
 K(\dots,z_{j},z_{-j},\dots,z_{-m},z_m,z_{-m-1},z_{m+1}),\\
\tau_j^\fl f(\dots,y_{-j},y_j,
\dots,y_{-m},y_m)\notag\\&\hspace{-15ex}=
 f(\dots,y_{j},y_{-j},\dots,y_{-m},y_m),\\
\tau_j^\sch h(\dots,y_{-j},y_j,\dots,t)&=
h(\dots,y_j,y_{-j},\dots,t).
\end{align}\ees

\noindent Permutations, $\pi\in S_{m-1}$,
\bes\begin{align}
\sigma_\pi K(\dots,z_{-m+1},z_{m-1},z_{-m},z_m,z_{-m-1},z_{m+1})&\notag\\
&\hspace{-35ex}=
 K(\dots,z_{-\pi_{m-1}},z_{\pi_{m-1}},z_{-m},z_m,z_{-m-1},z_{m+1}),\\
 \sigma_\pi^\fl f(\dots,y_{-m+1},y_{m-1},y_{-m},y_m)&\notag\\
& \hspace{-15ex}=
 f(\dots,y_{-\pi_{m-1}},y_{\pi_{m-1}},y_{-m},y_m),\\
 \sigma_\pi^\sch h(\dots,y_{-m+1},y_{m-1},t)&\notag\\
& \hspace{-8ex}
=h(\dots,y_{-\pi_{m-1}},y_{\pi_{m-1}},t).
\end{align}\ees
\noindent Special transformation $\kappa$
\bes\begin{align}
\kappa   K(\dots,z_{m-1},z_{-m},z_{m},z_{-m-1},z_{m+1})&\notag\\
&\hspace{-30ex}=  K(\dots,z_{m-1},-z_{-m-1},-z_{m+1},z_{-m},z_{m}),\\
\kappa^{\fl,\eta}f(\dots,y_{m-1},y_{-m},y_m)&\notag\\
&\hspace{-28ex}= y_{-m}^\eta
f\Big(\dots,
\frac{y_{m-1}}{y_{-m}},-\frac{1}{y_{-m}},
\frac{1}{2y_{-m}}\sum_{|j|\leq m}y_{-j}y_j\Big),
\\
\kappa^{\sch,\eta}h(\dots,y_{m-1},t)&\notag\\
&\hspace{-26ex}= t^\eta\exp\Big(\frac{1}{2t}\sum_{|j|\leq m-1}y_{-j}y_j\Big)
h\Big(\dots,\frac{y_{m-1}}{t},-\frac{1}{t}\Big).
\end{align}\ees
Square of $\kappa$
\bes\begin{align}
\kappa^2   K(\dots,z_{m-1},z_{-m},z_{m},z_{-m-1},z_{m+1})
&\notag\\&\hspace{-30ex}= K(\dots,z_{m-1},-z_{-m},-z_{m},-z_{-m-1},-z_{m+1}),\\
(\kappa^{\fl,\eta})^2  f(\dots,y_{m-1},y_{-m},y_{m})&\notag\\
&\hspace{-18ex}=f(\dots,-y_{m-1},y_{-m},y_{m}),\\
(\kappa^{\sch,\eta})^2h(\dots,y_{m-1},t)&=
h(\dots,-y_{m-1},t).
\end{align}\ees

\medskip

\noindent{\bf Laplacian/Laplacian / Heat operator}
\bes\begin{eqnarray}
\Delta_{n+2}
&=&\sum\limits_{|j|\leq m+1}\p_{z_{-j}}\p_{z_j},\\
\Delta_{n}
&=&\sum\limits_{|j|\leq m}\p_{y_{-j}}\p_{y_j},\\
\cL_{n-2}&=&\sum\limits_{|j|\leq m-1}\p_{y_{-j}}\p_{y_j}+2\p_t.
\end{eqnarray}\ees

\subsection{Special solutions of the heat equation}
\label{Special solutions of the heat equation}

Let us describe how to obtain solutions of the heat equation from solutions of the Laplace equation.

Consider first  a function on the level of $\rr^{n+2}$
\begin{align}   K(z)&=z_{-m}^{1-\frac{n}{2}}g\Big(\frac{z_{1}}{z_{-m}},\dots,\frac{z_{m-1}}{z_{-m}}\Big)\exp\Big(-\frac{z_{m+1}}{z_{-m}}\Big),
\end{align}
where $g$ is a harmonic function on $\rr^{n-2}$.
It is easy to see that $K$ is harmonic and satisfies
\begin{align}
  B_{m+1,m}K&=K.
\end{align}

Besides, $K$ is homogeneous of degree $1-\frac{n}{2}$. Therefore, we can descend on the level of dimension $n$, obtaining the function
\begin{align}   k(y)&=y_{-m}^{1-\frac{n}{2}}g\Big(\frac{y_{1}}{y_{-m}},\dots,\frac{y_{m-1}}{y_{-m}}\Big)  \exp\Big(\sum_{|i|\leq m-1}\frac{y_{-i}y_i}{y_{-m}}+y_{m} \Big).
\end{align}
It is harmonic and satisfies
\begin{align}
  B_{m+1,m}^\fl k&=k.
\end{align}  

Descending on the level of $\rr^{n-2}\oplus\rr$ we obtain
\beq
h(y,t)=
t^{1-\frac{n}{2}}g\Big(\frac{y_{1}}{t},\dots,\frac{y_{m-1}}{t}\Big)  \exp\Big(\sum_{|i|\leq m-1}\frac{y_{-i}y_i}{t} \Big).
\eeq
which solves the heat equation:
\beq
\cL_{n-2}h=0.
\eeq

\subsection{Wave packets for the heat equation}

Let us use the coordinates $(y,t)\in\rr^{n-2}\oplus\rr$.
    Recall that
    \begin{align}
      M^{\sch,\eta}&=\sum\limits_{|j|\leq m-1}y_j\p_{y_j}+2t\p_t-\eta.
      \end{align}
    The following proposition is proven by analogous arguments as Prop.
    \ref{cac0}. It allows us to form wave packets that are eigenfunctions of $M$:

    \bep
    Suppose that
 $]0,1[\ni s\overset{\gamma}\mapsto \tau(s)$ is a contour satisfying
    \beq
    f(\tau y,\tau^2 t)\tau^{-\nu}\Big|_{\tau(0)}^{\tau(1)}=0.\eeq
    Set
    \beq
    h_\nu(y,t):=\int_\gamma f(\tau y,\tau^2 t)\tau^{-1-\nu}\d\tau.\eeq
    Then
     \beq
    M^{\sch,\eta}
     h_\nu=(\nu-\eta) h_\nu.\eeq
\label{analo}     \eep

\section{Heat equation in 2 dimensions and the confluent equation}
\label{s9}
\init

The goal of this section is to derive the ${}_1\cF_1$  equation together with its symmetries from the heat equation in $2$ dimensions, which in turn comes from the Laplace equation in $6$ and $4$ dimensions.  Let us describe the main steps of this derivation:
\ben
\item\label{it1h}
  We start from the Schr\"odinger Lie algebra  $\sch(2)$ and group $\Sch(2)$
  considered as a subalgebra  of $\so(6)$, resp. a subgroup of $\mathrm{O}(6)$, acting in $6$ dimensions. The main initial operator is the Laplacian $\Delta_6$.
\item\label{it3h} We descend onto $4$ dimensions. The 6-dimensional Laplacian $\Delta_6$ becomes  the 4-dimensional Laplacian $\Delta_4$. 
\item We assume that the variable $y_2$ appears only in the exponential $\e^{y_2}$ and the variable $y_{-2}$ is renamed $t$. The Laplacian $\Delta_4$ becomes the heat operator $\cL_2$.
The  representations $B^{\sch,\eta}$ and $\alpha^{\sch,\eta}$ preserve our class of functions. With $\eta=-1$ and $\eta=-3$ they are generalized symmetries of the heat operator.
\item\label{it4h} We choose coordinates $w,s,u_1$, so that the Cartan operators are expressed in terms of  $s$, $u_1$. We compute
  $\cL_2$, $B^{\sch,\eta}$, and  $\alpha^{\sch,\eta}$ in the new coordinates.
\item\label{it5h} We make an ansatz that diagonalizes the Cartan  operators, whose eigenvalues,  denoted by $-\theta$ and $\alpha$,
  become parameters. The operators
  $\cL_2$,  $B^{\sch,\eta}$, and  $\alpha^{\sch,\eta}$ involve now only the single variable $w$. 
The operator  $\frac{s^2}{2}\cL_2$ becomes the  ${}_1\cF_1$  operator.
  Generalized symmetries of
  $\cL_2$ yield transmutation relations and discrete symmetries of the
 ${}_1\cF_1$  operator.
\een

The first part of this section is devoted to a description of the above steps, except for Step \ref{it3h}, discussed in detail in Sect. \ref{s8}.

The remaining part of this section is devoted to the theory of the ${}_1\cF_1$ equation 
and its solutions. Its organization is parallel
 to that of  Sect. \ref{s6} on the ${}_2\cF_1$ equation.
The main additional complication is the fact that besides
the ${}_1\cF_1$ equation
and the ${}_1F_1$ function,
it is useful to discuss  the closely related
${}_2\cF_0$ equation
and the  ${}_2F_0$ function.
In fact, some of the  standard solutions of the  ${}_1\cF_1$ equation are expressed in terms of
the  ${}_1F_1$ function, others in terms of the
  ${}_2F_0$ function.

\subsection{$\sch(2)$ in $6$ dimensions}

We again consider $\rr^6$ with the coordinates (\ref{sq0})
and the product given by  (\ref{sq1}):
\[
\langle z|z\rangle=2z_{-1}z_1+2z_{-2}z_2+2z_{-3}z_3.\]
We describe various object related to the Lie algebra
$\sch(2)$ treated as a subalgebra of $\so(6)$. We also list
 a few typical Weyl symmetries of $\Sch(2)$.

\noindent{\bf Lie algebra $\sch(2)$.}
Cartan algebra
\bes\begin{eqnarray}
  M&=&z_{-2}\p_{z_{-2}}-z_{2}\p_{z_{2}}-z_{-3}\p_{z_{-3}}+z_{3}\p_{z_{3}},\\
  N_1&=&-z_{-1}\p_{z_{-1}}+z_1\p_{z_1},\\
B_{3,2}&=&z_{-3}\p_{z_2}-z_{-2}\p_{z_{3}}.
\end{eqnarray}\ees
Root operators
\bes\begin{eqnarray}
   B_{3,-1}&=&z_{-3}\p_{z_{-1}}-z_1\p_{z_{3}},\\
 B_{2,1}&=&z_{-2}\p_{z_1}-z_{-1}\p_{z_{2}},\\
 B_{3,1}&=&z_{-3}\p_{z_{1}}-z_{-1}\p_{z_{3}},\\
 B_{2,-1}&=&z_{-2}\p_{z_{-1}}-z_1\p_{z_2},\\
 B_{3,-2}&=&z_{-3}\p_{z_{-2}}-z_2\p_{z_{3}},\\
  B_{-3,2}&=&z_{3}\p_{z_{2}}-z_{-2}\p_{z_{-3}}.
\end{eqnarray}\ees

\noindent{\bf Weyl symmetries}
\bes\begin{align}\id
K(z_{-1},z_1,z_{-2},z_2,z_{-3},z_3)&=K(z_{-1},z_1,z_{-2},z_2,z_{-3},z_3),\\
\tau_1
K(z_{-1},z_1,z_{-2},z_2,z_{-3},z_3)&=K(z_{1},z_{-1},z_{-2},z_2,z_{-3},z_3),\\
\kappa
K(z_{-1},z_1,z_{-2},z_2,z_{-3},z_3)&=
K(z_{-1},z_1,-z_{-3},-z_{3},z_{-2},z_{2}),\\
\tau_1\kappa
K(z_{-1},z_1,z_{-2},z_2,z_{-3},z_3)&=
K(z_1,z_{-1},-z_{-3},-z_{3},z_{-2},z_{2}).
\end{align}\ees

\noindent{\bf Laplacian}
\beq
\Delta_6=2\partial_{z_{-1}}\partial_{z_{1}}+
2\partial_{z_{-2}}\partial_{z_{2}}+
2\partial_{z_{-3}}\partial_{z_{3}}
.\label{sq2-}\eeq


\subsection{$\sch(2)$ in $4$ dimensions}

We descend on the level of $\rr^4$, with the coordinates $(y_{-1},y_1,y_{-2},y_2)$
  and the scalar product given by
  \[
\langle y|y\rangle=2y_{-1}y_1+2y_{-2}y_2.\]

\noindent{\bf Lie algebra $\sch(2)$.}
Cartan algebra
\bes\begin{eqnarray*}
M^{\fl ,\eta}&=&y_{-1}\p_{y_{-1}}+y_{1}\p_{y_{1}}+2y_{-2}\p_{y_{-2}}-
\eta,\\
N_1^\fl &=&-y_{-1}\p_{y_-1}+y_{1}\p_{y_{1}},\\
B_{3,2}^\fl&=&\p_{y_2}.
\end{eqnarray*}\ees

\noindent Root operators
\bes\begin{eqnarray*}
   B_{3,-1}^{\fl}&=&\p_{y_{-1}},\\
   B_{2,1}^{\fl}&=&y_{-2}\p_{y_1}-y_1\p_{y_{2}},\\
    B_{3,1}^{\fl}&=&\p_{y_{1}},\\
    B_{2,-1}^{\fl}&=&y_{-2}\p_{y_{-1}}-y_1\p_{y_2},\\
    B_{3,-2}^{\fl}&=&\p_{y_{-2}},\\
         B_{-3,2}^{\fl,\eta}&=&
y_{-2}(y_{-1}\p_{y_{-1}}+y_{1}\p_{y_{1}}+y_{-2}\p_{y_{-2}}-\eta)-y_{-1}y_1\p_{y_2}.
\end{eqnarray*}\ees

\noindent{\bf Weyl symmetries}
\bes\begin{align*}\id
  f(y_{-1},y_1,y_{-2},y_2)&=f(y_{-1},y_1,y_{-2},y_2),\\
  \tau_1^\fl f(y_{-1},y_{1},y_{-2},y_2)&=f(y_1,y_{-1},y_{-2},y_2),\\
  \kappa^{\fl,\eta}f(y_{-1},y_{1},y_{-2},y_2)
&=y_{-2}^\eta
f\Big(\frac{y_{-1}}{y_{-2}},\frac{y_1}{y_{-2}},-\frac1{y_{-2}},
\frac{y_{-1}y_1+y_{-2}y_2}{y_{-2}}\Big),\\
  \tau_1\kappa^{\fl,\eta}f(y_{-1},y_{1},y_{-2},y_2)
&=y_{-2}^\eta
f\Big(\frac{y_{1}}{y_{-2}},\frac{y_{-1}}{y_{-2}},-\frac1{y_{-2}},
\frac{y_{-1}y_1+y_{-2}y_2}{y_{-2}}\Big).
\end{align*}
\ees

\subsection{$\sch(2)$ in $2+1$ dimensions}

We apply the ansatz involving the exponential $\e^{y_2}$.
We rename $y_{-2}$ to $t$.

\noindent{\bf Lie algebra $\sch(2)$.}
Cartan algebra
\begin{subequations}
\begin{eqnarray}
M^{\sch,\eta}&=&y_{-1}\p_{y_{-1}}+y_{1}\p_{y_{1}}+2t\p_{t}-
\eta,\label{heat1}\\
N_1^{\sch}&=&-y_{-1}\p_{y_{-1}}+y_{1}\p_{y_{1}}
,\\B_{32}^\sch&=&1.\label{heat2}
\end{eqnarray}
\end{subequations}

\noindent
Root operators
\bes\begin{eqnarray}
  B_{3,-1}^{\sch}&=&\p_{y_{-1}},\\
   B_{2,1}^{\sch}&=&t\p_{y_1}-y_{-1},\\
  B_{3,1}^{\sch}&=&\p_{y_{1}},\\
B_{2,-1}^{\sch}&=&t\p_{y_{-1}}-y_1,\\
B_{3,-2}^{\sch}&=&\p_{t},\\
B_{-3,2}^{\sch,\eta}&=&
 t(y_{-1}\p_{y_{-1}}+y_{1}\p_{y_{1}}+t\p_{t}-\eta)-y_{-1}y_1.
\end{eqnarray}\ees

\noindent
    {\bf Weyl symmetries}
\begin{subequations}
    \begin{align}
      \id g(y_{-1},y_1,t)&=g(y_{-1},y_1,t),\\
\tau_1^\sch h(y_{-1},y_{1},t)&=h(y_1,y_{-1},t),\\
\kappa^{\sch,\eta}h(y_{-1},y_1,t)
&=t^\eta\exp\Big(\frac{y_{-1}y_1}{t}\Big)
h\Big(\frac{y_{-1}}{t},\frac{y_1}{t},-\frac1{t}\Big),\\
\tau_1\kappa^{\sch,\eta}h(y_{-1},y_1,t)
&=t^\eta\exp\Big(\frac{y_{-1}y_1}{t}\Big)
h\Big(\frac{y_{1}}{t},\frac{y_{-1}}{t},-\frac1{t}\Big).
\end{align}
\end{subequations}

\noindent{\bf Heat operator}
\beq
\cL_2=2\p_{y_{-1}}\p_{y_1}+2\p_{t}.\label{heat0}\eeq

\subsection{$\sch(2)$ in the coordinates $w,s,u_1$}

We introduce new coordinates $w,s,u_1$
  \begin{align}
w =\frac{y_{-1}y_1}{t}\;,\ \ \ \ 
&u_1 =\frac{y_{1}}{\sqrt{t}} \;,\ \ \ \
s =\sqrt{t}\,,\phantom{\frac{1}{2}}\label{cor}\\
  \intertext{with the reverse transformations}
y_{-1} =\frac{ sw}{u_1} \;,\ \ \ \ 
&y_{1} =u_1s\;,\ \ \ \
t =s^2\,.\phantom{\frac{1}{2}}
\end{align}

\noindent{\bf Lie algebra $\sch(2)$.} Cartan algebra
\begin{subequations}
\begin{align*}
  M^{\sch,{\eta}}& =s\dds-\eta,\\
  N_1^\sch& =u_1\p_{u_1}
  ,\\B_{32}^\sch&=1.
\end{align*}
\end{subequations}

\noindent Root operators
\begin{subequations}
  \begin{align*}
        B_{3,-1}^\sch& =
    \frac{u_1}{s}\p_w,\\
    B_{2,1}^\sch& =\frac{s}{u_1}(w\p_w+u_1\p_{u_1}-w),\\
    B_{3,1}^\sch& =\frac{1}{{u_1}s}(w\p_w+{u_1}\p_{u_1}),\\
  B_{2,-1}^\sch& =
  s{u_1}(\p_w-1),\\
  B_{3,-2}^{\sch}& =\frac{1}{s^2}\big(-w\p_w-\frac12{u_1}\p_{u_1}+\frac12s\p_s\big),\\
 B_{-3,2}^{\sch,\eta}&=s^2\big(w\p_w+\frac12{u_1}\p_{u_1}+\frac12s\p_s-w-\eta\big).
\end{align*}
\end{subequations}

\noindent{\bf Weyl symmetries}

\bes\begin{align*}
  \id h(w,{u_1},s)&=h(w,{u_1},s),\\
  \tau_1^\sch h\Big(w,{u_1},s\Big)&=h\Big(w,\frac{w}{{u_1}},s\Big),\\
\kappa^{\sch,\eta}h(w,{u_1},s)&
=s^{2\eta}\e^w
h\Big(-w,-\i {u_1},\frac{\i}s\Big),\\
\tau_1\kappa^{\sch,\eta}h(w,{u_1},s)&
=s^{2\eta}\e^w
h\Big(-w,-\frac{\i w}{{u_1}},\frac{\i}s\Big).
\end{align*}\ees

\noindent{\bf Heat operator}
\beq
\label{con:1}
\cL_{2}
=\frac{2}{s^2}\Big(w\p_w^2+({u_1}\p_{u_1}+1-w)\p_w+\frac12(-{u_1}\p_{u_1}+s\p_s)\Big).
\eeq

\subsection{Confluent operator}

Let us make the ansatz
\beq h(w, {u_1},s)={u_1}^\alpha s^{-\theta-1} F(w).\label{ansatz}\eeq
Clearly,
\bes\begin{eqnarray}
  M^{\sch,-1} h&=&-\theta h,\\
  N_1^{\sch} h&=&\alpha h,\\
  u_1^{-\alpha}s^{\theta+1}\frac{s^2}{2}\cL_2  h&=&
  \cF_{\theta ,\alpha}(w,\p_w)F(w),
  \end{eqnarray}\ees
where we have introduced the {\em ${}_1\cF_1$ operator}
\begin{eqnarray}
\cF_{\theta ,\alpha}(w,\p_w)
&=&w\p_w^2+(1+\alpha-w)\p_w-\frac{1}{2}(1+\theta +\alpha)
.\end{eqnarray}

Let us also define the closely related {\em ${}_2\cF_0$ operator}
\begin{eqnarray}
\tilde\cF_{\theta,\alpha}(w,\p_w)&=&
w^2\p_w^2+(-1+(2+\theta)w)\p_w+\frac14(1+\theta)^2-\frac14\alpha^2.\label{closely}
\end{eqnarray}
It is equivalent to the ${}_1\cF_1$ operator.  In fact, if $z=-w^{-1}$, then
\beq
(-z)^{\frac{3+\alpha+\theta}{2}}\tilde\cF_{\theta,\alpha}(z,\p_z)(-z)^{-\frac{1+\alpha+\theta}{2}}
=\cF_{\theta,\alpha}(w,\p_w).
\label{g8a}\eeq
We will treat $\cF_{\theta ,\alpha}(w,\p_w)$ as the principal operator.

Traditionally, one uses the {\em classical parameters} $a,b,c$:
\bes\begin{align}
&\alpha:=c-1=a-b,&& \theta: =-c+2a=-1+a+b;&\\
&a=\frac{1+\alpha+\theta }{2},&& b=\frac{1 -\alpha+\theta}{2},\ \ \ &
 c=1+\alpha.
\end{align}\ees
Here are the traditional forms of the ${}_1\cF_1$ and ${}_2\cF_0$ operators:
\begin{align}
\cF(a;c;w,\p_w)&:=w\p_w^2+(c-w)\p_w-a,
\label{f1c}\\
 {\cF}(a,b;-;w,\partial_w)&:=w^2\partial_w^2 +( -1+(1+a+b)w)\partial_w +ab. \label{g8}\end{align}

\subsection{Transmutation relations and discrete symmetries}

The heat operator satisfies the following generalized symmetries:
\bes\begin{align}
  B^{\sch,-3}\cL_2&=\cL_2 B^{\sch,-1},\quad B\in \sch(2),\label{simmi1}\\
  \alpha^{\sch,-3}\cL_2&=\cL_2 \alpha^{\sch,-1},\quad \alpha\in \Sch(2).
\label{simmi2}  \end{align}\ees

Applying (\ref{simmi1}) to the roots of $\sch(2)$ we obtain the following transmutation relations
of the confluent operator:
\[\begin{array}{rrl}
&\p_w&\cF_{\theta ,\alpha}\\[0ex]
&=\ \ \ \cF_{\theta +1,\alpha+1}&\p_w,\\[0.6ex]
&(w\p_w+\alpha-w)&\cF_{\theta ,\alpha}\\[0ex]
&=\ \ \ \cF_{\theta -1,\alpha-1}&(w\p_w+\alpha-w),\\[0.6ex]
&(w\p_w+\alpha)&
\cF_{\theta ,\alpha}\\[0ex]
&=\ \ \ \cF_{\theta +1,\alpha-1}&(w\p_w+\alpha),\\[0.6ex]
& (\p_w-1)&\cF_{\theta ,\alpha},\\[0ex]
&=\ \ \ \cF_{\theta -1,\alpha+1}&(\p_w-1);\\[0.6ex]
&\big( w\p_w+\12(\theta + \alpha+1)\big)&w\cF_{\theta ,\alpha}\\[0ex]
&=\ \ \  w\cF_{\theta +2,\alpha}&\big( w\p_w+\12(\theta + \alpha+1)\big),\\[0.6ex]
&\big(w\p_w+\12(-\theta +\alpha+1)-w)
&w\cF_{\theta ,\alpha}\\[0ex]
&=\ \ \ w\cF_{\theta -2,\alpha}&\big(w\p_w+\12(-\theta +\alpha+1)-w\big).
\end{array}\]

Applying (\ref{simmi2})
 to the Weyl symmetries  of $\sch(2)$ yields discrete symmetries of the confluent operator, described below.

The following operators equal $\cF_{\theta ,\alpha}(w,\p_w)$ for the appropriate $w$:
\[\begin{array}{rrcl}
w=v:
&&\cF_{\theta ,\alpha}(v,\p_v),&\\[0.2ex]
&v^{-\alpha}&\cF_{\theta ,-\alpha}(v,\p_v)&v^{\alpha},\\
w=-v:
&-\e^{-v}&\cF_{-\theta ,\alpha}(v,\p_v)&\e^v,
\\[0.2ex]
&-\e^{-v}v^{-\alpha}&\cF_{-\theta ,-\alpha}(v,\p_v)&\e^vv^{\alpha}
.\end{array}\label{newnot1}\]
The third symmetry is sometimes called the {\em 1st Kummer transformation}.

\subsection{Factorizations of of the heat operator}
\label{Factorizations of of the heat operator}

Special role is played by three distinguished subalgebras in $\sch(2)$: two isomorphic to $\heis(2)$  and one isomorphic to $\so(3)$.

First note
 the commutation relations
\beq
    [B_{2,-1},B_{3,1}] = [B_{2,1}\,,B_{3,-1}] = B_{3,2}.\label{heis:1}\eeq
    Therefore, the following subalgebras  in $\sch(5)$ are isomorphic to $\heis(2)$:
\bes\begin{eqnarray}
  \heis_-(2)& \hbox{spanned by}&B_{2,-1},\ B_{3,1},\ B_{3,2},\\
  \heis_+(2)& \hbox{spanned by}&B_{2,1},\ B_{3,-1},\ B_{3,2}.
\end{eqnarray}\ees

Note that the flip of $(1,-1)$, denoted $\tau_1$,
belongs to $\Sch(5)$ and  satisfies
\beq\tau_1B_{2,-1}\tau_1=B_{2,1},\quad
\tau_1B_{3,1}\tau_1=B_{3,-1},\quad \tau_1B_{3,2}\tau_1=B_{3,2}.\label{heis:1-}\eeq
Hence, 
\beq \tau_1\heis_-(2)\tau_1=\heis_+(2).\eeq

Let us define
\begin{subequations}
\begin{align}
\cC_{-}& =2\,B_{2,-1} B_{3,1}+M+N_1-B_{3,2}\,\\
& =2\,B_{3,1} B_{2,-1}+M+N_1+B_{3,2},\\
\cC_{+}& =2\,B_{2,1} B_{3,-1}+M-N_1-B_{3,2}\,\\
      & =2\,B_{3,-1} B_{2,1}+M-N_1+B_{3,2}.
\end{align}\label{dire1}
\end{subequations}
$\cC_+$ and $\cC_-$ can be viewed as the Casimir operators
for $\heis_+(2)$, resp. for $\heis_-(2)$. Indeed, $\cC_+$, resp. $\cC_-$
commute with all operators in $\heis_+(2)$, resp. $\heis_-(2)$.
We also have
\beq\tau_1\cC_-\tau_1=\cC_+.\eeq

On  the level of $\rr^2\oplus\rr$, the two operators $\cC_+$ and $\cC_-$ coincide. Indeed, a
direct calculation yields
\begin{eqnarray}
\cC_{+}^{\sch,\eta}=\cC_{-}^{\sch,\eta}=2t(\p_{y_{-1}}\p_{y_1}+\p_{t})
-\eta-1.
\label{dire}\end{eqnarray}

Second,
note
 the commutation relations
\beq
    [B_{-3,2},B_{3,-2}] = N_2-N_3=-M.\label{commu}
    \eeq
Therefore, the following of $\sch(2)$ is isomorphic to $\so(3)$:
\beq \so_{23}(3)\quad\text{  spanned by }\quad
B_{-3,2}, B_{3,-2}, M.\eeq

The Casimir operator
for  $\so_{23}(3)$ is
\begin{subequations}
\begin{align}
\cC_{23}& =4\,B_{3,-2} B_{-3,2}-(M-1)^2+1\,\\
         & =4\,B_{-3,2} B_{3,-2}-(M+1)^2+1.
\end{align}\label{dire2}
\end{subequations}
By (\ref{deq3a}) we have
\beq (2z_{-2}z_2+2z_{-3}z_3)\Delta_6^\diamond=-1+\cC_{23}^{\diamond,-1}+(N_1^{\diamond,-1})^2.
\label{dire.}\eeq
Inserting (\ref{dire2}) 
into (\ref{dire.}) we obtain
\begin{subequations}\begin{align}\notag
& (2z_{-2}z_2+2z_{-3}z_3)\Delta_6^\diamond\\=&4B_{2,-3}B_{-2,3}
  -(N_1+M+1)(-N_1+M+1)\\
  =&4B_{-2,3}B_{2,-3}
  -(N_1+M-1)(-N_1+M-1),
\end{align}\label{suu1}\end{subequations}
where the $B$, $N_1$ and $M$ operators should be equipped with the superscript ${}^{\diamond,-1}$.

Let us sum up the factorizations 
in  the variables $y_{-1}y_1,t$
obtained with the help of the three subalgebras:
\begin{subequations}
\begin{align}
t\cL_2
& =2\,B_{2,-1} B_{3,1}+M+N_1-1\,\label{su1}\\
& =2\,B_{3,1} B_{2,-1}+M+N_1+1\label{su2}\\
& =2\,B_{2,1} B_{3,-1}+M-N_1-1\,\label{su3}\\
      & =2\,B_{3,-1} B_{2,1}+M-N_1+1,\label{su4}\\
2y_{-1}y_1\cL_2&=-4B_{2,-3}B_{-2,3}
  -(N_1+M+1)(N_1-M-1)\label{su5}\\
  &=-4B_{-2,3}B_{2,-3}
  -(N_1+M-1)(N_1-M+1),\label{su6}
  \end{align}\end{subequations}
where the $B$, $N_1$ and $M$ operators should be equipped with the superscript ${}^{\sch,-1}$.

Indeed, to obtain (\ref{su1})--(\ref{su4}) we insert (\ref{dire1}) into  (\ref{dire}). To obtain (\ref{su5})--(\ref{su6}) we rewrite (\ref{suu1}),   multiplying it by $-1$.

In the variables $w,u,s$, we need to make the replacements 
\begin{subequations}
\begin{eqnarray}
y_{-1}y_1&\quad\to\quad ws^2,
\label{faca-1.}\\
t&\quad\to\quad s^2.\label{faca-2.}
\end{eqnarray}
\end{subequations}

\subsection{Factorizations of the confluent operator}

Factorizations of $\cL_2$ described in Subsect.
\ref{Factorizations of of the heat operator}
yield the following factorizations of the confluent operator:

\begin{eqnarray*}
\cF_{\theta,\alpha}
&=&\Big(\partial_w-1\Big)\Big(w\partial_w+\alpha\Big)-\frac12(\theta-\alpha+1)\\
&=&\Big(w\partial_w+1+\alpha\Big)\Big(\partial_w-1\Big)-\frac12(\theta-\alpha-1)\\
&=&\partial_w\Big(w\partial_w+\alpha-w\Big)-\frac12(\theta+\alpha-1)\\
&=&\Big(w\partial_w+1+\alpha-w\Big)\partial_w-\frac12(\theta+\alpha+1),
\end{eqnarray*}
\begin{eqnarray*}
w\cF_{\theta,\alpha}
&=&\Big(
w\partial_w+\frac12(-\theta+\alpha-1)-w\Big)
\Big(w\partial_w+\frac12(\theta+\alpha+1)\Big)
\\&&-\frac14(-\theta+\alpha-1)(\theta+\alpha+1)\\
&=&\Big(w\partial_w+\frac12(\theta+\alpha-1)\Big)\Big(
w\partial_w+\frac12(-\theta+\alpha+1)-w\Big)
\\&&-\frac14(-\theta+\alpha+1)(\theta+\alpha-1).
\end{eqnarray*}

\subsection{The ${}_1F_1$ function}

The ${}_1\cF_1$ equation (\ref{f1c})
 has a regular singular point at $0$.
Its indices at $0$ are equal to $0$, $1-c$.
For $c\neq 0,-1,-2,\dots$, the unique solution of the confluent equation analytic
at $0$  and equal to 1 at 0 is called 
the {\em ${}_1F_1$  function} or
{\em  Kummer's confluent function}. It is equal to
\[F(a;c;w):=\sum_{n=0}^\infty
\frac{(a)_n}{
(c)_n}\frac{w^n}{n!}.\]
It is defined for
$c\neq0,-1,-2,\dots$.
Sometimes it is more convenient to consider
the functions
\begin{align*} {\bf F}  (a;c;w)&:=\frac{F(a;c;w)}{\Gamma(c)}=
\sum_{n=0}^\infty
\frac{(a)_n}{
\Gamma(c+n)}\frac{w^n}{n!},\\
 {\bf F}^\I   (a;c;w)&:=\frac{\Gamma(a)\Gamma(c-a)}{\Gamma(c)}F(a;c;w).
\end{align*}


In the Lie-algebraic parameters:
\begin{eqnarray*} 
F_{\theta ,\alpha}(w)&:=&F\Bigl(\frac{1+\alpha+\theta }{2};1+\alpha;w\Bigr)
,\\
 {\bf F}  _{\theta ,\alpha}(w)&:=&
 {\bf F}  \Bigl(\frac{1+\alpha+\theta }{2};1+\alpha;w\Bigr)
=
 \frac{F_{\theta ,\alpha}(w)}{\Gamma(\alpha+1)},\\
 {\bf F}^\I   _{\theta ,\alpha}(w)&:=&
 {\bf F}^\I   \Bigl(\frac{1+\alpha+\theta }{2};1+\alpha;w\Bigr)=
\frac{\Gamma(\frac{1+\alpha+\theta}{2})
\Gamma(\frac{1+\alpha-\theta}{2})F_{\theta ,\alpha}(w)}{\Gamma(\alpha+1)}
.\end{eqnarray*}

\subsection{The ${}_2F_0$ function}
\label{The ${}_2F_0$ function}

Recall from (\ref{g8a}) that in parallel with the ${}_1\cF_{1}$ operator
it is useful to consider the  ${}_2\cF_0$ operator.
The  ${}_2\cF_0$ operator
 does not have a regular singular point at zero, hence  to construct its solutions having a simple behavior at zero we cannot use the Frobenius method. One of such solutions is the {\em ${}_2F_0$ function}. For $w\in\cc\backslash[0,+\infty[$ it can be defined  by
\[F(a,b;-;w):=\lim_{c\to\infty}F(a,b;c;cw),\]
where $|\arg c-\pi|<\pi-\epsilon$, $\epsilon>0$.
It extends to an analytic function on the universal cover of
$\cc\backslash\{0\}$ 
with a branch point of an infinite order at 0.
It has the following asymptotic expansion:
\[
F(a,b;-;w)\sim\sum_{n=0}^\infty\frac{(a)_n(b)_n}{n!}w^n,
\ |\arg w-\pi|<\pi-\epsilon.
\]
Sometimes instead of ${}_2F_0$ it is useful to consider the function
\begin{eqnarray*}
  F^\I   (a,b;-;w)&:=&\Gamma(a)F(a,b;-;w).
\end{eqnarray*}

When we use the Lie-algebraic parameters, we denote the ${}_2F_0$ function by
$\tilde F$ and $\tilde  F^\I$. The tilde is needed to  avoid the confusion with the ${}_1F_1$ functions:
\begin{eqnarray*}
\tilde F_{\theta ,\alpha}(w)&:=&F\Bigl(\frac{1+\alpha+\theta }{2},\frac{1-\alpha+\theta }{2};-;w\Bigr),\\
\tilde  F^\I  _{\theta ,\alpha}(w)&:=& F^\I  \Bigl(\frac{1+\alpha+\theta }{2},
\frac{1 -\alpha+\theta}{2};-;w\Bigr)
\,=\, \Gamma\Big(\frac{1-\alpha+\theta}{2}\Big)\tilde F_{\theta ,\alpha}(w)  .
\end{eqnarray*}

\subsection{Standard solutions}
\label{Standard solutions11}
The ${}_1F_1$ equation has two singular points. $0$ is a regular
singular point and with each of its two indices we can associate the
corresponding solution. $\infty$ is not a regular singular point. However  we can define two solutions with a simple behavior around $\infty$. Therefore, we obtain 4 {\em standard solutions}. 

The solutions that have a simple behavior at zero are expressed in terms of the function $F_{\theta,\alpha}$. Using $4$ discrete symmetries yields 4 distinct expressions. Taking into account  Kummer's identity we obtain $2$ pairs of standard solutions.

The solutions with a simple behavior at $\pm\infty$ are expressed in terms of $\tilde F_{\theta,\alpha}$. Again,  4 discrete symmetries yield 4 distinct expressions. Taking into account the trivial identity 
 $\tilde F_{\theta,\alpha}= \tilde F_{\theta,-\alpha}$ we obtain $2$ pairs of standard solutions.
\begin{align*}
  \text{  $\sim1$ at $0$}:\quad&\quad
  F_{\theta ,\alpha}(w)\\
  =&\e^wF_{-\theta ,\alpha}(-w);\\
\text{  $\sim w^{-\alpha}$ at $0$}:\quad&\quad
w^{-\alpha} F _{\theta ,-\alpha}(w)\\=&w^{-\alpha}
\e^w F _{-\theta ,-\alpha}(-w);\end{align*}\begin{align*}
\text{  $\sim w^{-a}$  at $+\infty$}:\quad&\quad
  w^{\frac{-1-\theta -\alpha}{2}}\tilde F_{\theta , \alpha}(-w^{-1})\\=& w^{\frac{-1-\theta -\alpha}{2}}\tilde F_{\theta , -\alpha}(-w^{-1});\\
\text{  $\sim  (-w)^{b-1}\e^w$ at $-\infty$}:\quad&\quad
\e^w(-w)^{\frac{-1+\theta -\alpha}{2}}\tilde F_{-\theta ,\alpha}(w^{-1})\\=&\e^w(-w)^{\frac{-1+\theta -\alpha}{2}}\tilde F_{-\theta ,-\alpha}(w^{-1})
.
\end{align*}
The solution $\sim w^{-a}$ at $+\infty$ is often called {\em Tricomi's confluent function}.

\subsection{Recurrence relations}
Recurrence relations for the confluent function
correspond to roots of the Lie algebra $\sch(2)$:
\begin{eqnarray*}
 \p_w {\bf F}  _{\theta ,\alpha}(w)&=&\frac{1+\theta +\alpha}{2} {\bf F}  _{\theta +1,\alpha+1}(w),
 \\
 \left(w\p_w+\alpha-w\right) {\bf F}  _{\theta ,\alpha}(w)&=& {\bf F}  _{\theta -1,\alpha-1}(w),\\[0.6ex]
 \left(w\p_w+\alpha\right) {\bf F}  _{\theta ,\alpha}(w)&=& {\bf
   F}  _{\theta +1,\alpha-1}(w),\\
 \left(\p_w-1\right) {\bf F}  _{\theta ,\alpha}(w)&=&\frac{-1+\theta -\alpha}{2} {\bf F}  _{\theta -1,\alpha+1}(w),
\\[0.6ex]
 \left(w\p_w+\frac{1+\theta +\alpha}{2}\right) {\bf F}  _{\theta ,\alpha}(w)&=&\frac{1+\theta +\alpha}{2} {\bf F}  _{\theta +2,\alpha}(w),
\\
 \left(w\p_w+\frac{1-\theta +\alpha}{2}-w\right) {\bf F}  _{\theta ,\alpha}(w)&
=&\frac{1-\theta +\alpha}{2} {\bf F}  _{\theta
  -2,\alpha}(w).
\end{eqnarray*}

\subsection{Wave packets for the heat equation in 2 dimensions}

Consider the space $\rr^2\oplus\rr$ and the heat equation given by the operator $\cL_2=2\partial_{y_{-1}}\partial_{y_1}+2\partial_t$.
  Recall that  \begin{align*}
  M^{\sch,-1}&= y_{-1}\p_{y_{-1}}+y_1\p_{y_1}+2t\p_t+1,\\
      N_1^\sch&=-y_{-1}\p_{y_{-1}}+y_1\p_{y_1}.\end{align*}

Set
\bes  \begin{align}\notag&
    G^a_{\theta,\alpha}(y_{-1},y_1,t)\\:=& \int_{\gamma^a}\tau^{-\alpha-1}t^{\frac{-1-\theta+\alpha}{2}}(\tau^{-1}y_{-1}-1)^{\frac{-1+\theta-\alpha}{2}}\exp\Big(\frac{(y_{-1}-\tau)y_1}{t}\Big)\d\tau,\\
    \notag   & G_{\theta,\alpha}^b(y_{-1},y_1,t)\\
    := &\int_{\gamma^b}\tau^{-\alpha-1}t^{\frac{-1-\theta-\alpha}{2}}(\tau y_1-1)^{\frac{-1+\theta+\alpha}{2}}\exp\Big(\frac{y_{-1}(y_1-\tau^{-1})}{t}\Big)\d\tau.  \end{align}\ees
(The superscripts $a$ and $b$ denote two kinds of wave packets, and not parameters $a$, $b$).

  \bep If the  contours  $\gamma^a$ and $\gamma^b$ are appropriately chosen, then
\begin{align} 
  \cL_2 G_{\theta,\alpha}^a&=0,\quad&\cL_2G^b_{\theta,\alpha}&=0,
  \label{ko1}\\
 M^{\sch,-1} G_{\theta,\alpha}^a&=-\theta
 G_{\theta,\alpha}^a,\quad& M^{\sch,-1} G_{\theta,\alpha}^b&=-\theta G_{\theta,\alpha}^b,
 \label{ko2}\\
    N_1 G_{\theta,\alpha}^a&=\alpha
    G_{\theta,\alpha}^a,\quad    & N_1 G_{\theta,\alpha}^b&=\alpha
    G_{\theta,\alpha}^b.\label{ko3}
    \end{align}
\eep

\proof
By the analysis of Subsect.
\ref{Special solutions of the heat equation},
the following functions
\bes\begin{align}
g_\nu^a(y_{-1},y_1,t)&:=
t^{-1-\nu}y_{-1}^\nu\exp\Big(\frac{y_{-1}y_1}{t}\Big),\label{sd1}\\
  g_\nu^b(y_{-1},y_1,t)&:=
  t^{-1-\nu}y_1^\nu\exp\Big(\frac{y_{-1}y_1}{t}\Big)\label{sd2}
  \end{align}\ees
solve the heat equation. They still solve the heat equation after translating and rotating. Therefore, 
\bes\begin{align}
  G_{\theta,\alpha}^a(y_{-1},y_1,t)
  &=\int_{\gamma^a} g_{\frac{-1+\theta-\alpha}{2}}^a(\tau^{-1}(y_{-1}-1),\tau y_1,t)\tau^{-\alpha-1}\d\tau,
  \label{ka3}\\
  G_{\theta,\alpha}^b(y_{-1},y_1,t)
  &=\int_{\gamma^b} g_{\frac{-1+\theta+\alpha}{2}}^b(\tau^{-1}y_{-1},\tau (y_1-1),t)\tau^{-\alpha-1}\d\tau
  \label{ka1}
\end{align}\label{kaa}\ees
also solve the heat equation. This proves (\ref{ko1}).

If the contours satisfy the requirements of Prop. \ref{cac0}, then   (\ref{kaa})
imply  (\ref{ko3}).

We can rewrite  (\ref{kaa}) in a somewhat different way:
\bes\begin{align}
(\ref{ka3})
  &=\int_{\gamma^a}  g_{\frac{-1+\theta-\alpha}{2}}^a(\tau^{-1}( y_{-1}-1),\tau^{-1} y_1,\tau^{-2}t)(\tau^{-1})^\theta\d(\tau^{-1}),
  \label{ka4}\\
(\ref{ka1})  &=\int_{\gamma^b}  g_{\frac{-1+\theta+\alpha}{2}}^b(\tau y_{-1},\tau (y_1-1),\tau^2t)\tau^\theta\d\tau.
  \label{ka2}
\end{align}\label{kab}\ees
If the contours satisfy the requirements of Prop. \ref{analo}, then  (\ref{kab}) imply (\ref{ko2}).
 \qed

Now we express the above wave packets in the coordinates $w,s,{u_1}$:
\bes\begin{align}
(\ref{ka4})&=
  \int s^{-1-\theta+\alpha}\Big(\frac{ws}{\tau {u_1}}-1\Big)^{\frac{-1+\theta-\alpha}{2}}\exp\Big(w-\frac{\tau {u_1}}{s}\Big)\tau^{-\alpha-1}\d\tau,\label{ki4}\\
   (\ref{ka2})&=\int s^{-1-\theta-\alpha}(\tau {u_1}s-1)^{\frac{-1+\theta+\alpha}{2}}\exp\Big(w\big(1-\frac{1}{\tau {u_1}s}\big)\Big)\tau^{-\alpha-1}\d\tau.
  \label{ki1}
\end{align}\ees

In (\ref{ki4}) we make the substitution $\sigma:=w-\frac{\tau {u_1}}{s}$, or $\tau=\frac{s}{{u_1}}(w-\sigma)$.
In (\ref{ki1}) we make the substitution $\sigma:=\frac{1}{1-\frac{1}{\tau {u_1}s}}$, or $\tau=\frac{\sigma}{{u_1}s(\sigma-1)}$.
We obtain
\bes\begin{align}
  G_{\theta,\alpha}^a(w,s,{u_1})
    &=
  s^{-1-\theta}{u_1}^\alpha
  F_{\theta,\alpha}^a(w),  \label{ki6}\\
    G_{\theta,\alpha}^b(w,s,{u_1})
  &= s^{-1-\theta}{u_1}^\alpha
  F_{\theta,\alpha}^b(w),
    \label{ki3}
\end{align}\ees
where
\bes\begin{align}
    F_{\theta,\alpha}^a(w)&:=\int_{\gamma^a}\sigma^{\frac{-\alpha+\theta-1}{2}}
  (w-\sigma)^{\frac{-\alpha-\theta-1}{2}}\e^\sigma\d\sigma,  \label{ki6.}\\ F_{\theta,\alpha}^b(w)&:=\int_{\gamma^b}\exp\Big(\frac{w}{\sigma}\Big)\sigma^{-\alpha-1}(\sigma-1)^{\frac{\alpha+\rho-1}{2}}\d\sigma.     \label{ki3.}
\end{align}\ees
The above analysis shows that (for appropriate contours) the functions
(\ref{ki6.}) and (\ref{ki3.}) satisfy the confluent equation.

\subsection{Integral representations}

Let us prove directly that integral  (\ref{ki6.}) and (\ref{ki3.}) solve the confluent equation.
\bet\ben \item[a)]
Let $[0,1]\ni \tau\overset{\gamma}\mapsto t(\tau)$ satisfy
$t^{a-c+1}\e^t(t-w)^{-a-1}\Big|_{t(0)}^{t(1)}=0.$
Then
\beq
\cF(a;c;w,\p_w)
\int_\gamma t^{a-c}\e^t(t-w)^{-a}\d t=0.\label{dad1}\eeq
\item[b)]
Let  $[0,1]\ni \tau\overset{\gamma}\mapsto t(\tau)$ satisfy
$\e^{\frac{w}{t}}t^{-c}(1-t)^{c-a}\Big|_{t(0)}^{t(1)}=0.$
Then
\beq
\cF(a;c;w,\p_w)\int_\gamma\e^{\frac{w}{t}}t^{-c}(1-t)^{c-a-1}\d t
=0.
\label{dad}
\eeq
\een\label{dad4}\eet

\proof
We check that for any contour $\gamma$
\begin{eqnarray*}
\text{lhs of (\ref{dad1})}&=&
-a\int_\gamma\d t\,\p_t
t^{a-c+1}\e^t(t-w)^{-a-1},\\
\text{lhs of (\ref{dad}) }&=&-\int_\gamma\d t\,\p_t
\e^{\frac{w}{t}}t^{-c}(1-t)^{c-a}.
\end{eqnarray*}
 \qed

\subsection{Integral representations of standard solutions}

Using the integral representations of type a)
and  attaching contours to 
$-\infty$, $0$ and $w$
 we can obtain all standard solutions.

 Similarly, using the integral representations of type b)
and  attaching contours to 
$0-0$, $1$ and $\infty$  we can obtain all standard solutions.

Here is the list of contours:

\medskip

\noindent
\begin{tabular}{lll}&a)&b)\\
    $\sim1$ at $0$:\qquad&
$]-\infty,(0,w)^+,-\infty[$,&$[1,+\infty[$;\\
        \text{  $\sim w^{-\alpha}$ at $0$}:\qquad&
             $[0,w]$,&$(0-0)^+$;\\
             \text{  $\sim w^{-a}$  at $+\infty$}:\qquad&
     $ ]-\infty,0]$,&$]-\infty,0]$;\\
\text{  $\sim  (-w)^{b-1}\e^w$ at $-\infty$}:\qquad&
$[w,-\infty[$,&$[0,1]$.\\
\end{tabular}

\medskip

($0,w)^+$ means that we bypass $0$ and $w$ counterclockwise. $(0-0)^+$ means that the contour departs from $0$ on the negative side, encircles it and then comes back again from the negative side.

Here are the explicit formulas for a)-type integral representations:
\begin{subequations}
\begin{align}\text{all $\theta,\alpha$:}&\label{sub1}\\
\underset{]-\infty,(0,w)^+-\infty[}{\frac{1}{2\pi \i}\int}
t^{\frac{-1+\theta -\alpha}{2}}\e^t(t-w)^{\frac{-1-\theta -\alpha}{2}}\d t&\notag\\
= {\bf F}  _{\theta ,\alpha}(w),\notag\end{align}
\begin{align}
\Re(1-\alpha)>|\Re \theta |:&\label{sub2}\\
\int_0^w 
t^{\frac{-1+\theta -\alpha}{2}}\e^t(w-t)^{\frac{-1-\theta -\alpha}{2}}\d t\notag
&\\
=
w^{-\alpha} {\bf F}^\I   _{\theta ,-\alpha}(w),&\quad w\not\in]-\infty,0];\notag\\
\int_w^0
(-t)^{\frac{-1+\theta -\alpha}{2}}\e^t(t-w)^{\frac{-1-\theta -\alpha}{2}}\d t\notag
&\\
=
(-w)^{-\alpha} {\bf F}^\I   _{\theta ,-\alpha}(w),&\quad w\not\in[0,\infty[;\notag\end{align}
    \begin{align}
\Re(1+\theta -\alpha)>0:&\label{sub3}\\
\int_{-\infty}^0
(-t)^{\frac{-1+\theta -\alpha}{2}}\e^t(w-t)^{\frac{-1-\theta -\alpha}{2}}\d t
&\notag\\=
w^{\frac{-1-\theta -\alpha}{2}} \tilde F^\I  _{\theta ,\alpha}(-w^{-1}),&\quad
 w\not\in]-\infty,0];\notag\end{align}\begin{align}xs
\Re(1-\theta -\alpha)>0:&\label{sub4}\\
\int_{-\infty}^w
(-t)^{\frac{-1+\theta -\alpha}{2}}\e^t(w-t)^{\frac{{-}1-\theta -\alpha}{2}}\d t
&\notag\\=
\e^w(-w)^{\frac{-1+\theta -\alpha}{2}} \tilde F^\I  _{-\theta ,\alpha}(w^{-1}),&\quad w\not\in[0,\infty[\notag.\end{align}
\end{subequations}

We also present explicit formulas for b)-type integral representations:
\begin{subequations}
\begin{align}
  \Re(1+\alpha)>|\Re \theta |:&\label{subi1}\\
\int\limits_{[1,+\infty[}\e^{\frac{w}{t}}t^{-1-\alpha}(t-1)^{\frac{-1-\theta +\alpha}{2}}
\d t&\notag\\
= {\bf F}^\I   _{\theta ,\alpha}(w);&
\notag\end{align}\begin{align}
\text{    all $\theta,\alpha$:}&\label{subi2}\\
\frac{1}{2\pi\i}
\int\limits_{(0-0)^+}
\e^{\frac{w}{t}}t^{-1-\alpha}(1-t)^{\frac{-1-\theta +\alpha}{2}}
\d t&\notag\\
=w^{-\alpha} {\bf F}  _{\theta ,-\alpha}(w),&\quad \Re w>0;\notag\end{align}\begin{align}
 \Re(1+\theta +\alpha)>0:&\label{subi3}\\
\int_{-\infty}^0
\e^{\frac{w}{t}}(-t)^{-1-\alpha}(1-t)^{\frac{-1{ -}\theta +\alpha}{2}}
\d t&\notag\\=
w^{\frac{-1-\theta -\alpha}{2}}
\tilde F^\I  _{\theta ,-\alpha}(-w^{-1}),&\quad \Re w>0;\notag\end{align}\begin{align}
\Re(1-\theta +\alpha)>0:&\label{subi4}\\
\int_0^1
\e^{\frac{w}{t}}t^{-1-\alpha}(1-t)^{\frac{-1{ -}\theta +\alpha}{2}}
\d t&\notag\\
=\e^w(-w)^{\frac{ -1+\theta -\alpha}{2}}
\tilde F^\I  _{-\theta ,-\alpha}(w^{-1}),&\quad \Re w<0.\notag
\end{align}
\end{subequations}

\subsection{Connection formulas}

The two solutions with a simple behavior at infinity
can be expressed as linear combination of the solutions with a simple behavior at zero:
\bes\begin{eqnarray}\label{pqd1}
w^{\frac{-1-\theta -\alpha}{2}}\tilde F_{\theta ,\pm \alpha}(-w^{-1})
&=&\frac{\pi{\bf F}  _{\theta ,\alpha}(w)}{\sin{\pi(-\alpha)}\Gamma\left(\frac{1+\theta -\alpha}{2}\right)}
 \\
&&+\frac{\pi  w^{-\alpha} {\bf F}  _{\theta ,-\alpha}(w)}{\sin\pi \alpha\Gamma\left(\frac{1+\theta +\alpha}{2}\right)}
,
\quad w\not\in]-\infty,0];\notag
 \\\label{pqd2}
\e^w(-w)^{\frac{-1 +\theta -\alpha}{2}}\tilde F_{-\theta ,\pm\alpha}(w^{-1})
&=&\frac{\pi {\bf F}  _{\theta ,\alpha}(w)}{\sin{\pi(-\alpha)}
\Gamma\left(\frac{1-\theta -\alpha}{2}\right)}
\\
&&+\frac{\pi(-w)^{-\alpha} {\bf F}_{\theta ,-\alpha}(w)}
{\sin\pi \alpha\Gamma\left(\frac{1-\theta +\alpha}{2}\right)}
,\quad w\not\in[0,+\infty[.\notag
\end{eqnarray}\label{pqd0}
\ees

Note that (\ref{pqd1}) uses a different domain from (\ref{pqd2}).  This is natural, however it is inconvenient when we want to  rewrite
(\ref{pqd0})  in the matrix form, because on the right hand side of (\ref{pqd1}) and  (\ref{pqd2}) the second standard solutions differ by a phase factor.

Let us introduce the matrix
\begin{align*}
A_{\theta,\alpha}&:=
\frac{\pi}{\sin(\pi\alpha)}\left[\begin{array}{cc}
  \frac{-1}{\Gamma\left(\frac{1+\theta -\alpha}{2}\right)}
& \frac{\e^{-\frac{\i\pi}{2}\alpha}}{\Gamma\left(\frac{1+\theta +\alpha}{2}\right)}
\\[2.5ex]
\frac{-1}{\Gamma\left(\frac{1-\theta -\alpha}{2}\right)}
&\frac{\e^{\frac{\i\pi}{2}\alpha}}{\Gamma\left(\frac{1-\theta +\alpha}{2}\right)}
  \end{array}\right],\end{align*}
satisfying
\begin{align}
  A_{\theta,\alpha}^{-1}&=
  \frac{\i\e^{\frac{\i\pi}{2}\theta}}{2}
  \left[\begin{array}{cc}
  \frac{\e^{\frac{\i\pi\alpha}{2}}}{\Gamma\left(\frac{1-\theta +\alpha}{2}\right)}
& \frac{-\e^{-\frac{\i\pi\alpha}{2}}}{\Gamma\left(\frac{1+\theta +\alpha}{2}\right)}
\\[2.5ex]
\frac{1}{\Gamma\left(\frac{1-\theta -\alpha}{2}\right)}
&\frac{-1}{\Gamma\left(\frac{1+\theta -\alpha}{2}\right)}
    \end{array}\right],\\
  \det A_{\theta,\alpha}&=  -\frac{\i\pi\e^{-\frac{\i\pi}{2}\theta}}{2\sin(\pi\alpha)}.
\end{align}
Then we have for $\Im w>0$
\begin{align}&
  \left[\begin{array}{c}
      w^{\frac{-1-\theta -\alpha}{2}}\tilde F_{\theta ,\pm \alpha}(-w^{-1})\\[2ex]
      \e^w(-w)^{\frac{-1+\theta -\alpha}{2}}\tilde F_{-\theta ,\pm\alpha}(w^{-1})
      \end{array}\right]
  =
A_{\theta,\alpha}
  \left[\begin{array}{c}
 {\bf F}  _{\theta ,\alpha}(w)\\[2.5ex]
 (-\i w)^{-\alpha} {\bf F}_{\theta ,-\alpha}(w)
 \end{array}\right].
\end{align}


Let us show how to derive connection
formulas from integral representations of  type a).
We have
\bes    \begin{align}\notag
    &    \Bigg(
    \int\limits_{]-\infty,0-\i 0]}+\int\limits_{[0-\i0,w]}
-    \int\limits_{]-\infty,0+\i 0]}-\int\limits_{[0+\i0,w]}\Bigg) t^{\frac{-1+\theta -\alpha}{2}}\e^t(t-w)^{\frac{-1-\theta -\alpha}{2}}\d t\\
    =&
    \underset{]-\infty,(0,w)^+,-\infty[}{\int}t^{\frac{-1+\theta -\alpha}{2}}\e^t(t-w)^{\frac{-1-\theta -\alpha}{2}}\d t ,\qquad w\not\in[-\infty,0[;
    \\[1ex]\notag
        &    \Bigg(
    \int\limits_{]-\infty,w-\i 0]}+\int\limits_{[w-\i0,0]}
-    \int\limits_{]-\infty,w+\i 0]}-\int\limits_{[w+\i0,0]}\Bigg) t^{\frac{-1+\theta -\alpha}{2}}\e^t(t-w)^{\frac{-1-\theta -\alpha}{2}}\d t\\
    =&
    \underset{]-\infty,(0,w)^+,-\infty[}{\int}t^{\frac{-1+\theta -\alpha}{2}}\e^t(t-w)^{\frac{-1-\theta -\alpha}{2}}\d t ,
\qquad w\not\in]0,+\infty]
.\end{align}\ees
    We obtain
\bes    \begin{align}\notag
      &-\sin(\pi\alpha) w^{\frac{-1-\theta -\alpha}{2}} \tilde F^\I  _{\theta ,\alpha}(-w^{-1})
    +
    \cos\frac{\pi(\theta+\alpha)}{2}w^{-\alpha} {\bf F}^\I   _{\theta ,-\alpha}(w)\\=&\pi  {\bf F}  _{\theta ,\alpha}(w),\qquad w\not\in[-\infty,0[
    ;\label{exw2}\\[1ex]\notag
  &  -\sin(\pi\alpha) \e^w(-w)^{\frac{-1+\theta -\alpha}{2}} \tilde F^\I  _{-\theta ,\alpha}(w^{-1})+
      \cos\frac{\pi(\theta-\alpha)}{2}(-w)^{-\alpha} {\bf F}^\I   _{\theta ,-\alpha}(w)\\=&\pi  {\bf F}  _{\theta ,\alpha}(w),\qquad w\not\in]0,+\infty].\label{exw1}
    \end{align}\ees
    This  implies (\ref{pqd0}).

\section{Heat equation in 1 dimension and the Hermite equation}
\label{s10}
\init

The goal of this section is to derive the Hermite  equation together with its symmetries from the heat equation in $1$ dimension, which in turn comes from the Laplace equation in $5$ and $3$ dimensions.

The first part of this section describes main steps of the derivation of the Hermite equation. They are parallel to those of the derivation of the ${}_1\cF_1$ equation:
\ben
\item\label{it1hh}
  We start from the Schr\"odinger Lie algebra $\sch(1)$ and group $\Sch(1)$
  considered as a subalgebra  of $\so(5)$, resp. a subgroup of $\mathrm{O}(5)$, acting in $5$ dimensions. The main initial operator is the Laplacian $\Delta_5$.
\item\label{it3hh} We descend onto $3$ dimensions.  The 5-dimensional Laplacian $\Delta_5$ becomes  the 3-dimensional Laplacian $\Delta_3$. 
\item  We descend on $1+1$ dimensions. The Laplacian $\Delta_3$ becomes the heat operator $\cL_1$.
The  representations $B^{\sch,\eta}$ and $\alpha^{\sch,\eta}$ with $\eta=-\frac12$ and $\eta=-\frac52$  are generalized symmetries of $\cL_1$.
\item\label{it4hh} We choose coordinates $w,s$, so that the Cartan operator is expressed in terms of  $s$. We compute
  $\cL_1$,  $B^{\sch,\eta}$ and  $\alpha^{\sch,\eta}$ in the new coordinates.
\item\label{it5hh} We make an ansatz that diagonalizes the Cartan operator, whose eigenvalue
  becomes a parameter,  denoted by $\lambda$. 
  $\cL_1$,  $B^{\sch,\eta}$ and  $\alpha^{\sch,\eta}$ involve now only the single variable $w$. 
  $2s^2\cL_1$ turns out to be the  Hermite operator.
  The generalized symmetries of
  $\cL_1$ yield transmutation relations and discrete symmetries of the
 Hermite  operator.
\een

(As in the previous section, in our presentation we omit the step \ref{it3hh}).

In the remaining part of this section
we develop the theory of the Hermite equation 
and its solutions. Its organization  is parallel to that of all other sections on individual equations, and especially of Sect. \ref{s7} on the Gegenbauer equation.
In particular, the Gegenbauer equation can be derived by a quadratic relation from the ${}_2\cF_1$ equation in essentially the same way as
 the Hermite equation can be derived  from the ${}_1\cF_1$ equation.

\subsection{$\sch(1)$ in 5 dimensions}

We again consider $\rr^5$ with the coordinates
\beq z_0,z_{-2},z_2,z_{-3},z_3\label{geg1-}\eeq
and the scalar product given by
\beq \langle z|z\rangle=z_0^2+2z_{-2}z_2+2z_{-3}z_3.\label{geg5-}\eeq
We
keep the notation from $\so(5)$---remember that $\sch(1)$ is a subalgebra of $\so(5)$.

 \noindent
     {\bf Lie algebra $\sch(1)$.}
The Cartan algebra
\bes\begin{eqnarray}
M&=&
z_{-2}\p_{z_{-2}}-z_{2}\p_{z_{2}}-z_{-3}\p_{z_{-3}}+z_{3}\p_{z_{3}},\\
B_{3,2}& =&z_{-3}\p_{z_2}-z_{-2}\p_{z_{3}}.
\end{eqnarray}\ees
Root operators
\bes\begin{eqnarray}
 B_{3,0}&=&z_{-3}\p_{z_0}-z_0\p_{z_{3}},\\
 B_{2,0}&=&z_{-2}\p_{z_0}-z_0\p_{z_2},\\
 B_{3,-2}&=&z_{-3}\p_{z_{-2}}-z_2\p_{z_{3}},\\
 B_{-3,2}&=&z_{3}\p_{z_2}-z_{-2}\p_{z_{-3}}.
\end{eqnarray}\ees

 \noindent{\bf Weyl symmetries}
\bes\begin{align} \id K(z_0,z_{-2},z_2,z_{-3},z_3)
&=K(z_0,z_{-2},z_2,z_{-3},z_3),\\
 \kappa K(z_0,z_{-2},z_2,z_{-3},z_3)
 &=K(z_0,-z_{-3},-z_{3},z_{-2},z_{2}),\\
  \kappa^2 K(z_0,z_{-2},z_2,z_{-3},z_3)
  &=K(z_0,-z_{-2},-z_2,-z_{-3},-z_3),\\
   \kappa^3 K(z_0,z_{-2},z_2,z_{-3},z_3)
   &=K(z_0,z_{-3},z_{3},-z_{-2},-z_{2}).
   \end{align}\ees

\noindent{\bf Laplacian}
\beq \Delta_5=\partial_{z_0}^2+2\partial_{z_{-2}}\partial_{z_{2}}+
2\partial_{z_{-3}}\partial_{z_{3}}
.\eeq

\subsection{$\sch(1)$ in 3 dimensions}

We descend on the level of $\rr^3$, 
with the variables $ y_0,y_{-2},y_2$
and the scalar product given by
\[ \langle y|y\rangle=y_0^2+2y_{-2}y_2.\label{geg5--}\]

 \noindent{\bf Lie algebra $\sch(1)$.}
Cartan algebra
\begin{eqnarray*}
M^{\fl,\eta}&=&y_{0}\p_{y_0}+2y_{-2}\p_{y_{-2}}-\eta,\\
B_{3,2}^\fl& =&\p_{y_2}.
\end{eqnarray*}

 \noindent Root operators
\begin{eqnarray*}
 B_{3,0}^{\fl}&=&\p_{y_0},\\
 B_{2,0}^{\fl}&=&y_{-2}\p_{y_0}-y_0\p_{y_2},\\
 B_{3,-2}^{\fl}&=&\p_{y_{-2}},\\
 B_{-3,2}^{\fl,\eta}
&=&y_{-2}\big(y_0\p_{y_0}+y_{-2}\p_{y_{-2}}-\eta\big)-\frac12y_0^2\p_{y_2}.
\end{eqnarray*}

 \noindent{\bf Weyl symmetries}
 \begin{align*}
   \id^{\fl,\eta}f(y_0,y_{-2},y_2)&=f(y_0,y_{-2},y_2),\\
\kappa^{\fl,\eta}f(y_0,y_{-2},y_2)&
=y_{-2}^\eta
f\Big(\frac{y_0}{y_{-2}},-\frac1{y_{-2}},\frac{y_0^2+2y_{-2}y_2}{2y_{-2}}\Big),\\
   \big(\kappa^{\fl,\eta}\big)^2f(y_0,y_{-2},y_2)&=(-1)^\eta f(-y_0,y_{-2},y_2),\\
\big(\kappa^{\fl,\eta}\big)^3f(y_0,y_{-2},y_2)&
=(-y_{-2})^\eta
f\Big(-\frac{y_0}{y_{-2}},-\frac1{y_{-2}},\frac{y_0^2+2y_{-2}y_2}{2y_{-2}}\Big).
\end{align*}

\noindent{\bf Laplacian}
\[\Delta_5^\fl=\partial_{y_0}^2+2\partial_{y_{-2}}\partial_{y_{2}}
.\]

\subsection{$\sch(1)$ in $1+1$ dimensions}

We descend onto the level of $\rr\oplus\rr$, as described in
Subsect. \ref{subsec-sch}.
 We rename $y_{-2}$ to $t$.

\noindent{\bf Lie algebra $\sch(1)$.} Cartan algebra:
\bes\begin{align}
  M^{\sch,\eta}&={y_0}\p_{y_0}+2t\p_t-\eta,\\
  B_{3,2}&=1.
  \end{align}\ees
Root operators
\bes\begin{eqnarray}
 B_{3,0}^{\sch}&=&\p_{y_0},\\
 B_{2,0}^{\sch}&=&t\p_{y_0}-{y_0},\\
 B_{3,-2}^{\sch}&=&\p_t,\\
 B_{-3,2}^{\sch,\eta}&=&t({y_0}\p_{y_0}+t\p_t-\eta)-\frac12y_0^2.\end{eqnarray}\ees

\noindent{\bf Weyl symmetry}
\bes\begin{align}
   \id^{\sch,\eta}h({y_0},t)&=
  h({y_0},t),\\
 \kappa^{\sch,\eta}h({y_0},t)&= t^{\eta}\exp(\tfrac{y_0^2}{2t})
 h(\tfrac{{y_0}}{t},-\tfrac{1}{t}),\\
 \big(\kappa^{\sch,\eta}\big)^2h({y_0},t)&=
 (-1)^\eta h(-{y_0},t),\\
  \big(\kappa^{\sch,\eta}\big)^3h({y_0},t)&= (-t)^{\eta}\exp(\tfrac{y_0^2}{2t})
  h(-\tfrac{{y_0}}{t},-\tfrac{1}{t}).
  \end{align}\ees

 \noindent{\bf Heat operator}
\beq\Delta_5^\sch=\cL_1=\partial_{y_0}^2+2\partial_t.\eeq

 \subsection{$\sch(1)$ in the coordinates $w,s$}

 Let us define new coordinates
\begin{align}\label{coo2}
w =\frac{{y_0}}{\sqrt{2\,t}}\;,\ \ \ \ \ \ \ \ \ 
&s =\sqrt{t}\,,\phantom{\frac{1}{2}}
\intertext{with the reverse transformation}
  {y_0} =\sqrt{2}\,s\,w\;,\ \ \ \ \ \ \ \ 
&t =s^2\,.\phantom{\frac{1}{2}}
\end{align}

 \noindent{\bf Lie algebra $\sch(1)$.} Cartan algebra
\begin{eqnarray*}
  M^{\sch,{\eta}}&=&s\dds-\eta,\\
B_{3\,2}& =&1.
\end{eqnarray*}
Root operators

\begin{align*}
  B_{3,0}^\sch& =\frac{1}{\sqrt{2}\,s}\,\ddw,\\
  B_{2,0}^\sch& =\frac{s}{\sqrt{2}}\,\left(\ddw-2w\right),\\
B_{3,-2}^\sch& =\frac{1}{2\,s^2}\,\left(-w\ddw+s\p_s\right),\\
B_{-3,2}^{\sch,{\eta}}& =\frac{s^2}{2}\,\left(w\ddw+s\p_s-2\eta-2w^2\right).
\end{align*}


 \noindent{\bf Weyl symmetries}
\begin{align*}
   \id^{\sch,{\eta}} h(w,s)&= h(w,s),\\
   \kappa^{\sch,{\eta}} h(w,s)&=s^{2\eta}\e^{w^2}h(\i w,-\tfrac{\i}{s}),\\
   (\kappa^{\sch,{\eta}})^2 h(w,s)&=(-1)^\eta h(-w,s),\\
   (\kappa^{\sch,{\eta}})^3 h(w,s)&=(-s^2)^\eta\e^{w^2}h(-\i w,-\tfrac{\i}{s})   .\end{align*}
 \noindent{\bf Heat operator}
\beq
\cL_1=\frac{1}{2s^2}\left(\ddw^2-2w\ddw+2s\dds\right)\,.\eeq

\subsection{Hermite operator}

Let us set $\eta=-\frac12$ and use
the ansatz
\beq h(w,s)=s^{-\lambda-\frac12} S(w).\label{ansatz1}\eeq
Clearly,
\begin{eqnarray}
  M^{\sch,-\frac12} h&=&-\lambda h,\\
  s^{\lambda+\frac12}2s^2\cL_1 h&=&\cS_\lambda(w,\p_w)S(w),
\end{eqnarray}
where we have introduced the  {\em Hermite operator}
\beq
S_\lambda(w,\p_w):=\p_w^2-2w\p_w-2\lambda-1.
\eeq
We will also use an alternative notation
\beq
S(a;w,\p_w):=\p_w^2-2w\p_w-2a,
\eeq
so that
\beq
\lambda={a}-\frac12,\ \ \ a=\lambda+\frac12.\eeq


\subsection{Quadratic transformation}
\label{Quadratic transformation-}

Let us go back to 2+1 dimensions and the heat operator
\beq\cL_2=2\partial_{y_{-1}}\partial_{y_{1}}+
2\p_t
.\label{sq20=}\eeq
Let us use the reduction described in Subsect.
\ref{Dimensional reduction}, and then applied in Subsect. \ref{Quadratic transformation}:
\beq y_0:=\sqrt{2y_{-1}y_{1}},\quad
u:=\sqrt{\frac{y_{1}}{y_{-1}}}.\label{coo1=}\eeq
  In the new variables,
  \begin{align}
        N_1&=u\p_u,\\
    \cL_2
    &=\Big(\partial_{y_0}+\frac1{2y_0}\Big)^2
    -\frac1{y_0^2}\Big(u\partial_u-\frac12\Big)
    \Big(u\partial_u+\frac12\Big)+
 2\p_t.
  \end{align}
  Therefore,
\bes  \begin{align}
  (uy_0)^{\frac12}\cL_2  (uy_0)^{-\frac12}
   & =
    -\frac1{y_0^2}N_1
    \Big(N_1- 1\Big)+
    \cL_1,\\
      (u^{-1}y_0)^{\frac12}\cL_2  (u^{-1}y_0)^{-\frac12}
   & =
    -\frac1{y_0^2}N_1
    \Big(N_1+1\Big)+
    \cL_1.
    \end{align}\ees
  Compare the coordinates (\ref{cor}) for 2+1 dimensions and
  the coordinates (\ref{coo2}) for 1+1 dimensions.
  The coordinate $s$ are the same.
 This is not the case of
  $w$, so let us rename  $w$ from (\ref{coo2})
  as $v$. We then have $w=v^2$.
  We also have
  \[uy_0=\sqrt2su_1,\quad u^{-1}y_0=\sqrt 2w u_1^{-1}s.\]
  Hence on functions that do not depend on $u$ we obtain
\bes  \begin{align}
    s^{\frac12}u_1^{\frac12}\cL_2 s^{-\frac12}u_1^{-\frac12}&=
    \cL_1 ,\\
        s^{\frac12}u_1^{-\frac12}v\cL_2 s^{-\frac12}u_1^{\frac12}v^{-1}&=
        \cL_1 .
        \end{align}
\ees

Thus by a quadratic
 transformation we can transform the Hermite equation into a special case of the confluent equation:
\bes\begin{eqnarray}
\cS_\lambda(v,\p_v)&=&4\cF_{\lambda,-\frac12}(w,\p_w),\\
v^{-1}\cS_\lambda(v,\p_v)v&=&4\cF_{\lambda,\frac12}(w,\p_w),
\end{eqnarray}\label{ada1}\ees
where
\[w=v^2,\ \ \ \ v=\sqrt w.\]

\subsection{Transmutation relations and discrete symmetries}
\label{symcom5}

The heat operator satisfies the
generalized symmetries
\bes\begin{eqnarray}
B^{\sch,-\frac52}\cL_1&=&\cL_1B^{\sch,-\frac12},\quad B\in \sch(1);\label{torr1}\\
\alpha^{\sch,-\frac52}\cL_1&=&\cL_1\alpha^{\sch,-\frac12},\quad
\alpha\in \Sch(1).\label{torr3}
\end{eqnarray}\ees

Equation (\ref{torr1}) applied to the roots of $\sch(1)$ implies the transmutation relations of the Hermite operator:
\[\begin{array}{rrll}
\p_w&\cS_\lambda &
=\ \ \cS_{\lambda +1}&\p_w,\\[0.3ex]
(\p_w-2w)
&\cS_\lambda &
=\ \ \ \cS_{\lambda -1}&(\p_w-2w),\\[0.6ex]
(w\p_w+\lambda +\frac12)&w^2\cS_\lambda&
=\ \ w^2\cS_{\lambda +2}& (w\p_w+\lambda +\12),\\[0.3ex]
(w\p_w-\lambda +\frac12-2w^2)&w^2\cS_\lambda &
=\ \ w^2\cS_{\lambda -2}&(w\p_w-\lambda +\12-2w^2).
\end{array}\]

Relation (\ref{torr1}) applied to the Weyl symmetries of $\sch(1)$ implies the discrete symmetries of the Hermite operator, described below.

The following operators equal $\cS_\lambda (w,\p_w)$ for an
appropriate $w$: 
\begin{subequations} \begin{align}
w=\pm v:\qquad
&&&\cS_\lambda (v,\p_v),\\
w=\pm\i v:\qquad
&&-\exp(-v^2)&\cS_{-\lambda }(v,\p_v)\exp(v^2).\label{siu2}
\end{align}
\end{subequations}

\subsection{Factorizations of the heat operator}
\label{Factorizations of the heat operator2}

Special role is played by two distinguished subalgebras of $\sch(2)$.

First note the commutation relations
\beq [B_{2,0},B_{3,0}]=B_{3,2}.\eeq
    Therefore, we have the following distinguished subalgebra  in $\sch(1)$ isomorphic to $\heis(2)$:
\begin{eqnarray}
  \heis_0(2)& \hbox{spanned by}&B_{2,0},\ B_{3,0},\ B_{3,2}.
\end{eqnarray}
Let us define
\bes\begin{align}
\cC_0& =2\,B_{2,0} B_{3,0}+2M-B_{3,2}\,\\
      & =2\,B_{3,0} B_{2,0}+2M+B_{3,2}.
\end{align}\ees

We have the commutation relations
\begin{align*}
    [\cC_0,B_{2,0}]&=-2B_{2,0}(B_{3,2}-1),\\
    [\cC_0,B_{3,0}]&=2B_{3,0}(B_{3,2}-1),\\
    [\cC_0,B_{3,2}]&=0.\end{align*}
But $B_{3,2}^{\sch,\eta}=1$. Therefore,
 on the level  of $\rr\oplus\rr$  the operator   $\cC_0^{\sch,\eta}$
     can be treated as a kind of a Casimir operator of $\heis_0(2)$:
it    commutes with all elements of
    $\heis_0(2)$.
Note the identity
\beq 2t\cL_1=\cC_{0}^{\sch,{-\tfrac{1}{2}}}.\eeq

Second,
consider
 $B_{-3,2}$, $B_{3,-2}$, $M$.  They are contained both in $\sch(6)$ and in $\sch(5)$. Therefore, the subalgebra $\so_{23}(3)$, described
 in Sect. \ref{Factorizations of of the heat operator} in the context of
 $\sch(6)$, is also contained in $\sch(5)$.
Recall that its 
 Casimir operator is
\begin{subequations}
\begin{align}
\cC_{23}& =4\,B_{3,-2} B_{-3,2}-(M+1)^2+1\,\\
         & =4\,B_{-3,2} B_{3,-2}-(M-1)^2+1.
\end{align}\label{dire2.}
\end{subequations}
By (\ref{deq3a}) we have
\beq (2z_{-2}z_2+2z_{-3}z_3)\Delta_5^\diamond
=\cC_{23}^{\diamond,{-\tfrac{1}{2}}}-\frac{3}{4}.
\label{dire..}\eeq
Inserting (\ref{dire2.}) 
into (\ref{dire..}) we obtain
\begin{subequations}\begin{align}
& (2z_{-2}z_2+2z_{-3}z_3)\Delta_5^\diamond\notag\\=&4B_{2,-3}B_{-2,3}
  -\Big(M+\frac32\Big)\Big(M+\frac12\Big),\\
  =&4B_{-2,3}B_{2,-3}
  -\Big(M-\frac32\Big)\Big(M-\frac12\Big),
\end{align}\label{suu1-}\end{subequations}
where the $B$, $N_1$ and $M$ operators should be decorated with the superscript ${}^{\diamond,-\frac12}$.

Let us sum up the factorizations 
in  the variables $y_0,t$
obtained with the help of the two subalgebras:
\begin{subequations}
\begin{align}
  2t\cL_1& =2\,B_{2,0} B_{3,0}-(-2M+1)\,\\
      & =2\,B_{3,0} B_{2,0}-(-2M-1),\\
  -y_0^2\cL_1&=4B_{2,-3}B_{-2,3}
  -\Big(M+\frac32\Big)\Big(M+\frac12\Big)\\
  &=4B_{-2,3}B_{2,-3}
  -\Big(M-\frac32\Big)\Big(M-\frac12\Big),
\end{align}
\end{subequations}
where the $B$, $N_1$ and $M$ operators should be equipped with the superscript ${}^{\sch,-\frac12}$.

In the coordinates $w,s$ we need to make the replacements
\begin{subequations}
\begin{align}\label{facb-2-}
  t&\quad\to\quad s^2,\\
    y_0^2&\quad\to\quad 2w^2s^2
\label{facb-1-}.\end{align}
\end{subequations}

\subsection{Factorizations of the Hermite operator}
\label{symcom5a}
The factorizations
of $\cL_1$ described in Subsect.
\ref{Factorizations of the heat operator2}
yield the following factorizations of the Hermite operator:
\begin{eqnarray*}
\cS_\lambda&=&\big(\p_w-2w\big)\p_w-2\lambda-1\\
&=&\p_w\big(\p_w-2w\big)-2\lambda+1,\\[0.5ex]
w^2\cS_\lambda&=&
\Big(w\p_w+\lambda-\frac32\Big)\Big(w\p_w-\lambda+\frac12-2w^2\Big)+\Big(\lambda-\frac32\Big)\Big(\lambda-\frac12\Big)\\
&=&
\Big(w\p_w-\lambda-\frac32-2w^2\Big)
\Big(w\p_w+\lambda+\frac12\Big)
+\Big(\lambda+\frac32\Big)\Big(\lambda+\frac12\Big).
\end{eqnarray*}

\subsection{Standard solutions}

The Hermite equation has only one singular point, $\infty$. One can define two kinds of solutions with a simple asymptotics at $\infty$. They can be derived from the expressions of Subsect. \ref{Standard solutions11},
using (\ref{ada1}) and (\ref{siu2})

\begin{align*}
\text{  $\sim w^{-{a}}$ for $w\to+\infty$:}&&\quad
S_\lambda (w)&:=w^{-\lambda -\frac12}\tilde F_{\lambda,\frac12 }(-w^{-2})\\
&&=&
w^{-{a}}F\Big(\frac{{a}}{2},\frac{{a}+1}{2};-;-w^{-2}\Big),\\
\text{${\sim}({-}\i w)^{{a}-1}\e^{w^2}$ for $w\to{+}\i\infty$:}&&\quad
\e^{w^2}S_{-\lambda }(-\i w)&=
(-\i w)^{\lambda -\frac12}\e^{w^2}\tilde F_{-\lambda,\frac12}(w^{-2})\\
&&&
\hspace{-14ex}=(-\i w)^{{a}-1}\e^{w^2}F\Big(\frac{1-a}{2},\frac{2-a}{2};-;-w^{-2}\Big).
\end{align*}

\subsection{Recurrence relations}

Each of the following recurrence relations corresponds to a root of $\sch(1)$:
\begin{eqnarray*}
\p_w S_\lambda(w)&=&-\Big(\frac12+\lambda \Big) S_{\lambda +1}(w),\\
(\p_w -2w)S_\lambda (w)&=&-2S_{\lambda -1}(w),\\[0.6ex]
\Big(w\p_w+\frac12+\lambda \Big)S_\lambda (w)&=&\frac12\Big(\frac12+\lambda \Big)
\Big(\frac32+\lambda\Big) S_{\lambda +2}(w),\\
\Big(w\p_w+\frac12-\lambda -2w^2\Big) S_\lambda (w)&=&-2
S_{\lambda -2}(w).
\end{eqnarray*}

The first pair corresponds correspond to the celebrated annihilation and creation operators in the theory of quantum harmonic oscillator. The second pair are the double annihilation and creation operators.

\subsection{Wave packets for the heat equation in 1 dimensions}

Consider the space $\rr\oplus\rr$ and the heat equation given by the operator $\cL_1=\partial_{y}^2+2\partial_t$.
Recall that
\beq M^{\sch,-\frac12}=y\p_y+2t\p_t+\frac12.\eeq

Set
\bes\begin{align}
  G_\lambda^a(y,t)&:=\int_{\gamma^a}
  t^{-\frac12}\exp\Big(\frac{(y-\tau^{-1})^2}{2t}\Big)
  \tau^{-\frac32+\lambda}\d\tau,\\
  G_\lambda^b(y,t)&:=
  \int_{\gamma^b}\e^{-\sqrt{2}y\tau-t\tau^2}\tau^{-\frac12+\lambda}\d\tau.
\end{align}\ees

\bep We have
\bes\begin{align}
  \cL_1 G_\lambda^a&=0,&\quad \cL_1 G_\lambda^b&=0;\\
  M^{\sch,-\frac12} G_\lambda^a&=-\lambda G_\lambda^a,&\quad M^{\sch,-\frac12} G_\lambda^b &=-\lambda G_\lambda^b.\label{ko2.}
\end{align}\ees
\eep

\proof
Set
\bes\begin{align} g^a(y,t)&:= t^{-\frac12}\exp\frac{(y-1)^2}{2t},\\
g^b(y,t)&:=\e^{-\sqrt{2} y-t}.\end{align}\ees
We have
\bes\begin{align}
  G_\lambda^a&=\int_{\gamma^a}\tau^{-1+\frac12+\lambda}g^a(\tau y,\tau^2 t) \d\tau,\\
    G_\lambda^b&=\int_{\gamma^b}\tau^{-1+\frac12+\lambda}g^b(\tau y,\tau^2 t) \d\tau.\label{kaka}\end{align}\ees
Clearly, $g^a$ and $g^b$ solve the heat equation.
 By (\ref{kaka}), $G_{\lambda}^a$, resp.
 $G_{\lambda}^b$ are wave packets made out of rotated $g^a$, resp. $g^b$. Therefore, they also solve the heat equation.

If the contours satisfy the requirements of Prop. \ref{analo}, then  (\ref{kaka})  implies (\ref{ko2.}).
 \qed

Let us express these wave packets in the coordinates $w,s$:
\bes\begin{align}
  G_\lambda^a(w,s)&
  =\int s^{-1}\exp\Big(\Big(w-\frac{1}{\sqrt2\tau s}\Big)^2\Big)\tau^{-2+\frac12+\lambda}\d\tau,\label{poe1}\\
  G_\lambda^b(w,s)
  &=\int\e^{-2sw\tau-s^2\tau^2}\tau^{-1+\frac12+\lambda}\d\tau.\label{poe2}
  \end{align}\ees
In (\ref{poe1}) 
we set $\sigma:=w-\frac{1}{\sqrt2\tau s}$,
so that $\tau=\frac{1}{(w-\sigma)\sqrt2 s}$.
In (\ref{poe2})
we set $\sigma:= s\tau$, so that $\tau=\frac{\sigma}{s}$. We obtain
\bes\begin{align}
  G_\lambda^a(w,s)&
  =(\sqrt2)^{\frac12-\lambda} s^{-\frac12-\lambda}
  F_\lambda^a(w),\\
  G_\lambda^b(w,s) &=s^{-\frac12-\lambda}F_\lambda^b(w),
  \end{align}\ees
where
\bes\begin{align}
F_\lambda^a(w)&:=
\int_{\gamma^a}\e^{\sigma^2}(w-\sigma)^{-\frac12-\lambda}\d\sigma,
\label{poee1}\\
F_\lambda^b(w)&:=\int_{\gamma^b}\e^{-2\sigma w-\sigma^2}\sigma^{-1+\frac12+\lambda}\d\sigma.\label{poee2}
  \end{align}\ees
The above analysis  shows that for appropriate contours
(\ref{poee1}) and (\ref{poee2}) are solutions of the Hermite equation.

\subsection{Integral representations}

Below we directly describe the two kinds of integral representations of solutions,  without passing through additional variables.

\bet\ben\item[a)]
Let $[0,1]\ni \tau\overset{\gamma}\mapsto t(\tau)$ satisfy
$\e^{t^2}(t-w)^{-{a}-1}\Big|_{t(0)}^{t(1)}=0.$
Then
\beq\cS({a};w,\p_w)\int_\gamma\e^{t^2}(t-w)^{-{a}}\d t=0.\label{dad10}\eeq
\item[b)]
Let $[0,1]\ni \tau\mapsto t(\tau)$ satisfy
$\e^{-t^2-2wt}t^{{a}}\Big|_{t(0)}^{t(1)}=0.$
Then\beq
\cS({a};w,\p_w)
\int_\gamma
\e^{-t^2-2wt}t^{a-1}\d t=0.\label{dad11}\eeq
\een\label{dad12}\eet

\proof
We check that for any contour $\gamma$
\begin{eqnarray*}
\text{lhs of (\ref{dad10})} &=&-{a} \int_\gamma\d t\,\p_t
\e^{t^2}(t-w)^{-{a}-1},\\
\text{lhs of (\ref{dad11})}&=&-2\int_\gamma\d t\,
\p_t\e^{-t^2-2wt}t^{{a}}.
\end{eqnarray*}

We can also deduce the second  representation from the first
by the discrete symmetry (\ref{siu2}).
\qed

\subsection{Integral representations of standard solutions}

In type a) representations  the integrand has a singular point at $0$ and goes to zero as $t\to \pm \infty$. We can thus use  contours with such endpoints. We will see that they give all standard solutions.

In type b) representations  the integrand has a singular point at $w$ and goes
to zero as $t\to \pm \i\infty$. Using  contours with such
endpoints, we will also obtain all standard solutions.

\medskip

\noindent
\begin{tabular}{lll}&a)&b)\\
  $\sim w^{-{a}}$ for $w\to+\infty$:\qquad&$[0,\infty[$,&$]-\i\infty,w^-,-\i\infty[$;\\
     $\sim(-\i w)^{{a}-1}\e^{w^2}$ for $w\to+\i\infty$:\qquad&
$]-\infty,0^+,-\infty[,$&$[w,\i\infty[$.\\
\end{tabular}

\medskip

It is convenient to introduce
 alternatively normalized solutions:
\begin{eqnarray*}
  S_\lambda^\I (w)&:=&2^{-\lambda-\frac12}\Gamma\Big(\lambda+\frac12\Big)S_\lambda(w).
  \\
\end{eqnarray*}
Here are integral representations of type a):
\begin{align}\text{all $\lambda$:}&&\\
\underset{]-\i\infty,w^-,\i\infty[}{-\i\int}\e^{t^2}
(w-t)^{-\lambda -\frac12}
\d
  t 
&=\sqrt\pi S_\lambda (w),&\quad w\not\in]-\infty,0];\notag\\
\Re\lambda<\frac12:&&\\
\underset{[w,\i\infty[}{-\i\int}\e^{t^2}(-\i (t-w))^{-\lambda -\frac12}\d t
&= \e^{w^2} S_{-\lambda }^\I(-\i w),&\quad w\not\in[0,\infty[.\notag
\end{align}
And here are  integral representations of type b):
\begin{align} -\frac12<\Re\lambda:&&\\
\int\limits_0^\infty\e^{-t^2-2tw}t^{\lambda -\frac12}\d t
&= S_\lambda ^\I(w),&\quad w\not\in]-\infty,0];\notag\\
\text{all $\lambda$:}&&\\
\int\limits_{]-\infty,0^+,\infty[}\e^{-t^2-2tw}(\i t)^{\lambda -\frac12}\d t&=
\sqrt\pi\e^{w^2} S_{-\lambda }(-\i w),&\quad w\not\in[0,\infty[.\notag
\end{align}


\section{The Helmholtz equation in 2 dimensions and the ${}_0\cF_1$  equation}
\label{s11}
\init

The goal of this section is to derive the ${}_0\cF_1$ equation together with its symmetries from the Helmoltz equation in $2$ dimensions.
The symmetries of these equations, together with its derivation, are the simplest and the best known. In particular, we do not need to consider {\em generalized symmetries}.

Here are the main steps from the derivation:

\ben
\item\label{it1he}
  We start from the Helmholtz operator  $\Delta_2-1$.  The  Lie algebra  $\aso(2)$ and group $\ASO(2)$
acting  in $2$ dimensions, are the obvious symmetries of this operator.
\item\label{it4he} We choose coordinates $w,u$, so that the Cartan element is expressed in terms of  $u$. We compute
  $\Delta_2-1$ and the representations
  of $\aso(2)$ and  $\ASO(2)$ in the new coordinates.
\item\label{it5he} We make an ansatz  diagonalizing the Cartan element,
  whose eigenvalue $\alpha$
  becomes a parameter. 
The only variable left is $w$. 
The Helmholtz operator  $\Delta_2-1$ becomes the  ${}_0\cF_1$  operator.
  The  symmetries of
  $\Delta_2-1$ yield transmutation relations and discrete symmetries of the
 ${}_0\cF_1$  operator.
\een

The remaining part of this section is to a large extent parallel to their analogs in Sects
\ref{s6}, \ref{s7}, \ref{s9} and \ref{s10}. Essentially all subsections have their counterparts there.
The  only exception is Subsect.
\ref{Equivalence with a subclass of the confluent equation}
on the equivalence of the ${}_0\cF_1$ equation with a subclass of the ${}_1\cF_1$ equation, and its many-dimensional unravelling.
This equivalence is obtained by a quadratic transformation, which
is quite different from the quadratic transformations for the Gegenbauer and Hermite equation considered in
Subsects \ref{Quadratic transformation}, resp.
 \ref{Quadratic transformation-}.


\subsection{$\aso(2)$}

We consider $\rr^2$ with split coordinates $x_-,x_+$ and the scalar product
\beq
\langle x|x\rangle=2x_-x_+.
\eeq
\noindent{\bf Lie algebra $\aso(\cc^2)$.}
Cartan operator
\beq N=-x_{-}\partial_{x_{-}}+x_{+}\partial_{x_+}. \label{helm1}\eeq
Root operators
\bes\begin{align}
B_- &=\p_{x_-},\\
B_+ &=\p_{x_+}.
\end{align}\ees

\noindent{\bf Weyl symmetry} 
\beq
\tau f(x_-,x_+)=f(x_+,x_-).
\eeq

\noindent{\bf Helmholtz operator}
\beq\Delta_2-1=
2\p_{x_-}\p_{x_+}-1.\label{helm2}\eeq

\subsection{Variables $w,u$}
\label{Variables $w,u$}

We introduce the coordinates
\begin{align}
w =\frac{x_-x_+}{2}\,,\qquad
u =x_+.
\label{coor3}\end{align}

\noindent{\bf Lie algebra $\aso(2)$.}
Cartan operator
\[ N=u\p_u.\]
Root operators
\begin{align*}
B_+& =\frac{u}{2}\partial_w,\\
B_-& =\frac{1}{u}(w\partial_w+u\partial_u).
\end{align*}

\noindent{\bf Weyl symmetry} 
\[
\tau
f(w,u)=f\Big(w,\frac{w}{u}\Big).
\]

\noindent{\bf Helmholtz operator}
\beq
\Delta_2-1=w\partial_w^2+(1+u\partial_u)\partial_w-1.
\eeq

\subsection{The ${}_0\cF_1$ operator}

Let us make the ansatz
\beq f(w,u)=u^\alpha F(w).\eeq
Clearly,
\begin{align} Nf&=\alpha f,\\
u^{-\alpha}  (\Delta_2-1)f&=\cF_\alpha(w,\partial_w)F,\end{align}
where we have introduced the {\em ${}_0\cF_1$ operator}
\beq
\cF_\alpha(w,\partial_w):=w\partial_w^2+(1+\alpha)\partial_w-1.
\eeq
Instead of the {\em Lie-algebraic parameter} $\alpha$ one could also use the {\em classical parameter}~$c$
\beq \alpha:=c-1,\ \ \ c=\alpha+1,\label{newnot3}\eeq
so that the ${}_0\cF_1$ operator becomes
\beq \cF(c;w,\p_w):=w\p_w^2+c\p_w-1.\eeq

\subsection{Equivalence with a subclass of the confluent equation}
\label{Equivalence with a subclass of the confluent equation}

The ${}_0\cF_1$ equation is equivalent to a subclass of the ${}_1\cF_1$ equation by a quadratic transformation. This quadratic transformation is however quite different from  transformations described in Subsect.
\ref{Dimensional reduction}, and then applied to derive the Gegenbauer quation 
and the Hermite equation. In this subsection we derive this equivalence starting from the heat equation in 2 dimensions.

First let us recall some elements of our derivation of the  ${}_1\cF_1$ operator.
As described in Sect. \ref{s9}, it was obtained from the heat operator
(\ref{heat0})
  together with  Cartan operators (\ref{heat1}), (\ref{heat2}):
\bes  \begin{align}
\frac{t}{2}\cL_2&=\frac{t}{2}\big(2\partial_t+2\partial_{y_{-1}}\partial_{y_1}\big),
\label{hhet} \\
M&=y_{-1}\p_{y_{-1}}+y_{1}\p_{y_{1}}+2t\p_t+1,\\
    N_1&=-y_{-1}\p_{y_{-1}}+y_{1}\p_{y_{1}}.
  \end{align}\ees
  (We  set $\eta=-1$ and dropped the superscript ${}^{\sch,-1}$).
Recall that substituting the coordinates (\ref{cor})
\beq    w =\frac{y_{-1}y_1}{t}\;,\ \ \ \ 
u_1 =\frac{y_{1}}{\sqrt{t}} \;,\ \ \ \
s =\sqrt{t}\eeq
we obtain
\bes  \begin{align}
    \frac{t}{2}\cL_2&=    w\p_w^2+(u\p_u+1-w)\p_w+\frac12(-u\p_u+s\p_s),
    \label{pocor}\\
    M&=s\p_s+1,\\
N_1&=u_1\p_{u_1}.
  \end{align}\ees
  After we set $M=-\theta$, $N_1=\alpha$, (\ref{pocor}) becomes $\cF_{\theta,\alpha}(w,\p_w)$.

  Consider now
  \begin{align}\notag
   & \frac{2t^2}{y_{-1}y_{1}}\e^{-\frac{y_{-1}y_1}{2t}}
    \cL_2
    \e^{\frac{y_{-1}y_1}{2t}}\\
    =&\frac{2t}{y_{-1}y_1}\big(   y_{-1}\p_{y_{-1}}+y_{1}\p_{y_{1}}+2t\p_t+1\big)+\frac{4t^2}{y_{-1}y_1}\p_{y_{-1}}\p_{y_1}-1\notag\\
    =&   \frac{2t}{y_{-1}y_1}M+2\p_{x_-}\p_{x_+}-1,\label{helmo}\\
    & \e^{-\frac{y_{-1}y_1}{2t}}
    N_1
    \e^{\frac{y_{-1}y_1}{2t}}=N_1\notag
    =-2x_-\p_{x_-}+2x_+\p_{x_+},
  \end{align}
  where we
  introduced new variables
  \beq
  x_-=\frac{y_{-1}^2}{2\sqrt2 t},\quad x_+=\frac{y_{1}^2}{2\sqrt2 t}.
  \eeq
  Therefore, on the subspace $M=0$ we have
  \begin{align}\notag
    \frac{2t^2}{y_{-1}y_{1}}\e^{-\frac{y_{-1}y_1}{2t}}
    \cL_2
    \e^{\frac{y_{-1}y_1}{2t}}
    =&\Delta_2-1,\\
     \e^{-\frac{y_{-1}y_1}{2t}}
    N_1
    \e^{\frac{y_{-1}y_1}{2t}}&=2N,
  \end{align}
  where $\Delta_2-1$ is the Helmholtz operator
(\ref{helm2})
  and $N$ the Cartan operator (\ref{helm1}).
  Remember, that in Subsect. \ref{Variables $w,u$}
  we express these operators
in the coordinates (\ref{coor3}). To avoid a clash of symbols, we rename $w$ from (\ref{coor3}) into $v$:
  \beq
  v=\frac{y_-y_+}{2},\quad
  u=y_+.\eeq
  Recall that in the $v,u$ coordinates we have
\bes  \begin{align}
\Delta_2-1&=v\partial_v^2+(1+u\partial_u)\partial_v-1,\label{wqe}\\
N&=u\p_u,
  \end{align}\ees
  so that (\ref{wqe}) on $N=\alpha$ becomes $\cF_\alpha(v,\p_v)$.

  Now    we can compare the coordinates $w,u_1$ and $v,u$
   \beq
  v=\frac{y_{-1}^2y_1^2}{16 t^2}=\Big(\frac{w}{4}\Big)^2,\quad
  u=\frac{y_1^2}{2\sqrt2t}=\frac{u_1^2}{2\sqrt2}.\eeq
This leads to the so-called {\em Kummer's 2nd transformation}, which reduces
the  ${}_0\cF_1$ equation to a special class of 
the confluent equation by a quadratic transformation: 
\beq
\cF_\alpha (v,\p_v)
=\frac{4}{w}\e^{-w/2}\cF_{0,2\alpha }(w,\p_w)\e^{w/2},\label{gas1}\eeq
or, in classical parameters
\begin{eqnarray}
\cF(c;v,\p_v)
=\frac{4}{w}\e^{-w/2}\cF\Big(c-\12;2c-1;w,\p_w\Big)\e^{w/2},
\label{gas}\end{eqnarray}
where $w=\pm 4\sqrt{v}$, $v=\big(\frac{w}{4}\big)^2$.

\subsection{Transmutation relations and symmetries}
\label{symcom2}

The following symmetries of the Helmholtz operator are obvious:

\begin{subequations}
  \begin{align}
  B(\Delta_2-1)&=(\Delta_2-1)B;\quad B\in \aso(2);\label{lorr3+}\\
  \alpha (\Delta_2-1)&=(\Delta_2-1)\alpha;\quad \alpha\in\ASO(2).\label{lorr3}
  \end{align}
\end{subequations}

Applying (\ref{lorr3+}) to the roots of $\aso(2)$ we obtain the trasmutation relations
\[\begin{array}{rrll}
\p_w&\cF_\alpha &
=\ \cF_{\alpha +1}&\p_w,\\[3ex]
(w\p_w+\alpha )&\cF_\alpha
&=\ \ \cF_{\alpha -1}&(w\p_w+\alpha ).\end{array}\]

Applying (\ref{lorr3}) to the Weyl symmetry of $\aso(2)$ we obtain the symmetry
\[\begin{array}{lcr}w^{-\alpha }\ \cF_{-\alpha }\ w^\alpha  &=&\cF_\alpha .\end{array}\]

\subsection{Factorizations}

The factorizations
\bes\begin{align}
\Delta_2-1&=2B_-B_+-1\\
&=2B_+B_--1,
\end{align}\ees
are completely obvious. They yield the factorizations of the ${}_0\cF_1$
operator:
\begin{eqnarray*}
\cF_\alpha&=&
\big(w\partial_w+\alpha+1\big)\partial_w-1\\
&=&\p_w
\big(w\partial_w+\alpha\big)-1.
\end{eqnarray*}

\subsection{The ${}_0F_1$ function}

The ${}_0\cF_1$ equation
 has a regular singular point at  $0$.
Its indices at $0$ are equal to $0$, $\alpha=1-c$.

If $c\neq0,-1,-2,\dots$, then the only solution of the ${}_0F_1$ equation
$\sim1$ at 0 
is called
the {\em ${}_0F_1$  function}.
It is 
\beq F(c;w):=\sum_{j=0}^\infty
\frac{1}{
(c)_j}\frac{w^j}{j!}.\eeq
It is defined for $c\neq0,-1,-2,\dots$.
Sometimes it is more convenient to consider
the function
\beq {\bf F}  (c;w):=\frac{F(c;w)}{\Gamma(c)}=
\sum_{j=0}^\infty
\frac{1}{
\Gamma(c+j)}\frac{w^j}{j!}\eeq
defined for all $c$.

Using (\ref{gas}), we can express the ${}_0F_1$ function in terms of the confluent function
\bes\begin{eqnarray}F(c;w)&=&
\e^{-2\sqrt{w}}F\Big(\frac{2c-1}{2};
2c-1;4\sqrt{w}\Big)\\&=&
\e^{2\sqrt{w}}F\Big(\frac{2c-1}{2};
2c-1;-4\sqrt{w}\Big).
\end{eqnarray}\ees

We will usually prefer to use the Lie-algebraic parameters:
\bes\begin{eqnarray}
F_\alpha (w)&:=&F(\alpha+1;w),\\
 {\bf F}  _\alpha (w)&:=& {\bf F} (\alpha+1;w)
.
\end{eqnarray}\ees

\subsection{Standard solutions}

We have two standard solutions corresponding to  two  indices
of the regular singular point $w=0$.
Besides, using Tricomi's function described in Subsect.  \ref{Standard solutions11}, we have an additional solution with a special behavior at $\infty$:
\begin{align*}
  \text{$\sim1$ at 0:}&&\quad
  F_\alpha (w)&=
\e^{-2\sqrt{w}}F_{0,2\alpha }\big(4\sqrt{w}\big)\notag\\&&&=
\e^{2\sqrt{w}}F_{0,2\alpha }\big(-4\sqrt{w}\big);\\[1ex]
\text{$\sim w^{-\alpha }$ at 0:}&&\quad
w^{-\alpha }F_{-\alpha} (w)&=
w^{-\alpha }\e^{-2\sqrt{w}}F_{0,-2\alpha }\big(4\sqrt{w}\big)\notag\\&&&=
w^{-\alpha }\e^{2\sqrt{w}}F_{0,-2\alpha }\big(-4\sqrt{w}\big);
\end{align*}
\begin{align*}
\text{
$\sim\e^{- 2\sqrt{w}} w^{-\frac{ \alpha }2-\frac14}$,
  $w\to+\infty$:}&&\quad
\tilde F_\alpha (w)&:=\e^{-2\sqrt w} w^{-\frac{\alpha}{2} -\frac14}
\tilde F_{0,2\alpha }\Big(-\frac{1}{4\sqrt w}\Big)\notag\\&&
&=\e^{-2\sqrt w} w^{-\frac{\alpha}{2} -\frac14}
\tilde F_{0,-2\alpha }\Big(-\frac{1}{4\sqrt w}\Big).
\end{align*}

Note that the third standard solution is a new function closely related to the MacDonald function. It satisfies
 the identity
\beq\tilde F_\alpha (w)=w^{-\alpha }\tilde F_{-\alpha }(w).\eeq
Its  asymptotics
 \beq
\tilde F_\alpha (w)\sim{\rm exp}(- 2w^\12) w^{-\frac{\alpha }2-\frac14}
\label{saddle1}\eeq
is valid
 in the sector
 $|\arg w|<\pi/2-\epsilon$ for
 $|w|\to\infty$.

\subsection{Recurrence relations}
The following recurrence relations follow from the transmutation relations
 \begin{eqnarray*}
 \p_w {\bf F}  _\alpha (w)&=& {\bf F}  _{\alpha +1}(w),
 \\[3mm]
 \left(w\p_w+\alpha\right) {\bf F}  _\alpha (w)&=& {\bf F}  _{\alpha -1}(w).
\end{eqnarray*}

\subsection{Wave packets}

Obviously, for any $\tau$
the function $\exp\Big(\frac{x_-}{\sqrt 2\tau}+\frac{\tau x_+}{\sqrt{2}}\Big)$ solves the Helmholtz equation. Therefore, for appropriate contours $\gamma$,
\beq
f(x_-,x_+):=\frac{1}{2\pi\i}\int_\gamma 
\exp\Big(\frac{x_-}{\sqrt 2\tau}+\frac{\tau x_+}{\sqrt{2}}\Big)\tau^{-\alpha-1}\d\tau\eeq
solves
\begin{align}
  (\Delta_2-1)f&=0,\\
  Nf&=\alpha f.
\end{align}
Substituting the coordinates $w,u$ we obtain
\begin{align}\notag
  f(w,u)&=
 \int_\gamma \exp\Big(\frac{w}{\tau u\sqrt2}+\frac{\tau u}{\sqrt2}\Big)\tau^{-\alpha-1}\d\tau
  \\
  &=u^\alpha 2^{-\frac{\alpha}{2}}\int_\gamma \exp\Big(\frac{w}{s}+s\Big)s^{-\alpha-1}\d s,\end{align}
where we made the substitution $s=\frac{\tau u}{\sqrt2}$.
Therefore,
\beq
F(w)=\int_\gamma \exp\Big(\frac{w}{s}+s\Big)s^{-\alpha-1}\d s.
\eeq
solves the ${}_0F_1$ equation.

\subsection{Integral representations}

There are three kinds of integral representations of solutions to the
${}_0F_1$ equation. The first is suggested by the previous subsection. Representations of the first kind 
 will be called {\em
   Bessel-Schl\"afli type representations}.
 The next two are inherited from the confluent equation by  2nd Kummer's identity. We will call them
 {\em
   Poisson-type representations}.

 \bet\label{schl}
 \ben\item[i)] {\bf Bessel-Schl\"afli type representations.}
Suppose that $[0,1]\ni t\mapsto\gamma(t)$ satisfies
\[
\e^t\e^{\frac{w}{t}}t^{-c}\Big|_{\gamma(0)}^{\gamma(1)}=0.\]
Then
\beq\cF(c;w,\p_w)
\int_\gamma\e^t\e^{\frac{w}{t}}t^{-c}\d t=0.\label{dad5}\eeq
\item[ii)] {\bf Poisson type a)  representations.} Let the contour $\gamma$ satisfy
\[(t^2-w)^{-c+3/2}\e^{2t}\Big|_{\gamma(0)}^{\gamma(1)}=0.\]
Then
\beq
\cF(c;w,\p_w)\int_\gamma(t^2-w)^{-c+1/2}\e^{2t}\d t=0.\label{in2}\eeq
\item[iii)] {\bf Poisson type b) representations.}
 Let the contour $\gamma$ satisfy
\[(t^2-1)^{c-1/2}\e^{2t\sqrt w}\Big|_{\gamma(0)}^{\gamma(1)}=0.\]
Then
\beq
\cF(c;w,\p_w)\int_\gamma(t^2-1)^{c-3/2}\e^{2t\sqrt w}\d t=0.\label{in3}\eeq
\een
\eet

\proof 
We check that for any contour $\gamma$
\[\text{lhs of  (\ref{dad5})}=-\int_\gamma\d t
\p_t\e^t\e^{\frac{w}{t}}t^{-c}.\]
This proves i).

To prove both Poisson type representations we use the quadratic relation
(\ref{gas}). Using the type a) representation for solutions of ${}_1\cF_1$  (\ref{dad1}), for appropriate contours $\gamma$ and $\gamma'$, we see that
\begin{eqnarray*}
&&\e^{-2\sqrt w}\int_\gamma\e^s s^{-c+\frac12}(s-4\sqrt w)^{-c+\frac12}\d s\\
&=&2^{-2c+2}\int_{\gamma'}\e^{2t}(t^2-w)^{-c+\frac12}\d t
\end{eqnarray*}
is annihilated by  $\cF(c)$, where we set
$t=\frac{s}{2}-\sqrt w$.  This proves ii).

Similarly, by the type b) representation  for solutions of ${}_1\cF_1$  (\ref{dad}), 
\begin{eqnarray*}
&&\e^{-2\sqrt w}\int_\gamma\e^{\frac{4\sqrt w}{s}} s^{-2c+1}(1-s)^{c-\frac32}\d s\\
&=&-2^{-2c+2}\int_{\gamma'}\e^{2t\sqrt w}(1-t^2)^{c-\frac32}\d t
\end{eqnarray*}
is annihilated by  $\cF(c)$, where we set
$t=\frac2s-1$.  This proves  iii).
\qed


\subsection{Integral representations of standard solutions}

In Bessel-Schl\"afli type representations the integrand goes to zero as $t\to-\infty$ and $t\to0-0$ (the latter for $\Re w>0$). Therefore, contours ending at these points yield solutions. We will see that in this way we can obtain all 3 standard solutions.

We can also obtain all solutions using Poisson type representations
(which are actually special cases of representations for solutions of the confluent equation).

\medskip

\noindent
\begin{tabular}{llll}
  &Bessel-Schl\"afli&Poisson type a) &Poisson type b)\\[2ex]
$  \sim1\text{ at }0$:&$]-\infty,0^+,\infty[$&
$[-1,1]$&
    \\
    $\sim w^{-\alpha }\text{ at }0$:&$(0-0)^+$&&$[-\sqrt{w},\sqrt{w}]$\\
$\sim\e^{- 2\sqrt{w}} w^{-\frac{ \alpha }2-\frac14}
\text{  for }w\to+\infty$:&$]-\infty,0]$&
$]-\infty,-1]$&$]-\infty,-\sqrt{w}]$\\
  \end{tabular}

\bigskip

Here are Bessel-Schl\"afli type representations. They are valid for all values of $\alpha$ and $\Re w>0$:
\begin{align}
\frac{1}{2\pi \i}\int\limits_{]-\infty,0^+,-\infty[}
\e^t\e^{\frac{w}{t}}t^{-\alpha -1}\d t
&= {\bf F}  _\alpha (w),\\
\frac{1}{2\pi \i}\int\limits_{[(0-0)^+]}\e^t\e^{\frac{w}{t}}t^{-\alpha -1}\d t
&=w^{-\alpha } {\bf F}  _{-\alpha }(w),\\
\int_{-\infty}^{0}\e^t\e^{\frac{w}{t}}(-t)^{-\alpha -1}\d t
&=\pi^{\12}\tilde F_\alpha (w).\end{align}

Next we give Poisson type representations, valid for
$w\not\in]-\infty,0]$:
\begin{align}
\Re \alpha >-\frac12:&\\
\int_{-1}^1(1-t^2)^{\alpha -\frac12}\e^{
  2t\sqrt w}\d t&=\Gamma(\alpha +\frac12)\sqrt\pi  {\bf F}  _\alpha (w),\notag\\
\frac12>\Re\alpha:&\\
\int_{-\sqrt w}^{\sqrt w}(w-t^2)^{-\alpha -\frac12}\e^{
2t}\d t&=\Gamma\Big({-}\alpha +\frac12\Big)\sqrt\pi w^{-\alpha }
    {\bf F}  _{-\alpha }(w);\notag\end{align}
    \begin{align}
 \Re \alpha >-\frac12:&\\
\int_{-\infty}^{-1}(t^2-1)^{\alpha -\frac12}\e^{
2t\sqrt w}\d t&=
\frac12\Gamma\Big(\alpha +\frac12\Big) \tilde F_\alpha (w),\notag\\
 \Re \alpha <\frac12:&\\
\int_{-\infty}^{-\sqrt w}(t^2-w)^{-\alpha -\frac12}\e^{
2t}\d t&=\frac12\Gamma\Big(-\alpha +\frac12\Big) \tilde F_\alpha (w).\notag
\end{align}

\subsection{Connection formulas}

From integral representations we easily obtain connection formulas.
As the basis we can use the solutions with a simple behavior at zero:
\begin{eqnarray*}
\tilde F_\alpha (w)
&=&\frac{\sqrt\pi}{\sin\pi (-\alpha )} {\bf F}  _\alpha (w)
+\frac{\sqrt \pi}{\sin\pi \alpha }
w^{-\alpha } {\bf F}  _{-\alpha }(w).
\end{eqnarray*}
Alternatively, we can use the basis conisting of the $\tilde F$ function and its clockwise or anti-clockwise analytic
continuation around $0$:
\begin{eqnarray*}
{\bf F}_\alpha(w)&=&\frac{1}{ 2\sqrt\pi}\left(\e^{\pm\i\pi
 ( \alpha+\12)}
\tF_\alpha(w)+
\e^{\mp\i\pi (\alpha+\12)}\tF_\alpha(\e^{\mp\i 2\pi}w)\right),\\
w^{-\alpha}{\bf F}_{-\alpha}(w)
&=&\frac{1}{ 2\sqrt\pi}\left(\e^{\mp\i\pi (\alpha-\12)}\tF_\alpha(w)-
\e^{\mp\i\pi( \alpha-\12)}\tF_\alpha(\e^{\mp\i 2\pi}w)\right).
\end{eqnarray*}

\end{document}